\documentclass[11pt]{article}
\usepackage{times}
\usepackage{amsmath} 
\usepackage{amssymb,latexsym}
\numberwithin{equation}{section} 
\usepackage{color} 
\usepackage{cancel} 
\usepackage{mathptmx}
\catcode`\@=11

\catcode`\@=11
\renewcommand\subsubsection{\@startsection{subsubsection}{2}{\z@}%
                                     {-3.25ex\@plus -1ex \@minus -.2ex}%
                                     {-0.01 mm}
                                     {\normalfont\bfseries}}
\newtheorem{Thm}{Theorem}[section]
\newtheorem{Lem}[Thm]{Lemma}
\newtheorem{Cor}[Thm]{Corollary}
\newtheorem{Prop}[Thm]{Proposition}
\newtheorem{Conj}[Thm]{Conjecture}
\newtheorem{example}[Thm]{Example}

\newtheorem{Def}[Thm]{Definition}

\definecolor{LawnGreen}{rgb}{0.49,0.99,0.00} 
\definecolor{ForestGreen}{rgb}{0.13,0.55,0.13}
\definecolor{DarkSeaGreen}{rgb}{0.56,0.74,0.56}
\definecolor{MediumSeaGreen}{rgb}{0.24,0.70,0.44}

\newcommand{\A}{\mathbb{A}}

\newcommand{\Z}{\mathbb{Z}}
\newcommand{\fF}{\mathbb{F}}
\newcommand{\Q}{\mathbb{Q}}

\newcommand{\N}{\mathbb{N}}
\newcommand{\C}{\mathbb{C}}
\newcommand{\R}{\mathcal{R}}

\newcommand{\rk}{\operatorname{rk}}

\newcommand{\md}{\operatorname{mod}}

\newcommand{\Hom}{\operatorname{Hom}} 
\newcommand{\Ext}{\operatorname{Ext}}
\newcommand{\End}{\operatorname{End}}

\newcommand{\dm}{\operatorname{dim}}  
\newcommand{\GL}{\operatorname{GL}} 
\newcommand{\dimv}{\underline{\dim}}

\newcommand{\Lam}{\Lambda}

\newcommand{\bsm}{\begin{smallmatrix}}
\newcommand{\esm}{\end{smallmatrix}}

\newcommand{\bi}{\mathbf{i}} 
\newcommand{\bj}{\mathbf{j}} 
 
\newcommand{\mm}{{\bf m}}

\def\proof{\medskip\noindent {\it Proof --- \ }}
\def\cqfd{\hfill $\Box$ \bigskip}

\def\CC{{\mathcal C}}
\def\resp{{\em resp. }}

\def\<{\langle\,}
\def\>{\,\rangle}

\def\eg{{\em e.g. }}

\def\T{{\mathcal T}}
\def\F{{\mathcal F}}
\def\P{{\mathbb P}}
\def\AA{{\mathcal A}}
\def\B{{\mathbf B}}
\def\BB{{\mathcal B}}

\def\V{{\mathcal V}}

\def\g{\mathfrak g}
\def\t{\mathfrak t}
\def\nn{\mathfrak n}
\def\bb{\mathfrak b}
\def\<{\langle}
\def\>{\rangle}
\def\n{{\mathfrak n}}
\def\si{\sigma}

\def\ra{\rightarrow}

\def\1{\mathbf 1}
\def\Gr{{\rm Gr}}

\def\aa{{\mathbf a}}
\def\bbf{{\mathbf b}}
\def\cc{{\mathbf c}}
\def\ii{{\mathbf i}}
\def\jj{{\mathbf j}}

\def\soc{{\rm soc}}
\def\P{{\mathcal P}}

\def\a{\alpha} 

\def\SC{\mathcal{S}}
\def\L{\Lambda}
\def\mod{{\rm mod}\,}

\def\Ext{{\rm Ext}}
\def\add{{\rm add}\,}
\def\End{{\rm End}}
\def\Hom{{\rm Hom}}

\def\de{\delta}
\def\De{\Delta}
\def\Sl{\mathfrak{sl}}
\def\la{\lambda}
\def\AA{\mathcal{A}}

\def\im{{\rm Im}\,}

\def\O{{\mathcal O}}
\def\1{{\mathbf 1}}
\def\ie{{\em i.e.\,}}

\def\G{\Gamma}
\def\eps{\epsilon}
\def\b{\beta}
\def\ga{\gamma}
\def\vphi{\varphi}

\def\vpi{\varpi}
\def\tpsi{\widetilde{\psi}}
\def\tvphi{\widetilde{\vphi}}

\def\De{\Delta}
\def\Oint{O_{\mathrm{int}}}
\def\r{{\mathrm{r}}}

\def\tf{{\widetilde f}}
\def\te{{\widetilde e}}
\def\bvphi{\overline{\vphi}}
\def\tB{\widetilde{B}}
\def\hI{\widehat{I}}
\def\Fac{\mathrm{Fac}}
\def\k{\kappa}
\def\Mat{\mathrm{Mat}}

\begin{document}

\title{\bf Cluster structures on quantum coordinate rings}
\author{C. Gei\ss, B. Leclerc, J. Schr\"oer}

\date{}

\maketitle

\begin{abstract}
We show that the quantum coordinate ring of the unipotent subgroup $N(w)$ of a
symmetric Kac-Moody group $G$ associated with a Weyl group element $w$ 
has the structure of a quantum cluster algebra. 
This quantum cluster structure arises naturally from a subcategory
$\CC_w$ of the module category of the corresponding  
preprojective algebra.
An important ingredient of the proof is a system of quantum determinantal identities
which can be viewed as a $q$-analogue of a $T$-system.
In case $G$ is a simple algebraic group of type $A, D, E$, we deduce 
from these results that the quantum coordinate ring of an open cell of
a partial flag variety attached to $G$ also has a cluster
structure.  
\end{abstract}

\bigskip
\setcounter{tocdepth}{1}
\tableofcontents

\section{Introduction}

Let $\g$ be the Kac-Moody algebra associated with a symmetric 
Cartan matrix.
Motivated by the theory of integrable systems in statistical mechanics
and quantum field theory, Drinfeld and Jimbo have introduced its 
quantum enveloping algebra $U_q(\g)$.
Let $\n$ denote the nilpotent subalgebra arising from
a triangular decomposition of $\g$. In the case when $\g$ is finite-dimensional,
Ringel \cite{Ri} showed that the positive part $U_q(\n)$ of $U_q(\g)$ can be
realized as the (twisted) Hall algebra of the category of representations
over $\fF_{q^2}$ of a quiver $Q$, obtained by orienting the Dynkin diagram of $\g$.
This was a major inspiration for Lusztig's geometric realization of $U_q(\n)$
in terms of Grothendieck groups of categories of perverse sheaves over
varieties of representations of $Q$, which is also valid when $\g$ is infinite-dimensional \cite{Lu0}.  

The constructions of Ringel and Lusztig involve the choice of an orientation
of the Dynkin diagram. In an attempt to get rid of this choice, Lusztig
replaced the varieties of representations of $Q$ by the varieties of
nilpotent representations of its preprojective algebra $\L = \L(Q)$,
which depends only on the underlying unoriented graph. 
He showed that one can realize the enveloping algebra $U(\n)$ 
as an algebra of $\C$-valued constructible functions over these nilpotent varieties
\cite{Lu0,Lu4}.
The multiplication of $U(\n)$ is obtained as a convolution-type product 
similar to the product of the Ringel-Hall algebra, but using
Euler characteristics of complex varieties instead of
number of points of varieties over finite fields. 
Note that this realization of $U(\n)$ is only available when the Cartan
matrix is symmetric.

One of the motivations of this paper was to find a similar construction of
the quantized enveloping algebra $U_q(\n)$, as a kind of Ringel-Hall algebra
attached to a category of representations of $\L$.
Unfortunately, there seems to be no simple way of $q$-deforming Lusztig's
realization of $U(\n)$. In this paper we try to overcome this difficulty by switching
to the dual picture. 

The dual of $U(\n)$ as a Hopf algebra can be identified
with the algebra $\C[N]$ of regular functions on the pro-unipotent group $N$
attached to $\n$. Dualizing Lusztig's construction, one can obtain 
for each nilpotent representation $X$ of $\L$ a distinguished regular function 
$\vphi_X\in\C[N]$, and the product $\vphi_X \vphi_Y$ can be calculated
in terms of varieties of short exact sequences with end-terms $X$ and $Y$
\cite{GLSMult}.  
For each element $w$ of the Weyl group $W$ of $\g$, 
the group $N$ has a finite-dimensional
subgroup $N(w)$ of dimension equal to the length of $w$. 
In particular, when $\g$ is finite-dimensional, we have
$N = N(w_0)$, where $w_0$ is the longest element of $W$.
In \cite{GLS}, we have shown that the coordinate ring $\C[N(w)]$ 
is spanned by the functions $\vphi_X$ where $X$ goes over
the objects of a certain subcategory $\CC_w$ of $\md(\L)$.
This category was
introduced by Buan, Iyama, Reiten and Scott in \cite{BIRS},
and independently in \cite{GLSUni1} for adaptable $w$.
Moreover, we have proved that $\C[N(w)]$ has a 
cluster algebra structure in the sense of Fomin and Zelevinsky \cite{FZ2},
for which the cluster monomials are of the form $\vphi_T$
for rigid objects $T$ of $\CC_w$.

The algebra $\C[N]$ has a quantum deformation, which we denote by
$A_q(\n)$ (see below \S\,\ref{multAqn}), and it is well known that 
the algebras $U_q(\n)$ and $A_q(\n)$ are in fact isomorphic. 
By works of Lusztig and De Concini-Kac-Procesi, $A_q(\n)$ has a subalgebra $A_q(\n(w))$, 
which can be regarded
as a quantum deformation of $\C[N(w)]$.
On the other hand, Berenstein and Zelevinsky \cite{BZ} have introduced the  
concept of a quantum cluster algebra. They have conjectured that
the quantum coordinate rings of double Bruhat cells in semisimple algebraic groups should have
a quantum cluster algebra structure.

In this paper, we introduce for every $w$ an explicit quantum cluster algebra
$\AA_{\Q(q)}(\CC_w)$, defined in a natural way in terms of the
category $\CC_w$. In particular, for every reachable rigid object $T$ 
of $\CC_w$ there is a corresponding quantum cluster monomial $Y_T \in \AA_{\Q(q)}(\CC_w)$. 
Our main result is the following quantization of the above theorem of \cite{GLS}.
\begin{Thm}\label{ThIntro}
There is an algebra isomorphism $\k\colon \AA_{\Q(q)}(\CC_w) \stackrel{\sim}{\to} A_q(\n(w))$. 
\end{Thm}

Note that quantizations of coordinate rings and quantizations of cluster algebras
are defined in very different ways. For example $A_q(\n) \cong U_q(\n)$ is 
given by its Drinfeld-Jimbo presentation, obtained by $q$-deforming the Chevalley-Serre-type
presentation of $U(\n)$. In contrast, quantum cluster algebras are defined
as subalgebras of a skew field of rational functions in $q$-commuting variables,
generated by a usually infinite number of elements given by an inductive procedure.

As a matter of fact, there does not seem to be so many examples 
of ``concrete'' quantum cluster algebras in the literature.
Grabowski and Launois \cite{GL} have shown that the quantum coordinate rings of
the Grassmannians $\Gr(2,n) \ (n\ge 2)$, $\Gr(3,6)$, $\Gr(3,7)$, and $\Gr(3,8)$
have a quantum cluster algebra structure.
Lampe \cite{La,La2} has proved two particular instances of Theorem~\ref{ThIntro}, namely
when $\g$ has type $A_n$ or $A_1^{(1)}$ and $w=c^2$ is the square of a Coxeter element. 
Recently, the existence of a quantum cluster structure on every algebra $A_q(\n(w))$ was 
conjectured by Kimura \cite[Conj.1.1]{Ki}.

Now Theorem~\ref{ThIntro} provides a large class of such examples,
including all algebras $U_q(\n)$ for $\g$ of type $A,D,E$.
By taking $\g = \Sl_n$ and some special permutation $w_k\in S_n$,
one also obtains that the quantum coordinate ring $A_q(\Mat(k,n-k))$ 
of the space of $k\times(n-k)$-matrices
has a quantum cluster algebra structure for every $1\le k\le n$.
This may be regarded as a cluster structure
on the quantum coordinate ring of an open cell of $\Gr(k,n)$.
More generally, for any simply-laced simple algebraic group~$G$, and any
parabolic subgroup $P$ of $G$, we obtain a cluster structure
on the quantum coordinate ring of the unipotent radical $N_P$ of $P$,
which can be regarded as a cluster structure
on the quantum coordinate ring of an open cell of the partial flag
variety $G/P$.

Note also that, by taking $w$ equal to the square of a Coxeter element,
our result gives a Lie theoretic realization
of all quantum cluster algebras associated with an arbitrary acyclic quiver
(but with a particular choice of coefficients).

Our strategy for proving Theorem~\ref{ThIntro} can be summarized as
follows. 

Let $A_q(\g)$ be the quantum analogue of the coordinate ring
constructed by Kashiwara \cite{K}.
We first obtain a general quantum determinantal
identity in $A_q(\g)$ (Proposition~\ref{qmin-id1}). This is 
a $q$-analogue (and an extension to the Kac-Moody case) of a determinantal identity of
Fomin and Zelevinsky \cite{FZ}. We then transfer this identity
to $A_q(\n)$ (Proposition~\ref{unip_minor_identity}). 
Here a little care must be taken since the restriction
map from $A_q(\g)$ to $A_q(\n)$ is \emph{not} a ring homomorphism.
Appropriate specializations of this identity give rise, for every
$w\in W$, to a system $\Sigma_w$
of equations (Proposition~\ref{Tsystem}) 
allowing to calculate certain quantum minors depending on $w$ in a recursive
way. This system is a $q$-analogue of a $T$-system 
arising in various problems of mathematical physics and
combinatorics (see \cite{KNS}), and we believe it could be of independent
interest. It will turn out that all quantum minors involved in $\Sigma_w$ belong
to the subalgebra $A_q(\n(w))$, and that among them we find a set of
algebra generators.

On the other hand, we show that the generalized determinantal
identities of \cite[Theorem 13.1]{GLS}, which relate a distinguished 
subset of cluster variables, take exactly the same form as 
the quantum $T$-system $\Sigma_w$
when we lift them to the quantum cluster
algebra $\AA_{\Q(q)}(\CC_w)$.
Therefore, after establishing 
that the quantum tori consisting of the initial variables of the
two systems are isomorphic, we can construct an injective algebra homomorphism
$\k$ from $\AA_{\Q(q)}(\CC_w)$ to the skew field of fractions $F_q(\n(w))$
of $A_q(\n(w))$, and show that the image of $\k$ contains a set of
generators of $A_q(\n(w))$. Finally, using an argument of specialization
$q \to 1$ and our result of \cite{GLS}, 
we conclude that $\k$ is an isomorphism from $\AA_{\Q(q)}(\CC_w)$
to $A_q(\n(w))$.

Another motivation of this paper was the open orbit conjecture of 
\cite[\S18.3]{GLS}. This conjecture states that all functions $\vphi_T$
associated with a rigid $\L$-module $T$ belong to the dual canonical 
basis of $\C[N]$. 
It can be seen as a particular instance of the general principle of Fomin and 
Zelevinsky \cite{FZ2} according to which, in every cluster algebra coming from Lie theory,
cluster monomials should belong to the dual canonical basis.
Since the dual canonical basis of $\C[N]$ is obtained by specializing 
at $q=1$ the basis $\B^*$ of $A_q(\n)$ dual to Lusztig's canonical basis 
of $U_q(\n)$, it is natural to conjecture that, more precisely, 
every quantum cluster monomial $Y_T$ of $\AA_{\Q(q)}(\CC_w)$ is mapped by $\k$ to 
an element of $\B^*$.
In fact,
it is not too difficult to show that $\k(Y_T)$ always satisfies
one of the two characteristic properties of $\B^*$ (see below \S\ref{canonbasACw}).
But unfortunately, the second property remains elusive, although 
Lampe \cite{La,La2} has proved it for all cluster variables in the two special
cases mentioned above.  

Finally we note that it is well known that the algebras $A_q(\n(w))$ are
skew polynomial rings. Therefore, by Theorem~\ref{ThIntro}, 
all quantum cluster algebras of the form $\AA_{\Q(q)}(\CC_w)$ are also 
skew polynomial rings, which is far from obvious from their definition.
One may hope that, conversely, the existence of a cluster structure on
many familiar quantum coordinate rings will bring some new insights for studying
their ring-theoretic properties, a very active subject in recent years
(see \eg \cite{BG,GLL,MC,Y} and references therein).

\section{The quantum coordinate ring $A_q(\g)$}

\subsection{The quantum enveloping algebra $U_q(\g)$}
Let $\g$ be a symmetric Kac-Moody algebra 
with Cartan subalgebra $\t$.
We follow the notation of \cite[\S1]{K}. 
In particular, 
we denote by $I$ the indexing set of the simple roots $\a_i\ (i\in I)$
of $\g$, by $P\subset\t^*$ its weight lattice, by $h_i (i\in I)$ 
the elements of $P^*\subset\t$
such that 
$\< h_i, \a_j\>=a_{ij}$ are the entries of 
the generalized Cartan matrix of $\g$.
Since $\g$ is assumed to be symmetric, we also have a 
symmetric bilinear form $(\cdot,\cdot)$ on $\t^*$ such that
$(\a_i,\a_j) = a_{ij}$.

The Weyl group $W<\GL(\t^*)$ is the Coxeter group generated by the 
reflections $s_i$ for $i\in I$, where
\[
s_i(\ga):= \ga -\<h_i,\ga\>\a_i.
\]
The length of $w\in W$ is denoted by $l(w)$. 
We will also need the contragradient action of $W$ on $P^*$.

Let $U_q(\g)$ be the corresponding quantum enveloping algebra, 
a $\Q(q)$-algebra with generators $e_i, f_i \ (i\in I), q^h\ (h\in P^*)$.
We write $t_i=q^{h_i}$.
We denote by $U_q(\nn)$ (\resp $U_q(\nn_-)$) 
the subalgebra of $U_q(\g)$
generated by $e_i\ (i\in I)$ (\resp $f_i\ (i\in I))$.
For $i\in I$, let $U_q(\g_i)$ denote the subalgebra of $U_q(\g)$
generated by $e_i, f_i, q^h\ (h\in P^*)$.

Let $M$ be a (left) $U_q(\g)$-module. 
For $\la\in P$, 
let
$M_{\la} = \{m\in M \mid q^h m = q^{\<\la,\,h\>} m \mbox{ for every } h\in P^*\}$ 
be the corresponding weight space of $M$. 
We say that $M$ is \emph{integrable} if (i) $M = \oplus_{\la\in P}\ M_\la$,
(ii) for any $i$, $M$ is a direct sum of finite-dimensional $U_q(\g_i)$-modules,
and (iii) for any $m\in M$, there exists $l\ge 0$ such that 
$e_{i_1}\cdots e_{i_l} m = 0$ for any $i_1,\ldots, i_l\in I$. 
We denote by $\Oint(\g)$ the category of integrable $U_q(\g)$-modules.
This is a semisimple category, with simple objects the irreducible
highest weight modules $V(\la)$ with highest weight $\la\in P_+$, the 
monoid of dominant weights.

\subsection{Bimodules}
Let $\vphi$ and $*$ be the $\Q(q)$-linear anti-automorphisms of $U_q(\g)$
defined by
\begin{eqnarray}\label{vphi}
&&\vphi(e_i) = f_i,\quad \vphi(f_i) = e_i, \quad \vphi(q^h) = q^h, \\[2mm]
&&e_i^* = e_i,\quad f_i^* = f_i, \quad (q^h)^* = q^{-h}. \label{star}
\end{eqnarray}
A right $U_q(\g)$-module $N$ gives rise to a left $U_q(\g)$-module $N^\vphi$
by defining 
\begin{equation}
x\cdot n = n\cdot\vphi(x), \qquad (n\in N, \ x\in U_q(\g)). 
\end{equation}
We say that $N$ is an integrable right module if $N^\vphi$ is an integrable
left module.
In particular, for $\la\in P_+$, we have an irreducible integrable right
module $V^\r(\la)$ such that $V^\r(\la)^\vphi = V(\la)$.
Let $m_\la$ be a highest weight vector in $V(\la)$, \ie $e_i\, m_\la = 0$
for any $i\in I$. Then $m_\la$ can be regarded as a vector $n_\la\in V^\r(\la)$,
which satisfies $n_\la \, f_i = 0$ for any $i\in I$.
Equivalently, $V^\r(\la)$ is isomorphic to the graded dual of $V(\la)$, endowed
with the natural right action of $U_q(\g)$. It follows that we have a natural
pairing $\<\cdot,\cdot\>_\la\colon V^\r(\la) \times V(\la) \to \Q(q)$, which satisfies
$\<n_\la,m_\la\>_\la = 1$, and
\begin{equation}
\<n\,x,\, m\>_\la = \<n,\, x\,m\>_\la, \qquad 
(m\in V(\la),\ n\in V^\r(\la),\ x\in U_q(\g)). 
\end{equation}
We denote by $\Oint(\g^{\rm op})$ the category of integrable right 
$U_q(\g)$-modules.
It is also semisimple, with simple objects $V^\r(\la)\ (\la\in P_+)$.

The tensor product of $\Q(q)$-vector spaces $V^\r(\la)\otimes V(\la)$ 
has the natural structure of a $U_q(\g)$-bimodule, via
\begin{equation}\label{bimod_matrix_coeff}
x\cdot(n\otimes m)\cdot y = (n\cdot y)\otimes (x\cdot m),\qquad 
(x,y \in U_q(\g),\ m\in V(\la),\ n \in V^\r(\la)). 
\end{equation}
 
\subsection{Dual algebra}\label{Uqstar}
Let $U_q(\g)^*=\Hom_{\Q(q)}(U_q(\g),\Q(q))$. This is a $U_q(\g)$-bimodule, via
\begin{equation}\label{bimod_Uqstar}
(x\cdot\psi\cdot y)(z) = \psi(yzx),\qquad 
(x,y,z\in U_q(\g),\ \psi\in U_q(\g)^*). 
\end{equation}
On the other hand, $U_q(\g)$ is a Hopf algebra, with comultiplication 
$\De\colon U_q(\g) \to U_q(\g)\otimes U_q(\g)$
given by
\begin{equation}\label{comult}
\De(e_i) = e_i\otimes 1 + t_i \otimes e_i,\quad
\De(f_i) = f_i\otimes t_i^{-1} + 1 \otimes f_i,\quad 
\De(q^h) = q^h\otimes q^h, 
\end{equation}
with counit $\eps\colon U_q(\g)\ra\Q(q)$ given by
\begin{equation}
\eps(e_i)=\eps(f_i)=0,\quad\eps(q^h)=1,
\end{equation}
and with antipode $S$ given by
\begin{equation}\label{antipode}
S(e_i)=-t_i^{-1}e_i,\quad S(f_i) = -f_i\,t_i,\quad S(q^h) = q^{-h}. 
\end{equation}

Dualizing $\De$, we obtain a multiplication on $U_q(\g)^*$, defined by
\begin{equation}  \label{multA}
(\psi\,\theta)(x) = (\psi\otimes\theta)(\De( x)), \qquad 
(\psi, \theta \in U_q(\g)^*,\ x\in U_q(\g)). 
\end{equation}
Later on it will often be convenient to use Sweedler's notation 
$\De(x) = \sum x_{(1)}\otimes x_{(2)}$ for
the comultiplication. Using this notation, (\ref{multA}) reads 
$(\psi\,\theta)(x) =\sum \psi(x_{(1)})\theta(x_{(2)})$. 
Combining~(\ref{bimod_Uqstar}) and~(\ref{multA}) we obtain
\begin{equation} \label{bimod_alg}
x\cdot(\psi\,\theta)\cdot y =
\sum (x_{(1)}\cdot\psi\cdot y_{(1)})\,(x_{(2)}\cdot\theta\cdot y_{(2)}).
\end{equation}

\subsection{Peter-Weyl theorem}
Following Kashiwara \cite[\S7]{K}, we define $A_q(\g)$ as the subspace of 
$U_q(\g)^*$ consisting of the linear forms $\psi$ such that the left 
submodule $U_q(\g) \psi$ belongs to $\Oint(\g)$,
and the right submodule $\psi U_q(\g)$ belongs to $\Oint(\g^{\rm op})$.
It follows from the fact that the categories $\Oint(\g)$ and 
$\Oint(\g^{\rm op})$
are closed under tensor product that $A_q(\g)$ is a subring of $U_q(\g)^*$.

The next proposition of Kashiwara can be regarded as a $q$-analogue of the Peter-Weyl 
theorem for the Kac-Moody group $G$ attached to $\g$ (see \cite{KP}).
We include a proof for the convenience of the reader.

\begin{Prop}[{\cite[Proposition 7.2.2]{K}}] \label{PW-thm}
We have an isomorphism $\Phi$ of $U_q(\g)$-bimodules
\[
\bigoplus_{\la\in P_+} V^\r(\la)\otimes V(\la) 
\stackrel{\sim}{\longrightarrow} A_q(\g)
\]
given by
\[
\Phi(n\otimes m)(x) = \<n\,x,\,m\>_\la,\qquad 
(m \in V(\la),\ n \in V^\r(\la),\ x\in U_q(\g)).  
\]
\end{Prop}
\proof
It follows from (\ref{bimod_matrix_coeff}) and 
(\ref{bimod_Uqstar}) that $\Phi$ defines a homomorphism
of $U_q(\g)$-bimodules from $\oplus_{\la\in P_+} V^\r(\la)\otimes V(\la)$ to 
$U_q(\g)^*$.  Since $V(\la)$ and $V^\r(\la)$ are integrable for all 
$\la\in P_+$, we see that $\im\Phi \subseteq A_q(\g)$.

Let us show that $\Phi$ is surjective. 
Let $\psi\in A_q(\g)$. We want to show that $\psi\in\im\Phi$.
Since $V:=U_q(\g)\psi$ is integrable, it decomposes as a (finite) direct sum of
irreducible integrable modules. Thus, without loss of generality, we may
assume that $V$ is isomorphic to $V(\la)$ for some $\la\in P_+$.
We may also assume that $\psi$ is a weight vector of $V$ (otherwise
we can decompose it as a sum of weight vectors).
Since $V:=U_q(\g)\psi$ and $W:=\psi U_q(\g)$ are both integrable, we
see that there exist $k$ and $l$ such that 
$e_{i_1}\cdots e_{i_k}\cdot\psi = \psi\cdot f_{j_1}\cdots f_{j_l}=0$
for every $i_1, \ldots, i_k, j_1, \ldots, j_l \in I$. 
Hence $\psi(xe_{i_1}\cdots e_{i_k})=\psi(f_{j_1}\cdots f_{j_l}x) = 0$
for every $x\in U_q(\g)$.
It follows that the linear form $a\in V^*$ defined by $a(\vphi)=\vphi(1)$
takes nonzero values only on a finite number of weight spaces of $V$.
Hence $a$ is in the graded dual of $V$, which we can identify to 
$V^r(\la)$. Moreover, 
\[
 \Phi(a\otimes \psi)(x) = \<a,\,x\psi\>_\la = a(x\psi) = x\psi(1) = \psi(x)
\]
for every $x\in U_q(\g)$. Therefore $\psi = \Phi(a\otimes \psi)$ belongs
to $\im\Phi$.

Now, $\Phi$ is also injective. Indeed, if for $n\otimes m\in 
V^\r(\la)\otimes V(\la)$ 
we have $\Phi(n\otimes m) = 0$, then for every $x\in U_q(\g)$, $\<nx,\,m\>_\la=0$. 
If $n\not = 0$, since $\<\cdot,\,\cdot\>_\la$
is a pairing between $V^\r(\la)$ and $V(\la)$ and $n\,U_q(\g) = V^\r(\la)$,
we get that $m=0$ and $n\otimes m = 0$. Hence the restriction of 
$\Phi$ to $V^\r(\la)\otimes V(\la)$
is injective. 
Finally, since the bimodules $V^\r(\la)\otimes V(\la)$ are simple and
pairwise non-isomorphic, $\Phi$ is injective. \cqfd 

In view of Proposition~\ref{PW-thm}, we can think of $A_q(\g)$ as a $q$-analogue
of the coordinate ring of the Kac-Moody group attached to $\g$ in \cite{KP}. 
We therefore call $A_q(\g)$ the \emph{quantum coordinate ring}.
When $\g$ is a simple finite-dimensional Lie algebra, 
$A_q(\g)$ is the quantum coordinate ring $\O_q(G)$ of the simply-connected   
simple Lie group with Lie algebra $\g$ studied by many authors,
see \eg \cite[\S 9.1.1]{J}.

\subsection{Gradings} Let $Q\subset P$ be the root lattice. It follows from the defining relations
of $U_q(\g)$ that it is a $Q$-graded algebra:
\begin{equation}
U_q(\g)=\bigoplus_{\alpha\in Q} U_q(\g)_\alpha,
\end{equation}
where
\begin{equation}
U_q(\g)_\alpha=\{x\in U_q(\g)\mid q^h x q^{-h} =q^{\< h,\a\>} x\ \mbox{ for all } h\in P^*\}.
\end{equation}
By Proposition~\ref{PW-thm}, we have 
\begin{equation}
A_q(\g)=\bigoplus_{\ga,\de\in P} A_q(\g)_{\ga,\de},
\end{equation}
where
\begin{equation}
A_q(\g)_{\ga,\de}=\{\psi\in A_q(\g)\mid 
q^l\cdot\psi\cdot q^r=q^{\< r,\ga\>+\< l,\de\>}\psi\ \mbox{ for all } r,l\in P^*\}.
\end{equation}

\begin{Lem}\label{gradings}
\begin{itemize} 
\item[(a)] 
With the above decomposition  $A_q(\g)$ is a $P\times P$-graded algebra.
\item[(b)]
For $x\in U_q(\g)_\a$, $\psi\in A_q(\g)_{\ga,\de}$ and $y\in U_q(\g)_\b$, we have
$x\cdot\psi\cdot y\in A_q(\g)_{\ga-\b,\de+\a}$.
\item[(c)] For $x\in U_q(\g)_\a$, $\psi\in A_q(\g)_{\ga,\de}$, we have $\psi(x)\not = 0$
only if $\a = \ga - \de$.
\end{itemize}
\end{Lem}

\section{Determinantal identities for quantum minors}
\subsection{Quantum minors} \label{ssec:qmin-def}
For our convenience we reproduce \emph{mutatis mutandis} 
a part of~\cite[\S9.2]{BZ}.
Using the isomorphism $\Phi$ from Proposition~\ref{PW-thm} we define for
each $\la\in P_+$ the element 
\begin{equation} \label{DefPrincMin}
\De^\la:=\Phi(n_\la\otimes m_\la)\in A_q(\g)_{\la,\la}.
\end{equation}
This is a $q$-analogue of a (generalized) principal minor, in the sense
of \cite[\S1.4]{FZ}.
An easy calculation shows that
\begin{equation} \label{PrincMin_Prp} 
\De^\la(f\,q^h\,e)=\eps(f)\,q^{\< h,\la\>}\eps(e),\qquad
(f\in U_q(\nn_-),\ h\in P^*,\ e\in U_q(\n)).
\end{equation}

For $(u,v)\in W\times W$, we choose reduced expressions
$\bi=(i_{l(u)},\ldots,i_2,i_1)$ and $\bj=(j_{l(v)},\ldots,j_2,j_1)$ so that
$u = s_{i_{l(u)}}\cdots s_{i_2}s_{i_1}$ and 
$v = s_{j_{l(v)}}\cdots s_{j_2}s_{j_1}$.
Next, we introduce positive roots 
\begin{equation}\label{defbeta}
\b_k=s_{i_1}s_{i_2}\cdots s_{i_{k-1}}(\a_{i_k}),\quad
\ga_l=s_{j_1}s_{j_2}\cdots s_{j_{l-1}}(\a_{j_l}),\qquad
(1\leq k\leq l(u),\  1\leq l\leq l(v)).
\end{equation}
Finally, for $\la\in P_+$, we set 
\begin{equation}\label{eq-bc}
b_k = (\b_k,\la),\quad c_l = (\ga_l,\la),
\qquad
(1\leq k\leq l(u),\  1\leq l\leq l(v)),
\end{equation}
and we define the (generalized) \emph{quantum minor}
$\De_{u(\la),v(\la)}\in A_q(\g)$ by
\begin{equation}
\De_{u(\la),v(\la)}= \left(f_{j_{l(v)}}^{(c_{l(v)})}\cdots f_{j_1}^{(c_1)}\right)
\cdot\De^\la\cdot \left(e_{i_1}^{(b_1)}\cdots e_{i_{l(u)}}^{(b_{l(u)})}\right).
\end{equation}
Here, as usual, we denote by $e_i^{(k)}$ (\resp $f_i^{(k)}$) the $q$-divided powers
of the Chevalley generators.
Equivalently, for $x\in U_q(\g)$ we have
\begin{equation}\label{DeltaValue}
\De_{u(\la),v(\la)}(x)=
\De^\la\left(e_{i_1}^{(b_1)}\cdots e_{i_{l(u)}}^{(b_{l(u)})} x\,
f_{j_{l(v)}}^{(c_{l(v)})}\cdots f_{j_1}^{(c_1)}\right).
\end{equation}
It follows from the quantum Verma relations~\cite[Proposition~39.3.7]{Lu}
that $\De_{u(\la),v(\la)}$ depends only on the pair of weights $(u(\la),v(\la))$,
and not on the choice of $u$ and $v$, or of their reduced expressions $\bi$ and $\bj$.
Moreover it is immediate that $\De_{u(\la),v(\la)}\in A_q(\g)_{u(\la),v(\la)}$.

We have the following direct consequence of the definition 
of quantum minors.
\begin{Lem}\label{cons_def}
If  $l(s_iu)=l(u)+1$ and $l(s_jv)=l(v)+1$ then
\[
\De_{s_iu(\la),s_jv(\la)}=f_j^{(c)}\cdot\De_{u(\la),v(\la)}\cdot e_i^{(b)},
\] 
where $b:=( \a_i,u(\la))\geq 0$ (\resp $c:=(\a_j,v(\la))\geq 0$) is the
maximal natural number such that
$\De_{u(\la),v(\la)}\cdot e_i^{(b)}\neq 0$ 
(\resp $~f_j^{(c)}\cdot\De_{u(\la),v(\la)}\neq 0$).
\end{Lem}

It is convenient to identify $\De_{u(\la),v(\la)}$ with an
extremal  vector of weight $(u(\la),v(\la))$ in the
simple $U_q(\g)\otimes U_q(\g)$ highest weight module $V^r(\la)\otimes V(\la)$.
Thus, we have
\begin{align}
f_i\cdot\De_{\ga,\de} &=0 \text{ if }(\a_i,\de)\leq 0, &  
e_i\cdot\De_{\ga,\de} &=0 \text{ if }(\a_i,\de)\geq 0, \label{eq:leftz}\\
\De_{\ga,\de}\cdot e_i&=0 \text{ if }(\a_i,\ga)\leq 0, &  
\De_{\ga,\de}\cdot f_i&=0 \text{ if }(\a_i,\ga)\geq 0. \label{eq:rightz}
\end{align}
In particular, we have 
$f_{j_{l(v)}}\cdot\De_{u(\la),v(\la)}=0=\De_{u(\la),v(\la)}\cdot e_{j_{l(u)}}$.

One may also identify $\De_{u(\la),v(\la)}$ with a matrix coefficient in $V(\la)$.
To do that, let us first denote by
\begin{equation}
m_{v(\la)} :=  f_{j_{l(v)}}^{(c_{l(v)})}\cdots f_{j_1}^{(c_1)} m_\la
\end{equation}
the extremal weight vector of $V(\la)$ with weight $v(\la)$.
Next, let $(\cdot,\,\cdot)_\la$ be the nondegenerate bilinear form on $V(\la)$ defined
by
\begin{equation}
(xm_\la,\,ym_\la)_\la  := \<n_\la\vphi(x),\,ym_\la\>_\la,\qquad (x,y\in U_q(\g)). 
\end{equation}
Then, using (\ref{DefPrincMin}) and (\ref{DeltaValue}) we easily get
\begin{equation}\label{matrixCoeff}
\De_{u(\la),v(\la)}(x)= (m_{u(\la)},\,xm_{v(\la)})_\la,\qquad (x\in U_q(\g)). 
\end{equation}

\subsection{A family of identities for quantum minors} 
Let $(\vpi_i)_{i\in I}$ be the fundamental weights, \ie
we have $\<h_j,\vpi_i\>=(\a_j,\vpi_i)=\de_{ij}$. We note that the fundamental weights
are only determined up to a $W$-invariant element. Moreover, it
is useful to observe that 
\begin{equation}\label{sivpi}
s_i(\vpi_j)=\begin{cases} 
\vpi_j-\a_j &\text{ if } i=j,\\
\vpi_j &\text{ otherwise.}
\end{cases}
\end{equation}
We can now state the main result of this section. This is a $q$-analogue
of~\cite[Theorem~1.17]{FZ}.
\begin{Prop} \label{qmin-id1}
Suppose that for $u,v\in W$ and $i\in I$ we have $l(us_i)=l(u)+1$ and
$l(vs_i)=l(v)+1$. Then 
\begin{equation}\label{eq:2.10}
\De_{us_i(\vpi_i),vs_i(\vpi_i)}\,\De_{u(\vpi_i),v(\vpi_i)}=
q^{-1}\De_{us_i(\vpi_i),v(\vpi_i)}\,\De_{u(\vpi_i),vs_i(\vpi_i)}+
\prod_{j\neq i}\De_{u(\vpi_j),v(\vpi_j)}^{-a_{ji}}
\end{equation}
holds in $A_q(\g)$. 
\end{Prop}
The proof of this proposition will be given after some preparation in 
Section~\ref{pf:qmin-id1} below, by following essentially 
the strategy of~\cite[\S2.3]{FZ}.
We first continue to review material from~\cite[\S9.2]{BZ}.

\begin{Lem} \label{lem:MaxDer}
Let $\psi\in A_q(\g)_{\ga,\de}$ and $\psi'\in A_q(\g)_{\ga',\de'}$. For
$i\in I$, assume that $a=(\a_i,\de)$ and $a'=(\a_i,\de')$ are the maximal
non-negative integers such that $f_i^a\cdot\psi\neq 0$ and 
$f_i^{a'}\cdot\psi'\neq 0$. Then 
\[
(f_i^{(a)}\cdot \psi)\,(f_i^{(a')}\cdot\psi')= f_i^{(a+a')}\cdot(\psi\,\psi').
\]
\end{Lem}
This  follows from the definition of the comultiplication
of $U_q(\g)$ and from~\eqref{bimod_alg}.
Note that we have an analogous result where $f_i$ is replaced by $e_i$ and acts
from the right.
The next lemma is an immediate consequence of Lemma~\ref{cons_def} and Lemma~\ref{lem:MaxDer}.
\begin{Lem} \label{lem:QMinIdRed}
Let $\la',\la''\in P_+$, $u', v', u'', v''\in W$, and $i,j\in I$, be such that
\[
l(s_j v')=l(v')+1,\quad l(s_j v'')=l(v'')+1,\quad l(s_i u')=l(u')+1,\quad l(s_i u'')=l(u'')+1.
\]
Then, putting $a = (\a_j,\, v'(\la')+v''(\la''))$,
and $b = (\a_i,\, u'(\la')+u''(\la''))$, we have 
\begin{align*}
f_j^{(a)} \cdot 
(\De_{u'(\la'),v'(\la')}\De_{u''(\la''),v''(\la'')}) &=
\De_{u'(\la'),s_jv'(\la')}\De_{u''(\la''),s_jv''(\la'')},\\
(\De_{u'(\la'),v'(\la')}\De_{u''(\la''),v''(\la'')}) \cdot
e_i^{(b)} &=
\De_{s_i u'(\la'),v'(\la')}\De_{s_i u''(\la''), v''(\la'')}.
\end{align*}
\end{Lem}
The next lemma follows easily from Lemma~\ref{lem:QMinIdRed} by induction on the length of $u$ and $v$.
\begin{Lem} \label{rem:MultQMin}
The quantum minors have the following multiplicative property:
\[
\De_{u(\la),v(\la)}\,\De_{u(\mu),v(\mu)}=\De_{u(\la+\mu),v(\la+\mu)},\quad
(u,v\in W,\ \la,\mu\in P_+).
\]
In particular, the factors of the second summand in the right hand side of (\ref{eq:2.10})
pairwise commute. 
\end{Lem}
The factors of the first summand in the right hand side of (\ref{eq:2.10})
also commute with each other, as shown by the following lemma.

\begin{Lem}\label{commute}
Suppose that for $u,v\in W$ and $i\in I$ we have $l(us_i)=l(u)+1$ and
$l(vs_i)=l(v)+1$. Then 
\[
\De_{us_i(\vpi_i),\,v(\vpi_i)}\,\De_{u(\vpi_i),\,vs_i(\vpi_i)}
=\De_{u(\vpi_i),\,vs_i(\vpi_i)}\,\De_{us_i(\vpi_i),\,v(\vpi_i)}. 
\]
\end{Lem}
\proof
For $x\in U_q(\nn)$ we have 
\[
x\cdot \De_{s_i(\vpi_i),\,\vpi_i} = x\cdot \De^{\vpi_i}\cdot e_i 
= \eps(x) \De^{\vpi_i}\cdot e_i = \eps(x) \De_{s_i(\vpi_i),\,\vpi_i}.
\]
Similarly, for $y\in U_q(\nn_-)$ we have 
\[
\De_{\vpi_i,\,s_i(\vpi_i)}\cdot y = \eps(y) \De_{\vpi_i,\,s_i(\vpi_i)}.
\] 
By \cite[Lemma~10.2]{BZ} we can then conclude that
\[ 
\De_{s_i(\vpi_i),\,\vpi_i}\,\De_{\vpi_i,\,s_i(\vpi_i)}=
\De_{\vpi_i,\,s_i(\vpi_i)}\,\De_{s_i(\vpi_i),\,\vpi_i}.
\] 
(In \cite{BZ}, $\g$ is generally assumed to be finite-dimensional, but this
assumption plays no role in the proof of \cite[Lemma~10.2]{BZ}.)
This proves the lemma for $u=v=e$. The general case then follows from 
successive applications of Lemma~\ref{lem:QMinIdRed}.
\cqfd

\subsection{Proof of Proposition~\ref{qmin-id1}} \label{pf:qmin-id1} 
With the help of Lemma~\ref{lem:QMinIdRed}, we see by an easy induction on 
the length of $u$ and $v$  that it is sufficient to verify the 
special case $u=v=e$, \ie
\begin{equation}
\underbrace{\vphantom{\prod_{j\neq i}}%
\De_{s_i(\vpi_i),s_i(\vpi_i)}\,\De_{\vpi_i,\vpi_i}-
q^{-1}\De_{s_i(\vpi_i),\vpi_i}\,\De_{\vpi_i,s_i(\vpi_i)}}_{\psi_1}=
\underbrace{\prod_{j\neq i}\De_{\vpi_j,\vpi_j}^{-a_{ji}}}_{\psi_2}.
\end{equation}
Note that
\[
\ga_i:=-\sum_{j\neq i} a_{ji}\vpi_j = 2\vpi_i-\a_i=\vpi_i+s_i(\vpi_i)\in P_+.
\]
Thus, by Lemma~\ref{rem:MultQMin} we have
\begin{equation} \label{eq:gamma_i}
\psi_2=\De^{\ga_i}=\De_{\ga_i,\ga_i},
\end{equation}
in particular the factors of $\psi_2$ commute. It follows from the definition
of $\De^{\ga_i}$ that
\begin{itemize}
\item[(1)] $\psi_2(1_U)=1$,
\item[(2)] $\psi_2\in A_q(\g)_{\ga_i,\ga_i}$,
\item[(3)] $e_j\cdot \psi_2=0=\psi_2\cdot f_j$ for all $j\in I$.
\end{itemize}
Here $1_U$ stands for the unit of $U_q(\g)$.
By Proposition~\ref{PW-thm} these properties characterize $\psi_2$ uniquely,
so it is sufficient to verify the  properties (1) -- (3) for $\psi_1$.
\medskip

{\bf Property (1).\ } We use that for any $\psi,\phi\in A_q(\g)$ we have
$(\psi\,\phi)(1_U)=\psi(1_U)\phi(1_U)$ since $\De(1_U)=1_U\otimes 1_U$, and
note that $q^0=1_U$.
Now, $\De_{\vpi_i,s_i(\vpi_i)}=f_i\cdot\De^{\vpi_i}$, thus
$\De_{\vpi_i,s_i(\vpi_i)}(1_U)=\De^{\vpi_i}(f_i)=0$.
Similarly, $\De_{s_i(\vpi_i),\vpi_i}(1_U)=0$ and $\De_{\vpi_i,\vpi_i}(1_U)=1$.
Finally, 
\[
\De_{s_i(\vpi_i),s_i(\vpi_i)}(1_U)=\De^{\vpi_i}(e_i\, f_i)=
{\De^{\vpi_i}(f_i\,1_U e_i)}+\De^{\vpi_i}\left(\frac{q^{h_i}-q^{-h_i}}{q-q^{-1}}\right)=0+1.
\]
Thus $\psi_1(1_U)=1$.
\medskip

{\bf Property (2).\ } Recall that $\De_{\ga,\de}\in A_q(\g)_{\ga,\de}$. Since
$A_q(\g)$ is $P\times P$-graded it follows from~\eqref{eq:gamma_i} that
both summands of $\psi_1$ belong to $A_q(\g)_{\ga_i,\ga_i}$.
\medskip

{\bf Property (3).\ } Since $\De(e_j)=e_j\otimes 1 +q^{h_j}\otimes e_j$ the
operator $e_j\cdot -$ acts on $A_q(\g)$ as a ``graded $q$-derivation'', \ie
for $\psi\in A_q(\g)_{\ga,\de}$ and $\psi'\in A_q(\g)_{\ga',\de'}$ we have
\begin{equation}
e_j\cdot (\psi\,\psi')= (e_j\cdot\psi)\,\psi'+ 
q^{\< h_j,\de\>}\psi\,(e_j\cdot\psi').
\end{equation}
Similarly,
\begin{equation}
(\psi\,\psi')\cdot f_j= \psi\ (\psi'\cdot f_j)+
q^{\< -h_j,\ga'\>}(\psi\cdot f_j)\,\psi'.
\end{equation}

Now, for $j\neq i$ we have $\< h_j,\vpi_i\>=0$ and
$\< h_j,s_i(\vpi_i)\>=\< h_j,\vpi_i -\a_i\>=-a_{ij}\geq 0$. Thus, 
by~\eqref{eq:leftz}, 
the left multiplication by $e_j$ annihilates 
all the minors appearing in 
$\psi_1$, and as a consequence $e_j\cdot\psi_1=0$.

It remains to show that $e_i\cdot\psi_1=0$. To this end we observe first 
that $e_i\cdot\De_{u(\vpi_i),s_i(\vpi_i)}=\De_{u(\vpi_i),\vpi_i}$ for any $u\in W$.
In fact, 
\begin{multline*}
e_i\cdot\De_{u(\vpi_i),s_i(\vpi_i)}=(e_i\,f_i)\cdot\De_{u(\vpi_i),\vpi_i}=
{(f_i\, e_i)\cdot\De_{u(\vpi_i),\vpi_i}}\ +\ 
\frac{q^{h_i}-q^{-h_i}}{q-q^{-1}}\cdot\De_{u(\vpi_i),\vpi_i}\\
=0+\frac{q^{\< h_i,\vpi_i\>}-q^{\< -h_i,\vpi_i\>}}{q-q^{-1}}\De_{u(\vpi_i),\vpi_i}
=\De_{u(\vpi_i),\vpi_i}.
\end{multline*}
So, we can now calculate
\begin{multline} \label{eq:eikillpsi}
e_i\cdot\psi_1= 
(e_i\cdot\De_{s_i(\vpi_i),s_i(\vpi_i)})\,\De_{\vpi_i,\vpi_i}-
q^{-1}(e_i\cdot\De_{s_i(\vpi_i),\vpi_i})\,\De_{\vpi_i,s_i(\vpi_i)}\\
+q^{\<h_i,s_i(\vpi)\>} \De_{s_i(\vpi_i),s_i(\vpi_i)})\,(e_i\cdot\De_{\vpi_i,\vpi_i})-
q^{\<h_i,\vpi_i\>-1}\De_{s_i(\vpi_i),\vpi_i}\,(e_i\cdot\De_{\vpi_i,s_i(\vpi_i)})\\
=\De_{s_i(\vpi_i),\vpi_i}\De_{\vpi_i,\vpi_i}-
q^0\De_{s_i(\vpi_i),\vpi_i}\De_{\vpi_i,\vpi_i}=0.
\end{multline}

Finally, we have to show that $\psi_i\cdot f_j=0$ for all $j\in I$.
Again, for $j\neq i$ we see by~\eqref{eq:rightz} that the map $x\mapsto x\cdot f_j$ annihilates
all the minors occuring in $\psi_1$. In order to see that
$\psi_1\cdot f_i=0$, we note that by Lemma~\ref{commute},
we have 
$\De_{s_i(\vpi_i),\vpi_i}\,\De_{\vpi_i,s_i(\vpi_i)}=
\De_{\vpi_i,s_i(\vpi_i)}\,\De_{s_i(\vpi_i),\vpi_i}$. 
Then we can proceed as in
\eqref{eq:eikillpsi}.
Proposition~\ref{qmin-id1} is proved.

\section{The quantum coordinate rings $A_q(\bb)$ and $A_q(\nn)$}

\subsection{The quantum coordinate ring $A_q(\bb)$}
Let $U_q(\bb)$ be the subalgebra of $U_q(\g)$ generated by $e_i\ (i\in I)$ and
$q^h\ (h\in P^*)$. We have
\begin{equation}
\De(U_q(\bb)) \subset U_q(\bb)\otimes U_q(\bb), \qquad
S(U_q(\bb)) \subseteq U_q(\bb), 
\end{equation}
hence $U_q(\bb)$ is a Hopf subalgebra. 
Therefore, as in (\ref{multA}), $U_q(\bb)^*$ has a multiplication
dual to the comultiplication of $U_q(\bb)$.
Clearly, the map $\rho\colon U_q(\g)^* \to U_q(\bb)^*$ 
given by restricting linear forms from $U_q(\g)$ to $U_q(\bb)$ is an
algebra homomorphism.

We define $A_q(\bb):= \rho(A_q(\g))$.
Let $Q_+ = \bigoplus_{i\in I} \N \a_i$.
For $\ga, \de \in P$, let $A_q(\bb)_{\ga,\de} = \rho(A_q(\g)_{\ga,\de})$.
Since $U_q(\bb)\subseteq \bigoplus_{\a \in Q_+} U_q(\g)_\a$, we have by Lemma~\ref{gradings},
\begin{equation}
A_q(\bb) = \bigoplus_{\ga-\de\in Q_+} A_q(\bb)_{\ga,\de}.
\end{equation}

\subsection{The quantum coordinate ring $A_q(\nn)$}\label{multAqn}
Recall that $U_q(\nn)$ is the subalgebra of $U_q(\g)$ generated by $e_i\ (i\in I)$.
Because of (\ref{comult}), this is \emph{not} a Hopf subalgebra. Nevertheless, we can endow $U_q(\nn)^*$
with a multiplication as follows. Recall that every $y\in U_q(\bb)$ can be written
as a $\Q(q)$-linear combination of elements of the form
$x\,q^h$ with $x\in U_q(\nn)$ and $h\in P^*$.
Given $\psi\in U_q(\nn)^*$, we define the linear form $\tpsi\in U_q(\bb)^*$ by 
\begin{equation}\label{tilde}
\tpsi(x\,q^h) = \psi(x),\qquad (x\in U_q(\nn),\ h\in P^*).
\end{equation}
Clearly, $\iota\colon \psi \mapsto \tpsi$ is an injective linear map.
Moreover, since $\De(x\,q^h) = \sum x_{(1)}q^h\otimes x_{(2)}q^h$, it follows 
immediately from (\ref{multA}) that
\begin{equation}
\left(\tpsi\cdot\tvphi\right)(x\,q^h) = \left(\tpsi\cdot\tvphi\right)(x),\qquad
(\psi, \vphi\in U_q(\nn)^*,\ x\in U_q(\nn),\ h\in P^*). 
\end{equation}
Therefore $\iota(U_q(\nn)^*)$ is a subalgebra of $U_q(\bb)^*$, and we can define
\begin{equation}\label{defmult}
\psi\cdot\vphi = \iota^{-1}(\tpsi\cdot\tvphi). 
\end{equation}

We have $U_q(\nn) = \bigoplus_{\a\in Q_+} U_q(\nn)_{\a}$, where 
$U_q(\nn)_{\a} = U_q(\nn) \cap U_q(\g)_{\a}$ is finite-dimensional for every $\a\in Q_+$.
Let 
\begin{equation}
A_q(\nn) = \bigoplus_{\a\in Q_+} \Hom_{\Q(q)}(U_q(\nn)_{\a},\Q(q))
= \bigoplus_{\a\in Q_+} A_q(\nn)_\a \subset U_q(\nn)^*
\end{equation}
denote the graded dual of $U_q(\nn)$.
It is easy to see that $A_q(\nn)$ is a subalgebra of $U_q(\nn)^*$ for the multiplication 
defined in (\ref{defmult}). Moreover, $\iota(A_q(\nn)) \subset A_q(\bb)$, and more precisely
(\ref{tilde}) shows that
\begin{equation}
\iota(A_q(\nn)_\a) = A_q(\bb)_{\a,0},\quad (\a\in Q_+). 
\end{equation}
To summarize, $A_q(\nn)$ can be identified with
the subalgebra $\iota(A_q(\nn))=\bigoplus_{\a\in Q_+} A_q(\bb)_{\a,0}$ of $A_q(\bb)$.

\subsection{The algebra isomorphism between $A_q(\nn)$ and $U_q(\nn)$}\label{isomAU}
For $i\in I$, let $\de_i \in \End_{\Q(q)}(U_q(\nn))$ be the $q$-derivation defined by
$\de_i(e_j)=\de_{ij}$ and
\begin{equation}\label{qderiv}
\de_i(xy)= \de_i(x)\,y + q^{\< h_i,\a\>}x\,\de_i(y),
\qquad
(x \in U_q(\nn)_\a,\ y\in U_q(\nn)). 
\end{equation}
It is well known that there exists a unique nondegenerate symmetric
bilinear form on $U_q(\nn)$ such that $(1,1)=1$ and
\begin{equation}
(\de_i(x),\,y) = (x,\, e_iy),\qquad (x \in U_q(\nn),\ y\in U_q(\nn),\ i\in I). 
\end{equation}
Denote by $\psi_x$ the linear form on $U_q(\nn)$ given by
$\psi_x(y) = (x,\,y)$. 
Then, the map $\Psi\colon x\mapsto \psi_x$ is an isomorphism of
graded vector spaces from $U_q(\nn)$ to $A_q(\nn)$. 

\begin{Prop}\label{isomPhi}
Let $A_q(\nn)$ be endowed with the multiplication (\ref{defmult}). 
Then $\Psi$ is an isomorphism of algebras from $U_q(\nn)$ to $A_q(\nn)$. 
\end{Prop}
\proof
We need to show that
\begin{equation}
\tpsi_{xy}(zq^h) = \left(\tpsi_x\cdot\tpsi_y\right)(zq^h),
\quad
(x,y,z\in U_q(\nn),\ h\in P^*). 
\end{equation}
By linearity, we can assume that $z=e_{i_1}\cdots e_{i_m}$ for some $i_1,\ldots,i_m\in I$.
By definition, 
\begin{equation}\label{calcul1}
\tpsi_{xy}(zq^h) = (xy,\,z) = (\de_{i_m}\cdots \de_{i_1}(xy),\,1).
\end{equation}
Let us assume, without loss of generality, that $x\in U_q(\nn)_\a$.
It follows from (\ref{qderiv}) that
\begin{equation}\label{prodder}
\de_{i_m}\cdots \de_{i_1}(xy) =  \sum_{K} q^{\si(J,K)} \de_{i_{j_s}}\cdots \de_{i_{j_1}}(x)\,
\de_{i_{k_r}}\cdots \de_{i_{k_1}}(y) 
\end{equation}
where the sum is over all subsets $K=\{k_1 < \cdots < k_r\}$ of $[1,m]$,
$J= \{j_1<\cdots < j_s\}$ is the complement of $K$ in $[1,m]$, and
\[
\si(J,K) =  \left\<h_{k_1},\,\a - \sum_{j\in J,\, j<k_1}\a_j\right\>
+ \cdots + 
\left\<h_{k_r},\,\a - \sum_{j\in J,\, j<k_r}\a_j\right\>.
\]
Moreover, a summand of the r.h.s. of (\ref{prodder}) can give a nonzero contribution
to (\ref{calcul1}) only if $\a = \sum_{j\in J} \a_j$. In this case we have
\begin{equation}\label{AJK}
\si(J,K) =  \left\<h_{k_1},\,\sum_{j\in J,\, j>k_1}\a_j\right\>
+ \cdots + 
\left\<h_{k_r},\,\sum_{j\in J,\, j>k_r}\a_j\right\>,
\end{equation}
and 
\begin{equation}
(xy,\,z) = \sum_{K} q^{\si(J,K)} (e_{i_{j_1}}\cdots e_{i_{j_s}},\,x)\,
(e_{i_{k_1}}\cdots e_{i_{k_r}},\,y). 
\end{equation}
On the other hand, it is easy to deduce from (\ref{comult}) that
\[
\De(z) = \sum_{K} q^{\si(J,K)} e_{i_{j_1}}\cdots e_{i_{j_s}}t_{i_{k_1}}\cdots t_{i_{k_r}}
\otimes
e_{i_{k_1}}\cdots e_{i_{k_r}} 
\]
where $\si(J,K)$ is again given by (\ref{AJK}).
It then follows from the definition of $\tpsi_x$ and $\tpsi_y$ that 
\[
\left(\tpsi_x\cdot\tpsi_y\right)(zq^h) = \left(\tpsi_x\cdot\tpsi_y\right)(z)
=
\sum_{K} q^{\si(J,K)} (e_{i_{j_1}}\cdots e_{i_{j_s}},\,x)\,
(e_{i_{k_1}}\cdots e_{i_{k_r}},\,y) =  \tpsi_{xy}(zq^h). 
\]
\cqfd

\section{Determinantal identities for unipotent quantum minors}

\subsection{Principal quantum minors}

Let $\la\in P_+$, and $u,v\in W$.
The quantum minor $\De_{u(\la),v(\la)}$ is called \emph{principal} when $u(\la)=v(\la)$.
In this case we have by Lemma~\ref{gradings}~(c) that $\De_{v(\la),v(\la)}(x) = 0$
if $x\not\in U_q(\g)_0$.
Therefore the restriction $\rho(\De_{v(\la),v(\la)}) \in  A_q(\bb)$ is given by
\begin{equation}
\rho(\De_{v(\la),v(\la)})(xq^h) = \eps(x)q^{\<h,v(\la)\>},\qquad (x\in U_q(\n),\ h\in P^*). 
\end{equation}
Define 
\begin{equation}
\De_{v(\la),v(\la)}^* := \De_{v(\la),v(\la)}\circ S \in U_q(\g)^*, 
\end{equation}
where $S$ is the antipode of $U_q(\g)$. 

\begin{Lem}\label{qcentral}
\begin{itemize}
\item[(a)] The principal quantum minors $\rho(\De_{v(\la),\,v(\la)})\ (v\in W,\, \la\in P_+)$
are invertible in $A_q(\bb)$, with inverse equal to $\rho(\De_{v(\la),v(\la)}^*)$.
\item[(b)] The principal quantum minors $\rho(\De_{v(\la),\,v(\la)})$
are $q$-central in $A_q(\bb)$. More precisely, for $\psi\in A_q(\bb)_{\ga,\de}$ we have
\[
\rho(\De_{v(\la),\,v(\la)})\cdot \psi = q^{(v(\la),\ga-\de)} \psi \cdot \rho(\De_{v(\la),\,v(\la)}). 
\]
\end{itemize}
\end{Lem}
\proof (a) First, we note that $\rho(\De_{v(\la),v(\la)}^*)$ belongs to $A_q(\bb)$.
Indeed, let $U_q(\bb_-)$ be the subalgebra of $U_q(\g)$ generated by $f_i\ (i\in I)$ and
$q^h\ (h\in P^*)$, and $U_q^0$ the subalgebra generated by 
$q^h\ (h\in P^*)$.
Since $S$ is an anti-automorphism of $U_q(\g)$ which stabilizes 
$U_q(\bb)$, $U_q(\bb_-)$, and $U_q^0$, it is clear that $\De_{v(\la),v(\la)}^*$
generates an integrable left (\resp right) submodule of $U_q(\g)^*$.
So $\De_{v(\la),v(\la)}^*\in A_q(\g)$, and $\rho(\De_{v(\la),v(\la)}^*)\in A_q(\bb)$. 
Now, it follows easily from the definition that 
\[
\De_{v(\la),v(\la)}^*(xq^h) = \eps(x) q^{-\<h,v(\la)\>}, \qquad (x\in U_q(\n),\ h\in P^*), 
\]
so that
\[
\left(\De_{v(\la),\,v(\la)}\cdot \De_{v(\la),v(\la)}^*\right)(xq^h)
=
\sum \De_{v(\la),\,v(\la)}(x_{(1)}q^h) \De_{v(\la),v(\la)}^*(x_{(2)}q^h)
=
\eps(xq^h),
\]
which proves (a).

(b) As for (a), it is enough to evaluate each side of the equation at a typical element $xq^h$ of $U_q(\bb)$.
Since the equation relates two elements of $A_q(\bb)_{\ga+v(\la),\de+v(\la)}$, we may assume that
$x\in U_{\ga - \de}$. By linearity, we may further assume that 
$x=e_{i_1}\cdots e_{i_k}$, where $\a_{i_1}+\cdots +\a_{i_k} = \ga-\de$, without loss of generality.
Using (\ref{comult}), we have
\[
\left(\De_{v(\la),\,v(\la)}\cdot \psi\right)(xq^h) = \sum \De_{v(\la),\,v(\la)}(x_{(1)}q^h) \psi(x_{(2)}q^h)
= \De_{v(\la),\,v(\la)}(t_{i_1}\cdots t_{i_k}q^h) \psi(xq^h),
\]
because $\De_{v(\la),\,v(\la)}(x_{(1)}q^h) \not = 0$ only if $x_{(1)} \in U_q(\g)_0$.
Hence
\[
\left(\De_{v(\la),\,v(\la)}\cdot \psi\right)(xq^h) =  q^{(v(\la),\,\ga-\de)}\De_{v(\la),\,v(\la)}(q^h) \psi(xq^h).
\]
On the other hand it also follows from (\ref{comult}) that
\[
\left(\psi\cdot\De_{v(\la),\,v(\la)} \right)(xq^h) = \sum \psi(x_{(1)}q^h)\De_{v(\la),\,v(\la)}(x_{(2)}q^h)
= \psi(xq^h) \De_{v(\la),\,v(\la)}(q^h),
\]
hence the result.
\cqfd

\subsection{Unipotent quantum minors}

The quantum minor $\De_{u(\la),\,v(\la)}$  
belongs to $A_q(\g)_{u(\la),v(\la)}$,
hence its restriction $\rho(\De_{u(\la),\,v(\la)})$ belongs to
$A_q(\bb)_{u(\la),\,v(\la)}$. Therefore, if $v(\la)\not = 0$, it does not belong to $\iota(A_q(\nn))$. 
But a slight modification does, as we shall now see.

We define the \emph{unipotent quantum minor} $D_{u(\la),\,v(\la)}$ by
\begin{equation}
D_{u(\la),\,v(\la)} = \rho\left(\De_{u(\la),\,v(\la)}\cdot \De_{v(\la),\,v(\la)}^*\right). 
\end{equation}
This is an element of $A_q(\bb)_{u(\la)-v(\la),\,0} \subset \iota(U_q(\nn)^*)$.
In fact, the same calculation as in the proof of Lemma~\ref{qcentral} shows that
for $x\in U_q(\nn)$ and $h\in P^*$, we have
\begin{equation}
D_{u(\la),\,v(\la)}(xq^h) = \De_{u(\la),\,v(\la)}(xq^h)\, \De_{v(\la),\,v(\la)}(q^{-h}) 
= \De_{u(\la),\,v(\la)}(x) = D_{u(\la),\,v(\la)}(x). 
\end{equation}
Thus, we can regard $D_{u(\la),\,v(\la)}$ as the restriction to $U_q(\nn)$ of the quantum
minor $\De_{u(\la),\,v(\la)}$. 
But we should be aware that this restriction is \emph{not} an algebra homomorphism.
For example, by Lemma~\ref{rem:MultQMin} we have
\[
\De_{u(\la),v(\la)}\cdot\De_{u(\mu),v(\mu)}= \De_{u(\mu),v(\mu)}\cdot\De_{u(\la),v(\la)} ,\quad
(u,v\in W,\ \la,\mu\in P_+).
\]
The corresponding commutation relation for unipotent quantum minors
is given by the following
\begin{Lem}\label{qcommute2}
For $u,v\in W$ and $\la,\mu\in P_+$ we have
\[
D_{u(\la),v(\la)}\cdot D_{u(\mu),v(\mu)}= q^{(v(\mu),\,u(\la))-(v(\la),\,u(\mu))}
D_{u(\mu),v(\mu)}\cdot D_{u(\la),v(\la)}.
\]
\end{Lem}
\proof Let us write for short $\rho(\De_{u(\la),\,v(\la)}) = \De_{u(\la),\,v(\la)}$
and $\rho(\De_{v(\la),\,v(\la)}^*)  = \De^*_{v(\la),\,v(\la)}$.
We have, by Lemma~\ref{qcentral},
\begin{eqnarray*}
D_{u(\la),v(\la)} D_{u(\mu),v(\mu)}&=& \De_{u(\la),\,v(\la)}\De^*_{v(\la),\,v(\la)} 
\De_{u(\mu),\,v(\mu)}\De^*_{v(\mu),\,v(\mu)}\\
&=&
q^{(v(\la),\,v(\mu)-u(\mu))}
\De_{u(\la),\,v(\la)}\De_{u(\mu),\,v(\mu)}
\De^*_{v(\la),\,v(\la)}\De^*_{v(\mu),\,v(\mu)}\\
&=&
q^{(v(\la),\,v(\mu)-u(\mu))}
\De_{u(\mu),\,v(\mu)}\De_{u(\la),\,v(\la)}
\De^*_{v(\mu),\,v(\mu)}\De^*_{v(\la),\,v(\la)}\\
&=&
q^{(v(\mu),\,u(\la))-(v(\la),\,u(\mu))}
\De_{u(\mu),\,v(\mu)}\De^*_{v(\mu),\,v(\mu)} 
\De_{u(\la),\,v(\la)}\De^*_{v(\la),\,v(\la)}\\
&=& q^{(v(\mu),\,u(\la))-(v(\la),\,u(\mu))} D_{u(\mu),v(\mu)} D_{u(\la),v(\la)} 
\end{eqnarray*}
\cqfd

Similarly, Lemma~\ref{commute} implies:
\begin{Lem}\label{qcommute}
 Suppose that for $u,v\in W$ and $i\in I$ we have $l(us_i)=l(u)+1$ and
$l(vs_i)=l(v)+1$. Then 
\[
D_{us_i(\vpi_i),\,v(\vpi_i)}\,D_{u(\vpi_i),\,vs_i(\vpi_i)}
=q^{(vs_i(\vpi_i),\,us_i(\vpi_i))-(v(\vpi_i),\,u(\vpi_i))}
D_{u(\vpi_i),\,vs_i(\vpi_i)}\,D_{us_i(\vpi_i),\,v(\vpi_i)}. 
\]
\cqfd
\end{Lem}
We can also regard the quantum unipotent minors
as linear forms on $U_q(\nn)$ given by matrix coefficients of integrable 
representations. 
Indeed, using (\ref{matrixCoeff}), we have
\begin{equation}
D_{u(\la),v(\la)}(x) = (m_{u(\la)},\,xm_{v(\la)})_\la,\qquad (x\in U_q(\nn)). 
\end{equation}
In particular, when $u=e$ we get the same quantum flag minors
\begin{equation}\label{matrix_coeff_flag}
D_{\la,v(\la)}(x) = (m_{\la},\,xm_{v(\la)})_\la,\qquad (x\in U_q(\nn)). 
\end{equation}
as in \cite[\S6.1]{Ki} (up to a switch from $U_q(\n_-)$ to $U_q(\n)$).

\subsection{A family of identities for unipotent quantum minors}

We are now in a position to deduce from Proposition~\ref{qmin-id1} an algebraic
identity satisfied by unipotent quantum minors.
Later on, we will see that particular cases of this identity can be seen as
quantum exchange relations in certain quantum cluster algebras.

As in \S\ref{pf:qmin-id1}, let us write $\ga_i = \vpi_i+s_i(\vpi_i)$, so that
\begin{equation}\label{factorqminor}
\De_{u(\ga_i),v(\ga_i)} = \prod_{j\neq i}\De_{u(\vpi_j),v(\vpi_j)}^{-a_{ji}}, 
\end{equation}
where, by Lemma~\ref{rem:MultQMin}, the order of the factors in the 
right-hand side is irrelevant.

\begin{Prop}\label{unip_minor_identity}
Suppose that for $u,v\in W$ and $i\in I$ we have $l(us_i)=l(u)+1$ and
$l(vs_i)=l(v)+1$. Then 
\[
q^{A}\,
D_{us_i(\vpi_i),\,vs_i(\vpi_i)}\,D_{u(\vpi_i),\,v(\vpi_i)}=
q^{-1+B}\,
D_{us_i(\vpi_i),\,v(\vpi_i)}\,D_{u(\vpi_i),\,vs_i(\vpi_i)}
+
D_{u(\ga_i),v(\ga_i)}
\]
holds in $A_q(\nn)$, where
\[
A =  (vs_i(\vpi_i),\,u(\vpi_i)-v(\vpi_i)),\qquad
B =  (v(\vpi_i),\,u(\vpi_i)-vs_i(\vpi_i)).
\]
\end{Prop}
\proof
Again let us write for short $\rho(\De_{u(\la),\,v(\la)}) = \De_{u(\la),\,v(\la)}$
and $\rho(\De^*_{v(\la),\,v(\la)})  = \De^*_{v(\la),\,v(\la)}$.
We apply the restriction homomorphism $\rho$ to the equality of Proposition~\ref{qmin-id1},
and we multiply both sides from the right
by
\[
\De^*_{v(\vpi_i),\,v(\vpi_i)}\De^*_{v(s_i(\vpi_i)),\,v(s_i(\vpi_i))}
= \De^*_{v(\ga_i),v(\ga_i)}.  
\]
Note that all these minors commute by Lemma~\ref{rem:MultQMin}.
The result then follows directly from the definition of unipotent quantum minors, and from
Lemma~\ref{qcentral}~(b).
\cqfd

It is sometimes useful to write the second summand of the right-hand side of Proposition~\ref{unip_minor_identity}
as a product. It is straightforward to deduce from (\ref{factorqminor}) and Lemma~\ref{qcentral}~(b)
that we have
\begin{equation}\label{prodD}
D_{u(\ga_i),v(\ga_i)}
=
q^C\,
\prod_{j\not = i}^{\longrightarrow} \left(D_{u(\vpi_j),\,v(\vpi_j)}
\right)^{-a_{ij}}, 
\end{equation}
where
\begin{equation}\label{eqC}
C = 
\sum\limits_{\substack{j<k\\ j\not= i\not = k}}
a_{ij}a_{ik}(v(\vpi_j),\,u(\vpi_k)-v(\vpi_k))
+ \sum_{j\not = i} \left(\begin{matrix}-a_{ij} \cr 2\end{matrix}\right) (v(\vpi_j),\,u(\vpi_j)-v(\vpi_j)). 
\end{equation}

\subsection{A quantum $T$-system}\label{ssectTsystem}

Let $\ii = (i_1,\ldots,i_r)\in I^r$ be such that $l(s_{i_1}\cdots s_{i_r})=r$.
We will now deduce from Proposition~\ref{unip_minor_identity} a system of identities
relating the unipotent quantum minors
\begin{equation}
D(k,l;j) := D_{s_{i_1}\cdots s_{i_k}(\vpi_j),\ s_{i_1}\cdots s_{i_l}(\vpi_j)},
\qquad
(0\le k < l \le r,\ j\in I). 
\end{equation}
Here, we use the convention that $D(0,l;j)= D_{\vpi_j,\ s_{i_1}\cdots s_{i_l}(\vpi_j)}$,
a \emph{quantum flag minor}.
This system can be viewed as a $q$-analogue of a $T$-system, (see \cite{KNS}).
It will allow us to express every quantum minor $D(k,l;j)$ in terms of the
flag minors $D(0,m;i)$. 

Note that, because of (\ref{sivpi}), every quantum minor $D(k,l;j)$ is equal to a minor of
the form $D(b,d;j)$ where $i_b=i_d=j$.
When this is the case, we can simply write $D(b,d;j) = D(b,d)$. 
Note that in particular, $D(b,b)=1$ for every $b$.
By convention, we write $D(0,b) = D_{\vpi_{i_b},\ s_{i_1}\cdots s_{i_b}(\vpi_{i_b})}$.
We will also use the following shorthand notation:
\begin{eqnarray}
b^-(j) &:=& \max\left(\{s < b \mid i_s = j\}\cup\{0\}\right),\\
b^- &:=& \max\left(\{s < b \mid i_s = i_b\}\cup\{0\}\right),\\
\mu(b,j) &:=&  s_{i_1}\cdots s_{i_b}(\vpi_j).\label{equamu}
\end{eqnarray}
In (\ref{equamu}) we understand that $\mu(0,j) = \vpi_j$. 
Clearly, we have $D(b,d;j) = D(b^-(j),d^-(j))$.
\begin{Prop}\label{Tsystem}
Let $1\le b<d\le r$ be such that $i_b = i_d = i$. 
There holds
\begin{equation}\label{eqTsystem}
q^A\, D(b,d) D(b^-,d^-) = 
q^{-1+B}\, D(b,d^-)D(b^-,d)
\ +\ q^C
\prod_{j\not = i}^{\longrightarrow}D(b^-(j),d^-(j))^{-a_{ij}}
\end{equation}
where
\[
A = (\mu(d,i),\,\mu(b^-,i)-\mu(d^-,i)),\qquad
B = (\mu(d^-,i),\,\mu(b^-,i)-\mu(d,i)), 
\]
and 
\[
C 
= 
\sum\limits_{\substack{j<k\\ j\not= i\not = k}}
a_{ij}a_{ik}\left(\mu(d,j),\,\mu(b,k)-\mu(d,k)\right)
+\ \sum_{j\not = i} \left(\begin{matrix}-a_{ij} \cr 2\end{matrix}\right) 
\left(\mu(d,j),\ \mu(b,j)-\mu(d,j)\right). 
\]
\end{Prop}
\proof
This follows directly from Proposition~\ref{unip_minor_identity}, (\ref{prodD}), and (\ref{eqC}), by taking 
\[
u = s_{i_1}\cdots s_{i_{b-1}},\qquad
v = s_{i_1}\cdots s_{i_{d-1}},\qquad
i=i_b=i_d. 
\]
\cqfd 

\section{Canonical bases} 

\subsection{The canonical basis of $U_q(\nn)$}\label{def_can}

We briefly review  Lusztig's definition of a canonical basis of $U_q(\nn)$.

Recall the scalar product $(\cdot,\cdot)$ on $U_q(\n)$ defined in \S\ref{isomAU}.
In \cite[Chapter 1]{Lu}, Lusztig defines a similar scalar product $(\cdot,\cdot)_L$, 
using the same $q$-derivation $\de_i$ (denoted by ${}_ir$ in \cite[1.2.13]{Lu}) but with
a different normalization $(e_i,e_i)_L = (1-q^{-2})^{-1}$.
It it easy to see that $(x,y)=0$ if and only if $(x,y)_L=0$, and if $x,y\in U_q(\n)_\b$ then
\begin{equation}\label{scalarL}
(x,y)_L = (1-q^{-2})^{-\deg\b} (x,y), 
\end{equation}
where, for $\b = \sum_i c_i\a_i$, we set $\deg\b = \sum_ic_i$. This slight difference
will not affect the definition of the canonical basis below, 
and we will always use $(\cdot,\cdot)$ instead of $(\cdot,\cdot)_L$.

Let $\A = \Q[q,q^{-1}]$.
We introduce the $\A$-subalgebra $U_\A(\n)$ of $U_q(\n)$ generated by the divided powers
$e_i^{(k)}\ (i\in I,\ k\in \N)$. We define a ring automorphism $x \mapsto \overline{x}$ of 
$U_q(\g)$ by 
\begin{equation}
 \overline{q} = q^{-1},\quad \overline{e_i} = e_i,\quad  \overline{f_i}=f_i, \qquad (i\in I).
\end{equation}
This restricts to a ring automorphism of $U_q(\n)$.

The canonical basis $\B$ is an $\A$-basis of  $U_\A(\n)$ such that
\begin{equation}\label{bar}
 \overline{b}=b,\qquad (b\in \B).
\end{equation}
Moreover, for every $b, b'\in \B$ the scalar product $(b,b')\in \Q(q)$
has no pole at $q=\infty$, and
\begin{equation}\label{almost_ortho}
(b,\, b')|_{q=\infty} = \delta_{b,b'}. 
\end{equation}
By this we mean that 
\[
(b,b') = \frac{a_jq^j+\cdots a_1q+a_0}{a'_kq^k+\cdots+ a'_1q + a'_0}, 
\qquad (a_i,a'_i \in \Z,\ a_j\not = 0,\ a'_k\not = 0), 
\]
with $j<k$ when $b\not = b'$, and $j=k, a_j=a'_k$ when $b=b'$.
It is easy to see that if an $\A$-basis of   
$U_\A(\n)$ satisfies (\ref{bar}) and (\ref{almost_ortho}),
then it is unique up to sign (see \cite[14.2]{Lu}).
The existence of $\B$ is proved in \cite[Part 2]{Lu}, and a consistent choice of signs
is provided. 
Of course, $\B$ is also a $\Q(q)$-basis of~$U_q(\n)$.

\subsection{The dual canonical basis of $U_q(\nn)$}\label{dualcan}

Let $\B^*$ be the basis of $U_q(\n)$ adjoint to $\B$ with respect to the
scalar product $(\cdot,\,\cdot)$. We call it the \emph{dual canonical basis}
of $U_q(\n)$, since 
it can be identified via $\Psi$ with the dual basis of $\B$ in $A_q(\n)$.

Note that $\B^*$ is not invariant under the bar automorphism $x\mapsto \overline{x}$.
The property of $\B^*$ dual to (\ref{bar}) can be stated as follows.
Let $\si$ be the composition of the anti-automorphism $*$ and the
bar involution, that is, $\si$ is the ring \emph{anti}-automorphism 
of $U_q(\nn)$ such that
\begin{equation}
\si(q) = q^{-1},\qquad \si(e_i) = e_i. 
\end{equation}
For $\b \in Q_+$, define
\begin{equation}\label{defN(b)}
N(\b) := \frac{(\b,\,\b)}{2} - \deg\b. 
\end{equation}
Then, if $b \in U_q(\nn)_\b$
belongs to $\B^*$, there holds
\begin{equation}\label{eq:sigma_b}
\si(b) = q^{N(\b)}\, b, 
\end{equation}
(see \cite{Re,Ki}).

\subsection{Specialization at $q=1$ of $U_q(\n)$ and $A_q(\n)$}\label{specq1}

Recall that $U_\A(\n)$ is the $\A$-submodule of $U_q(\n)$ spanned by
the canonical basis $\B$. 
If we regard $\C$ as an $\A$-module via the homomorphism $q \mapsto 1$,
we can define 
\begin{equation}
 U_1(\n) := \C \otimes_\A U_\A(\n). 
\end{equation}
This is a $\C$-algebra isomorphic to the enveloping algebra $U(\n)$.

Similarly, let $A_\A(\n)$ be the $\A$-submodule of $A_q(\n)$ spanned by
the basis $\Psi(\B^*)$. 
Define 
\begin{equation}
 A_1(\n) := \C \otimes_\A A_\A(\n). 
\end{equation}
This is a $\C$-algebra isomorphic to the graded dual $U(\n)^*_{\rm gr}$.
This commutative ring can be identified with the coordinate ring $\C[N]$
of a pro-unipotent pro-group $N$ with Lie algebra the completion $\widehat{\n}$
of $\n$ (see \cite{GLS}).

\subsection{Global bases of $U_q(\n_-)$}

We shall also use Kashiwara's lower global basis $\B^{\rm low}$
of $U_q(\n_-)$, constructed in \cite{K1}. It was proved by Grojnowski and Lusztig that
$\vphi(\B) = \B^{\rm low}$, where $\vphi$ is the anti-automorphism
of (\ref{vphi}).

For $i\in I$, we introduce the $q$-derivations $e'_i$ and ${}_ie'$ of 
$U_q(\n_-)$, defined by $e'_i(f_j)={}_ie'(f_j) = \de_{ij}$
and, for homogeneous elements $x, y \in U_q(\n_-)$,
\begin{eqnarray}\label{qderiv-}
e'_i(xy)&=& e'_i(x)\,y + q^{\< h_i,{\rm wt}(x)\>}x\,e'_i(y),\\[2mm]
{}_ie'(xy)&=& q^{\< h_i,{\rm wt}(y)\>}{}_ie'(x)\,y + x\,\,{}_ie'(y).
\end{eqnarray}
Note that ${}_ie' = *\circ e'_i\circ *$.
Let us denote by $(\cdot,\,\cdot)_K$ the Kashiwara scalar product
on $U_q(\n_-)$. It is the unique symmetric bilinear form such that
$(1,\, 1)_K = 1$, and 
\begin{equation}
(f_i x,\,y)_K = (x,\,e'_i(y))_K,\qquad (x \in U_q(\nn_-),\ y\in U_q(\nn_-),\ i\in I). 
\end{equation}  
It also satisfies
\begin{equation}
(x f_i,\,y)_K = (x,\,{}_ie'(y))_K,\qquad (x \in U_q(\nn_-),\ y\in U_q(\nn_-),\ i\in I). 
\end{equation} 
Let $\bvphi$ be the composition of $\vphi$ and the bar involution, that is, 
$\bvphi$ is the ring \emph{anti}-automorphism 
of $U_q(\g)$ such that
\begin{equation}
\bvphi(q) = q^{-1},\quad \bvphi(e_i) = f_i, 
\quad \bvphi(f_i) = e_i,
\quad \bvphi(q^h) = q^{-h}. 
\end{equation}
The following lemma expresses the compatibility between the scalar products and 
$q$-derivations on $U_q(\n)$ and $U_q(\n_-)$.
\begin{Lem}\label{lem+-}
\begin{itemize}
 \item[(a)] For $i\in I$, we have ${}_ie'\circ \bvphi = \bvphi \circ \de_i$.
 \item[(b)] For $x, y \in U_q(\n)$ we have $\overline{(x,\,y)} = (\bvphi(x),\,\bvphi(y))_K$.
\end{itemize}
\end{Lem}
\proof
As ${}_ie'\circ \bvphi$ and $\bvphi \circ \de_i$ are both linear, it is enough to prove that 
${}_ie'\circ \bvphi(z) = \bvphi \circ \de_i(z)$ for
any homogeneous element $z$ of $U_q(\nn)$. 
We prove this by induction on the degree of $z$. If the degree is
1, then this follows easily from the definition of all these maps. 
Then assume that the degree
of $z$ is bigger than one. Then, without loss of generality, we can assume that 
$z = xy$ with degrees of $x$ and $y$ smaller
than the degree of $z$. Now
\[
{}_ie'(\bvphi(xy))= {}_ie'(\bvphi(y)\bvphi(x))
= q^{-\< h_i,\b\>}{}_ie'(\bvphi(y))\,\bvphi(x) + \bvphi(y)\,\,{}_ie'(\bvphi(x)).
\]
By induction on the degrees of $x$ and $y$ we can assume that 
${}_ie'(\vphi(x))=\bvphi(\de_i(x))$ and ${}_ie'(\bvphi(y))=\bvphi(\de_i(y))$,
so that
\[
{}_ie'(\bvphi(xy))=\bvphi\left(q^{\< h_i,\b\>}x\de_i(y) + \de_i(x)y)\right) = \bvphi(\de_i(xy)), 
\]
which proves (a). Then
\[
(\bvphi(e_ix),\bvphi(y))_K = (\bvphi(x)f_i,\bvphi(y))_K=(\bvphi(x),{}_ie'(\bvphi(y)))_K
=(\bvphi(x),\bvphi(\de_i(y))_K. 
\]
By induction on the degrees of $x$ and $y$ we can assume that 
$(\bvphi(x),\bvphi(\de_i(y)))_K=\overline{(x,\de_i(y))}=\overline{(e_ix,y)}$, which proves (b).
\cqfd

Let $\B^{\rm up}$ denote the upper global basis of $U_q(\n_-)$.
This is the basis adjoint to $\B^{\rm low}$ with respect to $(\cdot,\,\cdot)_K$.
By Lemma~\ref{lem+-}(b), we also have $\B^{\rm up} = \bvphi(\B^*)$.

Let us denote by $\BB(\infty)$ the crystal of $U_q(\n_-)$, and by $b_\infty$ its 
highest weight element.
The elements of $\B^{\rm low}$ (\resp $\B^{\rm up}$)
are denoted by $G^{\rm low}(b)\ (b\in \BB(\infty))$ (\resp $G^{\rm up}(b)$).
As usual, for every $i\in I$, one denotes by $\te_i$ and $\tf_i$ the Kashiwara
crystal operators of $\BB(\infty)$. We will need the following well known 
property of the upper global basis, (see \cite[Lemma 5.1.1]{K}, \cite[Corollary 3.16]{Ki}):
\begin{Lem}\label{rmax}
Let $b\in \BB(\infty)$, and put $k=\max\{j\in\N\mid ({}_ie')^j(G^{\rm up}(b)) \not = 0\}$.
Then, denoting by $({}_ie')^{(k)}$ the $k$th $q$-divided power of ${}_ie'$, we have 
$({}_ie')^{(k)}(G^{\rm up}(b))\in \B^{\rm up}$.
More precisely, 
\[
 ({}_ie')^{(k)}(G^{\rm up}(b)) = G^{\rm up}((\te_i^*)^k b),
\]
where $\te_i^*$ is the crystal operator obtained from $\te_i$ by conjugating with the involution $*$
of (\ref{star}).
\end{Lem}
The integer $k=\max\{j\in\N\mid (e'_i)^j(G^{\rm up}(b)) \not = 0\} = \max\{j\in\N\mid (\te_i)^j(b) \not = 0\}$
is denoted by $\eps_i(b)$.
Similarly, the integer 
$k=\max\{j\in\N\mid ({}_ie')^j(G^{\rm up}(b)) \not = 0\} = \max\{j\in\N\mid (\te_i^*)^j(b) \not = 0\}$
is denoted by $\eps_i^*(b)$.

\subsection{Unipotent quantum minors belong to $\B^*$}

Using the isomorphism $\Psi\colon U_q(\n) \to A_q(\n)$ of \ref{isomAU}, we can 
regard the unipotent quantum minors $D_{u(\la),v(\la)}$ as elements of $U_q(n)$.
More precisely, let $d_{u(\la),v(\la)} = \Psi^{-1}(D_{u(\la),v(\la)})$ be the
element of $U_q(\n)$ such that
\begin{equation}
D_{u(\la),v(\la)}(x) = (d_{u(\la),v(\la)},\,x),\qquad (x\in U_q(\n)). 
\end{equation}
By a slight abuse, we shall also call $d_{u(\la),v(\la)}$ a unipotent quantum minor.
In this section we show: 
\begin{Prop}\label{qmincan}
For every $\la\in P_+$ and $u,v \in W$ such that $u(\la)-v(\la)\in Q_+$, the
unipotent quantum minor $d_{u(\la),v(\la)}$ belongs to $\B^*$.
More precisely, writing $u=s_{i_{l(u)}}\cdots s_{i_1}$, $v=s_{j_{l(v)}}\cdots s_{j_1}$ 
and defining $b_k$ and $c_l$ as in (\ref{eq-bc}), we have
\[
d_{u(\la),v(\la)} = 
\bvphi\left(G^{\rm up}\left((\te_{i_{l(u)}}^*)^{b_{l(u)}}\cdots (\te_{i_1}^*)^{b_1}\tf_{j_{l(v)}}^{c_{l(v)}}\cdots \tf_{j_1}^{c_1}b_\infty\right)\right).
\]
\end{Prop}
\proof
We proceed in two steps, and first consider the case when $u=e$ is the unit of $W$,
that is, the case of unipotent quantum flag minors $d_{\la,v(\la)}$.
By (\ref{matrix_coeff_flag}), we have for $x\in U_q(\n)$,
\[
(d_{\la,v(\la)},\,x) = (m_\la,\,xm_{v(\la)})_\la = (\vphi(x)m_\la,\,m_{v(\la)})_\la. 
\]
It is well known that the extremal weight vectors $m_{v(\la)}$ belong to Kashiwara's
lower global basis of $V(\la)$, and also to the upper global basis. 
More precisely, 
we have 
\[
m_{v(\la)} = G^{\rm low}\left(\tf_{j_{l(v)}}^{c_{l(v)}}\cdots \tf_{j_1}^{c_1}b_\la\right)
= G^{\rm up}\left(\tf_{j_{l(v)}}^{c_{l(v)}}\cdots \tf_{j_1}^{c_1}b_\la\right),
\]
where $b_\la$ is the highest weight element of the crystal $\BB(\la)$ of $V(\la)$.
Hence we have
\begin{equation}\label{eq.5.14}
(d_{\la,v(\la)},\,x) = 
\left(\vphi(x)m_\la,\, G^{\rm up}\left(\tf_{j_{l(v)}}^{c_{l(v)}}\cdots \tf_{j_1}^{c_1}b_\la\right)\right)_\la.   
\end{equation}
It follows from (\ref{eq.5.14}) that $(d_{\la,v(\la)},\,x) = 1$ if 
$\vphi(x) = G^{\rm low}\left(\tf_{j_{l(v)}}^{c_{l(v)}}\cdots \tf_{j_1}^{c_1}b_\infty\right)$,
and $(d_{\la,v(\la)},\,x) = 0$ if $\vphi(x)$ is any other element of $\B^{\rm low}$.
Therefore, $d_{\la,v(\la)}$ is the element of $\B^*$ given by
\[
d_{\la,v(\la)} = \bvphi\left(G^{\rm up}\left(\tf_{j_{l(v)}}^{c_{l(v)}}\cdots \tf_{j_1}^{c_1}b_\infty\right)\right). 
\]
Now let us consider the general case when $u = s_{i_{l(u)}}\cdots s_{i_1}$ is non trivial.
For $k=1,\ldots, l(u)$, write $u_k = s_{i_{k}}\cdots s_{i_1}$. 
Using Lemma~\ref{cons_def}, we have for $x\in U_q(\n)$, 
\[
D_{u_k(\la),v(\la)}(x) = \De_{u_k(\la),v(\la)}(x) = \left(\De_{u_{k-1}(\la),v(\la)}\cdot e_{i_k}^{(b_k)}\right)(x)
\]
where $b_k=(\a_{i_k},u_{k-1}(\la)) = \max\{j \mid \De_{u_{k-1}(\la),v(\la)}\cdot e_{i_k}^{(j)} \not= 0\}$.
Now we can also write
\[
D_{u_k(\la),v(\la)}(x) = \left( d_{u_{k-1}(\la),v(\la)},\, e_{i_k}^{(b_k)}x\right) = \left( (\de_{i_k})^{(b_k)}d_{u_{k-1}(\la),v(\la)},\, x\right),
\]
hence
\[
d_{u_k(\la),v(\la)} = (\de_{i_k})^{(b_k)}d_{u_{k-1}(\la),v(\la)},
\]
where $(\de_i)^{(b)}$ means the $b$th $q$-divided power of the $q$-derivation $\de_i$.
Since $u_k(\lambda)-v(\la)\in Q_+$, the restriction $\rho(\De_{u_k(\la),v(\la)})$ is nonzero,
hence $b_k= \max\{j \mid (\de_{i_k})^{(b)}d_{u_{k-1}(\la),v(\la)}\not= 0\}$.
Applying $\bvphi$ and assuming by induction on $k$ that 
\[
d_{u_{k-1}(\la),v(\la)} = 
\bvphi\left(G^{\rm up}\left((\te_{i_{k-1}}^*)^{b_{k-1}}\cdots (\te_{i_1}^*)^{b_1}\tf_{j_{l(v)}}^{c_{l(v)}}\cdots \tf_{j_1}^{c_1}b_\infty\right)\right)
\]
we get by Lemma~\ref{lem+-}(a) that 
\[
\bvphi(d_{u_k(\la),v(\la)}) =  
({}_{i_{k}}e')^{(b_k)}\left(G^{\rm up}\left((\te_{i_{k-1}}^*)^{b_{k-1}}\cdots (\te_{i_1}^*)^{b_1}\tf_{j_{l(v)}}^{c_{l(v)}}\cdots \tf_{j_1}^{c_1}b_\infty\right)\right)
\]
and
\[
b_k= \eps_{i_k}^*\left((\te_{i_{k-1}}^*)^{b_{k-1}}\cdots (\te_{i_1}^*)^{b_1}\tf_{j_{l(v)}}^{c_{l(v)}}\cdots \tf_{j_1}^{c_1}b_\infty)\right).
\]
Thus, applying Lemma~\ref{rmax}, we get that
\[
\bvphi(d_{u_k(\la),v(\la)}) =  
G^{\rm up}\left((\te_{i_{k}}^*)^{b_{k}}\cdots (\te_{i_1}^*)^{b_1}\tf_{j_{l(v)}}^{c_{l(v)}}\cdots \tf_{j_1}^{c_1}b_\infty\right),
\]
and the statement follows by induction on $k$.
\cqfd

\section{Quantum unipotent subgroups}

In this section we provide a quantum version of the coordinate
ring $\C[N(w)]$ studied in \cite{GLS}, following \cite{Lu,S,Ki}. 

\subsection{The quantum enveloping algebra $U_q(\n(w))$}
We fix $w\in W$, and we denote by $\De_w^+$ the subset of positive roots $\a$
of $\g$ such that $w(\a)$ is a negative root. This gives rise to a finite-dimensional
Lie subalgebra
\[
\n(w) := \bigoplus_{\a\in \De_w^+} \n_\a 
\]
of $\n$, of dimension $l(w)$. The graded dual $U(\n(w))^*_{\rm gr}$ can be identified
with the coordinate ring $\C[N(w)]$ of a unipotent subgroup $N(w)$ of the Kac-Moody group
$G$ with ${\rm Lie}(N(w)) = \n(w)$. (For more details, see \cite{GLS}.) 

In order to define a $q$-analogue of $U(\n(w))$, one introduces Lusztig's braid group
operation on $U_q(\g)$ \cite{Lu}. For $i\in I$, Lusztig has proved the existence
of a $\Q(q)$-algebra automorphism $T_i$ of $U_q(\g)$ satisfying
\begin{eqnarray}
T_i(q^h) &=& q^{s_i(h)},\\
T_i(e_i) &=& -t_i^{-1}f_i,\\
T_i(f_i) &=& -e_it_i,\\
T_i(e_j) &=& \sum_{r+s=-\<h_i,\a_j\>}(-1)^rq^{-r}e_i^{(r)}e_je_i^{(s)}\qquad (j\not=i),\\
T_i(f_j) &=& \sum_{r+s=-\<h_i,\a_j\>}(-1)^rq^{r}f_i^{(s)}f_jf_i^{(r)}\qquad (j\not=i).
\end{eqnarray}
(This automorphism is denoted by $T'_{i,-1}$ in \cite{Lu}.)
For a fixed reduced decomposition $w=s_{i_r}\cdots s_{i_1}$, let us set, as in 
(\ref{defbeta}), 
\begin{equation}\label{def_beta}
\b_k = s_{i_1}\cdots s_{i_{k-1}}(\a_{i_k}),\qquad (1\le k \le r).
\end{equation}
Then $\De_w^+ = \{\b_1,\ldots,\b_r\}$.
We define following Lusztig, the corresponding quantum root vectors:
\begin{equation}
E(\b_k) := T_{i_1}\cdots T_{i_{k-1}}(e_{i_k}),\qquad (1\le k \le r). 
\end{equation}
It is known that $E(\b_k) \in U_q(\n)_{\b_k}$. For $\aa = (a_1,\ldots,a_r)\in \N^r$,
set 
\begin{equation}
E(\aa) := E(\b_1)^{(a_1)}\cdots E(\b_r)^{(a_r)}, 
\end{equation}
where $E(\b_k)^{(a_k)}$ denotes the $a_k$th $q$-divided power of $E(\b_k)$.
Lusztig has shown that the subspace of $U_q(\n)$ spanned by 
$\{E(\aa)\mid \aa\in\N^r\}$ is independent of the choice of the reduced word
$\ii=(i_r,\ldots,i_1)$ for~$w$. We denote it by $U_q(\n(w))$.
Moreover, 
\begin{equation}
\P_\ii:=\{E(\aa)\mid \aa\in\N^r\}
\end{equation}
is a basis of $U_q(\n(w))$, which we call the \emph{PBW-basis} attached to $\ii$.

In fact, $U_q(\n(w))$ is even a subalgebra of $U_q(\n)$. This follows from a formula
due to Leven\-dor\-skii-Soibelman (see \cite[4.3.3]{Ki}).

\subsection{The quantum coordinate ring $A_q(\n(w))$}

Using the algebra isomorphism $\Psi\colon U_q(\n) \to A_q(\n)$ of \ref{isomAU},
we can define $A_q(\n(w)) := \Psi(U_q(\n(w))$. This is a subalgebra of $A_q(\n)$.

Lusztig \cite[38.2.3]{Lu} has shown that $\P_\ii$ is orthogonal, that is $(E(\aa),E(\bbf))=0$ if 
$\aa \not = \bbf$. Moreover
\begin{equation}
(E(\b_k),E(\b_k))= (1-q^{-2})^{\deg\b_k-1},\qquad (1\le k\le r), 
\end{equation}
and
\begin{equation}
(E(\aa),E(\aa)) = \prod_{k=1}^r \frac{(E(\b_k),E(\b_k))^{a_k}}{\{a_k\}!}, 
\end{equation}
where by definition
\begin{equation}
\{a\}! = \prod_{j=1}^a \frac{1-q^{-2j}}{1-q^{-2}}. 
\end{equation}
Denote by $\P^*_\ii$ the basis  
\begin{equation}
E^*(\aa) = \frac{1}{(E(\aa),E(\aa))}\ E(\aa),\qquad (\aa \in \N^r)
\end{equation}
of $U_q(\n(w))$ adjoint to $\P_\ii$.
We call $\P^*_\ii$ the \emph{dual PBW-basis} of $U_q(\n(w))$ since 
it can be identified via $\Psi$ with the basis of $A_q(\n(w))$ dual to $\P_\ii$.
In particular we have the dual PBW generators:
\begin{equation}
E^*(\b_k) = (1-q^{-2})^{-\deg\b_k+1}E(\b_k),\qquad (1\le k\le r). 
\end{equation}

\subsection{Action of $T_i$ on unipotent quantum minors}

\begin{Prop}\label{Tminor}
Let $\la \in P^+$, and $u, v \in W$ be such that $u(\la) - v(\la) \in Q_+$,
and consider the unipotent quantum minor $d_{u(\la),v(\la)}$.
Suppose that $l(s_iu) = l(u)+1$, and $l(s_iv)= l(v)+1$. Then
\[
T_i\left( d_{u(\la),v(\la)}\right) 
= (1-q^{-2})^{(\a_i,v(\la)-u(\la))}\,d_{s_iu(\la),\,s_iv(\la)}. 
\]
\end{Prop}

The proof will use Proposition~\ref{qmincan} and the following lemmas.

\begin{Lem}\label{Tphi}
We have $T_i \circ \bvphi = \bvphi\circ T_i$. 
\end{Lem}
\proof This follows immediately from the definitions of $\bvphi$ and of $T_i$.
\cqfd

The next lemma is a restatement of a result of Kimura \cite[Theorem 4.20]{Ki},
based on previous results of Saito \cite{S} and Lusztig \cite{Lu2}.
Note that our $T_i$ is denoted by $T_i^{-1}$ in \cite{Ki}.
\begin{Lem}\label{lemKimura}
Let $b\in\BB(\infty)$ be such that $\eps_i(b) = 0$. Then
\[
T_i(G^{\rm up}(b)) = 
(1-q^2)^{(\a_i,{\rm wt}(b))}\,
G^{\rm up}\left(\tf_i^{\ \vphi_i^*(b)}(\te^*_i)^{\eps_i^*(b)} b\right),
\] 
where $\vphi_i^*(b) := \eps_i^*(b) + (\a_i, {\rm wt}(b))$.
\end{Lem}

\medskip\noindent {\it Proof of Proposition~\ref{Tminor} --- \ }
By Lemma~\ref{Tphi} and Proposition~\ref{qmincan}, we have
\[
\bvphi\left(T_i(d_{u(\la),v(\la)})\right) = 
T_i\left(G^{\rm up}\left((\te_{i_{l(u)}}^*)^{b_{l(u)}}\cdots (\te_{i_1}^*)^{b_1}\tf_{j_{l(v)}}^{\ c_{l(v)}}\cdots \tf_{j_1}^{\ c_1}b_\infty\right)\right).
\]
Let us write for short $b := (\te_{i_{l(u)}}^*)^{b_{l(u)}}\cdots (\te_{i_1}^*)^{b_1}\tf_{j_{l(v)}}^{\ c_{l(v)}}\cdots \tf_{j_1}^{\ c_1}b_\infty \in \BB(\infty)$.
Then 
\[
\eps_i(b) = \max\{s \mid e_i^s\cdot \De_{u(\la),v(\la)} \not = 0 \}. 
\]
The assumption $l(s_iv)=l(v)+1$ implies that $(\a_i, v(\la))\ge 0$, so 
by (\ref{eq:leftz}), we have $\eps_i(b) = 0$, and we
are in a position to apply Lemma~\ref{lemKimura}.
Because of the assumption $l(s_iu)=l(u)+1$
we have 
\[
\eps^*_i(b) = \max\{s \mid \De_{u(\la),v(\la)}\cdot e_i^s \not = 0 \} = (\a_i,u(\la)), 
\] 
and 
\[
\vphi_i^*(b) = \eps^*_i(b) + (\a_i,{\rm wt}(b)) = (\a_i,u(\la)) + (\a_i,v(\la)-u(\la)) = (\a_i,v(\la)). 
\]
Thus, again by Proposition~\ref{qmincan}, we have 
\begin{align*}
T_i(d_{u(\la),v(\la)})&= \bvphi\left(T_i(G^{\rm up}b)\right)\\
&=\bvphi\left((1-q^2)^{(\a_i,v(\la)-u(\la))} 
G^{\rm up}\left({\tf_i}^{\ \vphi_i^*(b)}(\te^*_i)^{\eps_i^*(b)} b\right)\right)\\
&= (1-q^{-2})^{(\a_i,v(\la)-u(\la))}
\bvphi\left(
G^{\rm up}\left((\te^*_i)^{\eps_i^*(b)}(\te_{i_{l(u)}}^*)^{b_{l(u)}}\cdots (\te_{i_1}^*)^{b_1}{\tf_i}^{\ \vphi_i^*(b)}\tf_{j_{l(v)}}^{\ c_{l(v)}}\cdots \tf_{j_1}^{\ c_1}b_\infty\right)\right)\\
&= (1-q^{-2})^{(\a_i,v(\la)-u(\la))}d_{s_iu(\la),\,s_iv(\la)}.
\end{align*}
\cqfd

\subsection{Dual PBW generators are unipotent quantum minors}

Recall from (\ref{def_beta}) the definition  of the roots $\b_k$.
\begin{Prop}\label{PBWminors}
For $k = 1,\ldots,r$, we have
\[
E^*(\b_k) = d_{s_{i_1}\cdots s_{i_{k-1}}(\vpi_{i_k}),\ s_{i_1}\cdots s_{i_{k}}(\vpi_{i_k})}. 
\]
\end{Prop}
\proof
We have $e_{i_k} = d_{\vpi_{i_k},\ s_{i_k}(\vpi_{i_k})}$,
hence
\[
 E(\b_k) = T_{i_1}\cdots T_{i_{k-1}}\left(d_{\vpi_{i_k},\ s_{i_k}(\vpi_{i_k})}\right).
\]
Applying $k-1$ times Proposition~\ref{Tminor}, we get
\[
 E(\b_k) = (1-q^{-2})^N d_{s_{i_1}\cdots s_{i_{k-1}}(\vpi_{i_k}),\ s_{i_1}\cdots s_{i_{k}}(\vpi_{i_k})},
\]
where
\[
 N = -(\a_{i_{k-1}},\a_{i_k}) - (\a_{i_{k-2}},\, s_{i_{k-1}}(\a_{i_k})) - \cdots
-(\a_{i_1},\, s_{i_2}\cdots s_{i_{k-1}}(\a_{i_k})).
\]
Now,
\begin{align*}
\b_k &=  s_{i_1}\cdots s_{i_{k-1}}(\a_{i_k})\\
& = \a_{i_k} + \sum_{n=1}^{k-1} (s_{i_n}\cdots s_{i_{k-1}}(\a_{i_k})-s_{i_{n+1}}\cdots s_{i_{k-1}}(\a_{i_k}))\\
& = \a_{i_k} - \sum_{n=1}^{k-1} (\a_{i_n},\, s_{i_{n+1}}\cdots s_{i_{k-1}}(\a_{i_k}))\a_{i_n},
\end{align*}
so that $\deg\b_k = 1 + N$. Hence,
\[
d_{s_{i_1}\cdots s_{i_{k-1}}(\vpi_{i_k}),\ s_{i_1}\cdots s_{i_{k}}(\vpi_{i_k})}
= (1-q^{-2})^{1-\deg\b_k}\, E(\b_k) = E^*(\b_k).
\]
\cqfd

\subsection{The skew field $F_q(\nn(w))$}\label{skewfield}

It is well known that $A_q(\nn(w))$ is an Ore domain.
This follows for example from the fact that it has polynomial growth
(see \cite[Appendix]{BZ}).
Hence $A_q(\nn(w))$ embeds in its skew field of fractions, which
we shall denote by $F_q(\nn(w))$.

Recall the shorthand notation of \S\ref{ssectTsystem} for unipotent quantum minors.

\begin{Prop}\label{flagFq}
Let $0\le b < d \le r$ be such that $i_b = i_d$. Then the unipotent quantum
minor $D(b,d)$
belongs to $F_q(\nn(w))$. In particular, the quantum flag minors
$D(0,d)= D_{\vpi_{i_d},\ s_{i_1}\cdots s_{i_{d}}(\vpi_{i_d})}$ belong to $F_q(\nn(w))$.
\end{Prop}

\proof 
If $b=d^-$, by Proposition~\ref{PBWminors}, we have $D(d^-,d)= \Psi(E^*(\b_{d}))$,
so 
\[
D(d^-,d)\in A_q(\n(w)) \subset F_q(\n(w)).
\]
Recall the determinantal identity (\ref{eqTsystem}).
Arguing as in \cite[Corollary 13.3]{GLS}, we can order the set of minors
$D(b,d)\ (0\le b < d \le r)$ so that (i) the minors $D(d^-,d)$
are the smallest elements, and (ii) the minor $D(b^-,d)$
is strictly bigger than all the other minors occuring in 
(\ref{eqTsystem}). 
This allows to express, recursively $D(b^-,d)$ as a rational 
expression in the minors $D(c^-,c)\ (1\le c \le r)$,
and shows that $D(b^-,d)\in F_q(\n(w))$.
\cqfd 

We will show later (see Corollary~\ref{minorpol}) that, in fact, all quantum minors $D(b,d)$ 
are \emph{polynomials} in the dual PBW-generators $D(c^-,c)$.
Hence, they belong to $A_q(\n(w))$.

\subsection{Specialization at $q=1$ of $A_q(\n(w))$}\label{specialq1}

Let $A_\A(\n(w))$ denote the free $\A$-submodule of $A_q(\n(w))$
with basis $\Psi(\P^*_\ii)$. This integral form of $A_q(\n(w))$ is
an $\A$-algebra, independent of the choice of the reduced word $\ii$. 
Moreover, if we regard $\C$ as an $\A$-module via the homomorphism $q \mapsto 1$,
we can define 
\begin{equation}
 A_1(\n(w)) := \C \otimes_\A A_\A(\n(w)). 
\end{equation}
This is a $\C$-algebra isomorphic to the coordinate ring $\C[N(w)]$, (see \cite[Theorem 4.39]{Ki}).
In particular, if $D_{u(\la),v(\la)}$ is a unipotent quantum minor
in $A_q(\n(w))$, the element $1\otimes D_{u(\la),v(\la)}$ can be identified
with the corresponding classical minor in $\C[N(w)]$.

\section{Quantum cluster algebras}

In this section we recall, following Berenstein and Zelevinsky \cite{BZ}, the definition 
of a quantum cluster algebra.

\subsection{Based quantum tori}
Let $L = (\la_{ij})$ be a skew-symmetric $r\times r$-matrix with integer entries.
The {\em based quantum torus} $\T(L)$
is the $\Z[q^{\pm 1/2}]$-algebra generated by
symbols $X_1,\ldots,X_r,X_1^{-1},\ldots,X_r^{-1}$
submitted to the relations
\begin{equation}\label{qT}
X_iX_i^{-1}=X_i^{-1}X_i=1,\qquad
X_i X_j = q^{\la_{ij}} X_j X_i, \qquad (1\le i,j\le r).
\end{equation} 
For $\aa=(a_1,\ldots,a_r)\in \Z^r$, set
\begin{equation}\label{eqX^a}
X^\aa:=q^{\frac{1}{2}\sum_{i>j}a_ia_j\la_{ij}}X_1^{a_1}\cdots X_r^{a_r}.
\end{equation} 
Then $\{X^\aa\mid \aa\in \Z^r\}$ is a $\Z[q^{\pm 1/2}]$-basis
of $\T(L)$, and we have for $\aa, \bbf\in\Z^r$,
\begin{equation}\label{multBZ}
X^\aa X^\bbf=q^{\frac{1}{2}\sum_{i>j}(a_ib_j-b_ia_j)\la_{ij}}X^{\aa+\bbf}
= q^{\sum_{i>j}(a_ib_j-b_ia_j)\la_{ij}}X^\bbf X^\aa.
\end{equation} 
Since $\T(L)$ is an Ore domain, it can be embedded in its skew field of fractions $\F$.

\subsection{Quantum seeds}
Fix a positive integer $n< r$.
Let $\tB = (b_{ij})$ be an $r\times (r-n)$-matrix with integer coefficients.
The submatrix $B$ consisting of the first $r-n$ rows of $\tB$ is called the
\emph{principal part} of $\tB$. We will require $B$ to be skew-symmetric.
We call $\tB$ an \emph{exchange matrix}.
We say that the pair $(L,\tB)$ is \emph{compatible} if we have
\begin{equation}\label{def_compatible}
\sum_{k=1}^r b_{kj}\la_{ki} = \de_{ij}d,\qquad (1\le j\le r-n,\ 1\le i\le r) 
\end{equation}
for some positive integer $d$.

If $(L,\tB)$ is compatible, the datum $\SC=((X_1,\ldots,X_r),L,\tB)$ is called a 
\emph{quantum seed} in $\F$. The set $\{X_1,\ldots,X_r\}$ is called the \emph{cluster}
of $\SC$, and its elements the \emph{cluster variables}. The cluster variables
$X_{r-n+1},\ldots, X_r$ are called \emph{frozen variables}, since they will not be
affected by the operation of mutation to be defined now.
The elements $X^\aa$ with $\aa\in\N^r$ are called \emph{quantum cluster monomials}.

\subsection{Mutations}
For $k=1,\ldots,r-n$, we define the \emph{mutation} $\mu_k(L,\tB)$ of 
a compatible pair $(L,\tB)$. 
Let $E$ be the $r\times r$-matrix with entries
\begin{equation}
e_{ij} = 
\begin{cases}
\de_{ij} &\text{ if } j\not = k,\\
-1 &\text{ if } i=j=k,\\
\max(0,-b_{ik}) &\text{ if } i\not = j = k.
\end{cases}
\end{equation}
Let $F$ be the $(r-n)\times (r-n)$-matrix with entries
\begin{equation}
f_{ij} = 
\begin{cases}
\de_{ij} &\text{ if } i\not = k,\\
-1 &\text{ if } i=j=k,\\
\max(0,b_{kj}) &\text{ if } i = k \not = j.
\end{cases}
\end{equation}
Then $\mu_k(L,\tB) = (\mu_k(L),\mu_k(\tB))$ where
\begin{equation}
\mu_k(L) := E^{T} L E, \qquad \mu_k(\tB) := E\tB F. 
\end{equation}
Note that the mutation $\mu_k(\tB)$ of the exchange
matrix is a reformulation of the classical one defined in \cite{FZ2}. 
It is easy to check that $\mu_k(L,\tB)$ is again a compatible pair,
with the same integer $d$ as in~(\ref{def_compatible}).
Define $\aa' = (a'_1,\ldots,a'_r)$ and $\aa'' = (a''_1,\ldots,a''_r)$ by
\begin{equation}\label{qmut1}
a'_i = 
\begin{cases}
-1 &\text{ if } i = k,\\
\max(0,b_{ik}) &\text{ if } i\not = k,  
\end{cases}
\qquad
a''_i = 
\begin{cases}
-1 &\text{ if } i = k,\\
\max(0,-b_{ik}) &\text{ if } i\not = k.
\end{cases}
\end{equation}
One then defines
\begin{equation}\label{qmut2}
\mu_k(X_i) = 
\begin{cases}
X^{\aa'} + X^{\aa''},&\text{ if } i = k,\\
X_i&\text{ if } i\not = k.
\end{cases}
\end{equation}
Berenstein and Zelevinsky show that the elements $X'_i := \mu_k(X_i)$
satisfy 
\begin{equation}
X'_i X'_j = q^{\la'_{ij}} X'_j X'_i, \qquad (1\le i,j\le r), 
\end{equation}
where $\mu_k(L) = (\la'_{ij})$. Moreover they form a \emph{free generating set}
of $\F$, that is, one can write $X'_i = \theta(X^{\cc^i})$ where $\theta$ is a 
$\Q(q^{1/2})$-linear automorphism of $\F$,
and $(\cc^1,\ldots,\cc^r)$ is a $\Z$-basis of $\Z^r$.  
Therefore 
\begin{equation}
\mu_k(\SC):=((\mu_k(X_1),\ldots,\mu_k(X_r)),\ \mu_k(L),\ \mu_k(\tB))
\end{equation}
is a new quantum seed in $\F$,
called the \emph{mutation of $\SC$ in direction $k$}.
Moreover,
the mutation operation is involutive, that is, $\mu_k(\mu_k(\SC))=\SC$.

\begin{Def}
The quantum cluster algebra $\AA_{q^{1/2}}(\SC)$ is the $\Z[q^{\pm 1/2}]$-subalgebra of
the skew field~$\F$ generated by the union of clusters of all quantum seeds
obtained from $\SC$ by any sequence of mutations. 
\end{Def}

The following basic result is called the \emph{quantum Laurent phenomenon}.
\begin{Thm}[\cite{BZ}]
The quantum cluster algebra $\AA_{q^{1/2}}(\SC)$ is contained in the based quantum torus
generated by the quantum cluster variables of any given quantum seed $\SC'$ mutation equivalent to $\SC$. 
\end{Thm}

In the next sections we are going to construct a class of quantum cluster
algebras attached to some categories of representations of preprojective algebras.

\section{The category $\CC_w$}

We recall the main facts about the category $\CC_w$ and its maximal rigid objects,
following \cite{BIRS, GLSUni1, GLS}.

\subsection{The preprojective algebra}

Let $Q$ be a finite connected quiver without 
oriented cycles, with vertex set $I$, and
arrow set $Q_1$. 
We can associate with $Q$ a symmetric generalized Cartan matrix $A=(a_{ij})_{i,j\in I}$,
where $a_{ij}=2$ if $i=j$, and otherwise $a_{ij}$ is minus the number of edges between
$i$ and $j$ in the underlying unoriented graph. 
We will assume that $A$ is the Cartan matrix attached to the Kac-Moody algebra $\g$, that is, 
\begin{equation}
a_{ij} = \<h_i,\a_j\> = (\a_i,\a_j),\qquad (i,j\in I).
\end{equation}

Let $\C\overline{Q}$ be the path algebra of the 
{\it double quiver } $\overline{Q}$ 
of $Q$, which is obtained
from $Q$ by adding to each arrow $a\colon i \to j$ in $Q_1$ an arrow
$a^*\colon j \to i$ pointing in the opposite direction.
Let $(c)$ be
the two-sided ideal of $\C\overline{Q}$ generated by the element
\[
c = \sum_{a \in Q_1} (a^*a - aa^*).
\]
The algebra
$$
\L := \C\overline{Q}/(c)
$$ 
is called the {\it preprojective algebra} of $Q$.
Recall that for all $X,Y \in \md(\L)$ we have
\begin{equation}
\dm \Ext_\L^1(X,Y) = \dm \Ext_\L^1(Y,X).
\end{equation}
This follows for example from the following important formula
\begin{equation}\label{CBform}
\dim\Ext^1_\L(X,Y)=\dim\Hom_\L(X,Y) + \dim\Hom_\L(Y,X) - (\dimv\, X,\ \dimv\, Y),
\end{equation}
where we identify the dimension vector $\dimv\, X$ with an element of $Q_+$
in the standard way.

We denote by $\hI_i\ (i\in I)$ the indecomposable injective $\L$-module 
with socle $S_i$.
Here, $S_i$ is the one-dimensional simple $\L$-module
supported on the vertex $i$ of $Q$.
Note that the modules $\hI_i$ are infinite-dimensional if $Q$ is not
a Dynkin quiver.

For a $\L$-module $M$, let
$\soc_{(j)}(M)$ be the sum of all submodules
$U$ of $M$ with $U \cong S_j$.
For $(j_1,\ldots,j_s)\in I^s$, there is a unique sequence
$$
0 = M_0 \subseteq M_1 \subseteq \cdots \subseteq M_s \subseteq M
$$
of submodules of $M$ such that $M_p/M_{p-1} = \soc_{(j_p)}(M/M_{p-1})$.
Define $\soc_{(j_1,\ldots,j_s)}(M) := M_s$.
(In this definition, we do not assume that $M$ is finite-dimensional.)

\subsection{The subcategory $\CC_w$}
Let $\ii = (i_r,\ldots,i_1)$ be a reduced expression of $w\in W$.
For $1 \le k \le r$, let 
\begin{equation}
V_k := V_{\ii,k} := \soc_{(i_k,\ldots,i_1)}\left(\hI_{i_k}\right), 
\end{equation}
and set
$V_\ii := V_1 \oplus \cdots \oplus V_r$.
The module $V_\ii$ is dual to the cluster-tilting object
constructed in \cite[Section II.2]{BIRS}.
Define
$$
\CC_\ii := \Fac(V_\ii) \subseteq \md(\L).
$$
This is the full subcategory of $\md(\L)$ whose objects are
quotient modules of a direct sum of a finite number of copies of $V_\ii$.
For $j\in I$, let $k_j := \max\{ 1 \le k \le r \mid i_k = j \}$.
Define 
$I_{\ii,j} := V_{\ii,k_j}$
and set
$$
I_\ii := I_{\ii,1} \oplus \cdots \oplus I_{\ii,n}.
$$
The category $\CC_\ii$ and the module $I_\ii$ depend only on $w$, and not
on the chosen reduced expression $\ii$ of $w$.
Therefore, we define 
$$
\CC_w := \CC_\ii, \qquad I_w := I_\ii.
$$

\begin{Thm}[{\cite[Theorem II.2.8]{BIRS}}]
\label{main1}
For any $w\in W$,
the following hold:
\begin{itemize}

\item[(a)]
$\CC_w$ is a Frobenius category. Its
stable category $\underline{\CC}_w$ is a 2-Calabi-Yau category.

\item[(b)]
The indecomposable $\CC_w$-projective-injective modules are
the indecomposable direct summands of $I_w$.

\item[(c)]
$\CC_w = \Fac(I_w)$.

\end{itemize}
\end{Thm}

Note that in the case when $Q$ is a Dynkin quiver, that is, $\g$ is a simple Lie algebra of type $A, D, E$,
and $w=w_0$ is the longest element in $W$, then $\CC_{w_0} = \md(\L)$.

\subsection{Maximal rigid objects}\label{maximal_rigid}

For a $\L$-module $T$, we denote by $\add(T)$ the additive closure
of $T$, that is, the full subcategory of $\md(\L)$ whose objects
are isomorphic to direct sums of direct summands of $T$.
A $\L$-module $T$ is called {\it rigid} if $\Ext_\L^1(T,T) = 0$.
Let $T\in \CC_w$ be rigid. We say that
\begin{itemize}
\item
$T$ is $\CC_w$-{\it maximal rigid} if 
$
\Ext_\L^1(T \oplus X,X) = 0
$
with $X \in \CC_w$ implies $X \in \add(T)$;

\item
$T$ is a $\CC_w$-{\it cluster-tilting module} if 
$
\Ext_\L^1(T,X) = 0
$
with $X \in \CC_w$
implies $X \in \add(T)$.
\end{itemize}

\begin{Thm}[{\cite[Theorem I.1.8]{BIRS}}]
\label{main9}
For a rigid $\L$-module $T$ in $\CC_w$ the following are equivalent:
\begin{itemize}

\item[(a)]
$T$ has $r$ pairwise non-isomorphic indecomposable direct summands;

\item[(b)]
$T$ is $\CC_w$-maximal rigid;

\item[(c)]
$T$ is a $\CC_w$-cluster-tilting module.

\end{itemize}
\end{Thm}
Note that, given Theorem~\ref{main1}, the proof of \cite[Theorem 2.2]{GLSRigid}
carries over to this more general situation (see \cite[Theorem 2.2]{GLSUni1} in the case when $w$ is adaptable).

\subsection{The quiver $\G_T$ and the matrix $\tB_T$}\label{subsecGammaT}
Let $T = T_1\oplus \cdots \oplus T_r$ be a $\CC_w$-maximal rigid module, with each
summand $T_i$ indecomposable. Clearly each indecomposable $\CC_w$-injective module 
$I_{\ii,j}$ is isomorphic to one of the $T_k$'s, so, up to relabelling, we can assume that 
$T_{r-n+j}\cong I_{\ii,j}$ for $j=1,\ldots, n$.
Consider the endomorphism algebra $A_T:=\End_\L(T)^{\rm op}$. This is a basic algebra,
with indecomposable projective modules 
\begin{equation}
P_{T_i} := \Hom_\L(T,T_i)\ (1\le i\le r).
\end{equation} 
The simple $A_T$-modules are the heads $S_{T_i}$ of the projectives $P_{T_i}$.
One defines a quiver $\G_T$ with vertex set $\{1,\ldots,r\}$, and 
$d_{ij}$ arrows from $i$ to $j$, where $d_{ij} = \dim \Ext^1_{A_T}(S_{T_i},S_{T_j})$. 
(This is known as the Gabriel quiver of $A_T$.) 
Most of the information contained 
in $\G_T$ can be encoded in an $r\times (r-n)$-matrix $\tB_T = (b_{ij})_{1\le i\le r,\ 1\le j\le r-n}$, 
given by
\begin{equation}
b_{ij} = (\text{number of arrows } j\to i \text{ in } \G_T)
- (\text{number of arrows } i\to j \text{ in } \G_T).
\end{equation}
Note that $\tB_T$ can be regarded as an exchange matrix, with skew-symmetric principal part.

The next theorem gives an explicit description of the quiver $\G_{T}$ 
(hence also of the matrix $\tB_{T}$) for certain $\CC_w$-maximal rigid modules.
Following \cite{BFZ}, we define a quiver
$\G_\ii$ as follows.
The vertex set of $\G_\ii$ is equal to $\{1,\ldots,r\}$.
For $1 \le k \le r$, let 
\begin{align}\label{kminus}
k^- &:= \max\left(\{0\}\cup \{1 \le s \le k-1 \mid i_s = i_k \}\right),\\
k^+ &:= \min\left(\{ k+1 \le s \le r \mid i_s = i_k \} \cup \{r+1\}\right).
\end{align}
For $1 \le s,\,t \le r$ such that $i_s \not = i_t$, 
there are $|a_{i_s,i_t}|$ arrows from $s$ to $t$ provided
$t^+ \ge s^+ > t > s$.
These are called the {\it ordinary arrows} of $\G_\ii$.
Furthermore, for each $1 \le s \le r$ there is an arrow $s \to s^-$
provided $s^- > 0$.
These are the {\it horizontal arrows} of $\G_\ii$.

The following result generalizes \cite[Theorem 1]{GLSAus} (see \cite[Theorem 2.3]{GLSUni1} in the case when $w$ is adaptable).

\begin{Thm}[{\cite[Theorem II.4.1]{BIRS}}]
\label{main6}
The $\L$-module $V_\ii$ is a $\CC_w$-maximal rigid module, 
and we have $\G_{V_\ii} = \G_\ii$.
\end{Thm}

\subsection{Mutations of maximal rigid objects}

We consider again an arbitrary $\CC_w$-maximal rigid module $T$, and
we use the notation of \ref{subsecGammaT}. 

\begin{Thm}[{\cite[Theorem I.1.10]{BIRS}}]\label{th_mutation}
Let $T_k$ be a non-projective indecomposable direct summand of $T$. 
\begin{itemize}
\item[(a)] There exists a unique indecomposable module $T_k^*\not\cong T_k$
such that $(T/T_k)\oplus T_k^*$ is $\CC_w$-maximal rigid.
We call $(T/T_k)\oplus T_k^*$ the \emph{mutation of $T$ in direction $k$},
and denote it by $\mu_k(T)$.  
\item[(b)]
We have $\tB_{\mu_k(T)} = \mu_k(\tB_T)$.
\item[(c)] We have $\dim\Ext^1_\L(T_k,T_k^*) = 1$.
Let $0\to T_k \to T'_k \to T_k^* \to 0$ and 
$0\to T_k^* \to T''_k \to T_k \to 0$ be non-split short exact sequences.
Then 
\[
 T'_k \cong \bigoplus_{b_{jk}<0} T_j^{-b_{jk}},\qquad
 T''_k \cong \bigoplus_{b_{jk}>0} T_j^{b_{jk}}.
\]
\end{itemize}
\end{Thm}
Note that again, given Theorem~\ref{main1}, the proof of \cite[\S2.6]{GLSRigid}
carries over to this more general situation.

\subsection{The modules $M[l,k]$}

In this section we assume that the reduced word $\ii$ is fixed.
So for simplicity we often omit it from the notation. 
Thus, for $k = 1,\ldots,r$, we may write $V_k$ instead of $V_{\ii,k}$. 
Moreover, we use the convention that $V_0 = 0$. 

For $1\le k\le l\le r$ such that $i_k=i_l=i$, we have a natural embedding 
of $\L$-modules $V_{k^-} \to V_{l}$.
Following \cite[\S9.8]{GLS}, we define $M[l,k]$ as the cokernel 
of this embedding, that is,
\begin{equation}
 M[l,k] := V_{l} / V_{k^-}.
\end{equation}
In particular, we set $M_k:= M[k,k]$, and
\begin{equation}
M=M_\ii := M_1\oplus \cdots \oplus M_r. 
\end{equation} 
We will use the convention that $M[l,k]=0$ if $k>l$.
Every module $M[l,k]$ is indecomposable and rigid.
But note that $M$ is \emph{not} rigid.
Define 
\begin{align}\label{kmin}
k_{\min} &:= \min\{1 \le s \le r \mid i_s = i_k \},\\
k_{\max} &:= \max\{ 1 \le s \le r \mid i_s = i_k \}.
\end{align}
Then $V_k = M[k,k_{\min}]$ corresponds to an \emph{initial interval}.
The direct sum of all modules $M[k_{\max},k]$ corresponding to 
\emph{final intervals} is also a $\CC_w$-maximal rigid module, denoted 
by $T_\ii$.

By \cite[\S10]{GLS}, for every module $X\in\CC_w$, there exists
a chain
\begin{equation}
0 = X_0 \subseteq X_1 \subseteq \cdots \subseteq X_r = X 
\end{equation}
of submodules of $X$ with $X_k/X_{k-1} \cong M_k^{m_k}$, for some 
uniquely determined non-negative integers $m_k$.
The $r$-tuple $\mm(X) := (m_1, \ldots, m_r)$ will be called the
\emph{$M$-dimension vector} of $X$.

\subsection{Dimensions of Hom-spaces}

Let $R=R_\ii=(r_{kl})$ the $r\times r$-matrix with entries
\begin{equation}
r_{kl} =
\begin{cases}
0& \text{ if } k<l,\\
1& \text{ if } k=l,\\
(\b_k,\b_l)& \text{ if } k>l,
\end{cases}
\end{equation}

\begin{Prop} \label{prp:HomDim}
Suppose that $X,Y\in \CC_w$ satisfy $\Ext^1_\L(X,Y)=0$. Then 
\[
\dim\Hom_\L(X,Y)= \mm(X)\, R\, \mm(Y)^T.
\]
\end{Prop}
\proof
By \cite[Proposition~10.5]{GLS},
for all $Y\in\CC_w$ with $M$-dimension vector $\mm(Y)$,
we have 
\begin{equation}
\dim\Hom_\L(V_k,Y)=\dim\Hom_\L\left(V_k,\ \bigoplus_k M_k^{(\mm(Y))_k}\right).
\end{equation}
Moreover, the $M$-dimension vector of $V_k$ is given by 
\begin{equation}\label{aaV}
(\mm(V_k))_j=\begin{cases} 1 &\text{ if } i_j=i_k\text{ and } j\leq k,\\
                               0 &\text{ else.}
\end{cases}
\end{equation}
Therefore, we can restate~\cite[Lemma~9.8]{GLS} as
\begin{equation} \label{eq:Lem}
\dim\Hom_\Lam(V_k,Y)= \mm(V_k)\, R\,\mm(Y)^T.
\end{equation}
Now, for $X\in\CC_w$ we can find a short exact sequence
$0\ra V''\ra V'\ra X\ra 0$  with  $V',V''\in\add(V_\ii)$.
Since $\Hom_\Lam(V_\ii,-)$ is exact on this sequence we conclude that
\begin{equation}\label{eq:xvv}
\mm(X)=\mm(V')-\mm(V'').
\end{equation}
Finally, since $\Ext^1_\L(X,Y)=0$, the sequence
\[
0\ra\Hom_\Lam(X,Y)\ra\Hom_\Lam(V',Y)\ra\Hom_\Lam(V'',Y)\ra 0
\]
is exact. Thus we can calculate
\begin{align*}
\dim\Hom_\Lam(X,Y) &=\dim\Hom_\Lam(V',Y)-\dim\Hom_\Lam(V'',Y)\\
&
=(\mm(V')-\mm(V''))\, R\, \mm(Y)^T\\
&
=\mm(X)\, R\, \mm(Y)^T.
\end{align*}
\cqfd

\begin{Lem} \label{lem8.6}
Let $1\le b < d \le r$ be such that $i_b = i_d=i$.
There holds
\[
\dim\Hom_\Lam(M[d,b^+],\,M[d^-,b])=
\mm(M[d,b])\, R\, \mm(M[d^-,b])^T.
\]
\end{Lem}

\proof
We have a short exact sequence
\[
0\ra M[d^-,b]\ra M[d,b]\oplus M[d^-,b^+]\ra M[d,b^+]\ra 0
\]
with
\begin{equation} \label{eq:MutAdd}
\mm(M[d^-,b])-(\mm(M[d,b])+\mm(M[d^-,b^+]))+
\mm(M[d,b^+])= 0.
\end{equation}
Since 
\[
\Ext^1_\L(M[d,b^+],M[d^-,b]) = 1, \quad 
\Ext^1_\L(M[d,b^+],M[d,b]) = \Ext^1_\L(M[d,b^+],M[d^-,b^+]) = 0,
\]
this yields an exact sequence
\begin{multline*}
0\ra\Hom_\Lam(M[d,b^+],M[d^-,b])\ra\Hom_\Lam(M[d,b^+],M[d,b]\oplus M[d^-,b^+])
\ra \\
\ra\Hom_\Lam(M[d,b^+],M[d,b^+])\ra \C\ra 0.
\end{multline*}
Applying Proposition~\ref{prp:HomDim}
to the second and third non-trivial terms
yields the required equality,
if we take into account \eqref{eq:MutAdd} and the fact that $r_{bb}=1$.
\cqfd

\begin{Prop}\label{propcalculhom}
Let $1\le b < d \le r$ be such that $i_b = i_d=i$.
Let $j,k\in I$ with $j\not = i$ and $k \not = i$.
There holds
\begin{align*}
&\dim\Hom_\Lam(M[d,b^+],M[d^-,b]) = (\dimv\, V_{d},\, \dimv\, M[d^-,b]) - (\dimv\, M[d^-,b])_i,\\[2mm]
&\dim\Hom_\Lam(M[d,b],M[d^-,b^+])= (\dimv\, V_{d},\, \dimv\, M[d^-,b^+]) -(\dimv\, M[d^-,b^+])_i,\\[2mm]
&\dim\Hom_\Lam(M[d^-(j),(b^-(j))^+],\,M[d^-(k),(b^-(k))^+])= (\dimv\, V_{d^-(j)},\, \dimv\, M[d^-(k),(b^-(k))^+])\\[2mm] 
&\hskip 7.8cm -\ (\dimv\, M[d^-(k),(b^-(k))^+])_j.
\end{align*}
\end{Prop}

\proof
By Lemma~\ref{lem8.6},
\[
\dim\Hom_\Lam(M[d,b^+],M[d^-,b]) =
\dim\Hom_\Lam(M[d,b],M[d^-,b]).
\]
Moreover, Proposition~\ref{prp:HomDim} shows that
\[
\dim\Hom_\Lam(M[d,b],M[d^-,b])
=
\dim\Hom_\Lam(V_{d},M[d^-,b]).
\]
Now, it follows from (\ref{CBform}) that 
\[
\dim\Hom(V_{d},M[d^-,b])= (\dimv\,V_{d},\,\dimv\, M[d^-,b]) - (\dimv\, M[d^-,b])_i.
\]
Indeed,  
$V_{d}$ and $M[d^-,b]$ belong to the subcategory $\CC_{s_{i_{d}}\cdots s_{i_1}}$,
and $V_{d}$ is projective-injective in this category, with socle $S_i$. Therefore 
\[
\dim\Hom(M[d^-,b],V_{d}) =(\dimv\, M[d^-,b])_i \quad \text{ and } \quad 
\Ext^1_\Lam(M[d^-,b],V_{d})=0.
\]
This proves the first equation.
For the remaining two equations we note that again, by Proposition~\ref{prp:HomDim},
we can replace in the left hand side 
$M[d,b]$ by $V_{d}$, and $M[d^-(j),(b^-(j))^+]$ by $V_{d^-(j)}$.
The claim follows.
\cqfd

\section{The quantum cluster algebra associated with $\CC_w$}

\subsection{The cluster algebra $\AA(\CC_w)$}

Following \cite{GLSUni1}, \cite{BIRS}, \cite{GLS}, we can associate with $\CC_w$
a (classical \ie not quantum) cluster algebra. This is given by the initial seed
\begin{equation}
\Sigma_{V_\ii} := ((x_1,\ldots,x_r), \tB_{V_\ii}).
\end{equation}
Although this seed depends on the choice of a reduced expression
$\ii$ for $w$, one can show that any two matrices $\tB_{V_\ii}$ and $\tB_{V_\jj}$
are connected by a sequence of mutations. Therefore this cluster algebra 
is independent of this choice, and we denote it by $\AA(\CC_w)$.
Moreover, every seed of $\AA(\CC_w)$ is of the form
\[
\Sigma_T = ((x_{T_1},\ldots,x_{T_r}),\tB_T),
\]
for a unique $\CC_w$-maximal rigid module $T=T_1\oplus\cdots \oplus T_r$, and some Laurent polynomials
$x_{T_1},\ldots,x_{T_r}$ in the variables $x_1 = x_{V_{\ii,1}},\ldots,x_r = x_{V_{\ii,r}}$. 
These modules $T$ are those which
can be reached from $V_\ii$ using a sequence of mutations, and we called them \emph{reachable}. 
(It is still an open problem whether every $\CC_w$-maximal rigid module is reachable or not.) 
If $\jj$ is another reduced expression for $w$, it is known that $V_\jj$ is reachable 
from $V_\ii$ \cite{BIRS}. Therefore, the collection of reachable $\CC_w$-maximal rigid modules does
not depend on the choice of $\ii$.

It was shown in \cite{GLSUni1,GLS} that there is a natural isomorphism from
$\AA(\CC_w)$ to the coordinate ring $\C[N(w)]$, mapping the cluster monomials
to a subset of Lusztig's dual semicanonical basis of $\C[N(w)]$.

\subsection{The matrix $L_T$}

Let $T=T_1\oplus\cdots\oplus T_r$ be a $\CC_w$-maximal rigid module  
as in \S\ref{subsecGammaT}. 
Let $L_T=(\la_{ij})$ be the $r\times r$-matrix with entries
\begin{equation}\label{def_lambda}
\la_{ij} := \dim\Hom_\L(T_i,T_j) - \dim\Hom_\L(T_j,T_i),\qquad (1\le i, j\le r). 
\end{equation}
Note that $L_T$ is skew-symmetric. 
From now on we will use the following convenient shorthand notation.
Given two $\L$-modules $X$ and $Y$, we will write 
\begin{equation}
[X,Y] := \dim\Hom_\L(X,Y), \qquad [X,Y]^1 := \dim\Ext^1_\L(X,Y).
\end{equation} 
Thus, we shall write $\la_{ij} = [T_i,T_j]-[T_j,T_i]$.
\begin{Prop}\label{Tcomp}
The pair $(L_T,\tB_T)$ is compatible. 
\end{Prop}
\proof
For $1\le j\le r-n$, and $1\le i,\ k\le r$, by Theorem~\ref{th_mutation} (c) we have:
\[
\sum_{k=1}^r b_{kj}\la_{ki} = [T_i,T'_j]+[T''_j,T_i]-[T'_j,T_i]-[T_i,T''_j]. 
\]
Let us first assume that $i\not = j$. 
Applying the functor $\Hom_\L(T_i,-)$ to the short exact sequence 
\[
0 \to T_j \to T'_j \to T_j^* \to 0,
\]
and taking into account that $[T_i,T_j]^1=0$, we get a short exact sequence
\[
0 \to \Hom_\L(T_i,T_j) \to \Hom_\L(T_i,T'_j) \to \Hom_\L(T_i,T_j^*) \to 0,
\]
therefore
\begin{equation}
[T_i,T'_j] = [T_i,T_j^*] + [T_i,T_j]. 
\end{equation}
Similarly, applying the functor $\Hom_\L(-,T_i)$ to the same short exact sequence 
and taking into account that $[T_j^*,T_i]^1=0$, we get a short exact sequence
\[
0 \to \Hom_\L(T^*_j,T_i) \to \Hom_\L(T'_j,T_i) \to \Hom_\L(T_j,T_i) \to 0,
\]
therefore
\begin{equation}
[T'_j,T_i] = [T^*_j,T_i] + [T_j,T_i]. 
\end{equation}
It follows that
\begin{equation}\label{1steq}
[T_i,T'_j] - [T'_j,T_i] = [T_i,T_j^*] - [T^*_j,T_i] + \la_{ij}. 
\end{equation}
Arguing similarly with the short exact sequence
\[
0 \to T^*_j \to T''_j \to T_j \to 0,
\]
we obtain
\begin{equation}\label{2ndeq}
[T_j'',T_i] - [T_i,T_j''] = [T_j^*,T_i] - [T_i,T^*_j] + \la_{ji}. 
\end{equation}
Hence
$\sum_{k=1}^r b_{kj}\la_{ki} = \la_{ij} + \la_{ji} = 0$. 

Assume now that $i=j$. Using that $[T_j,T_j^*]^1=1$ and $[T_j,T_j]^1=[T_j,T_j']^1=[T_j,T_j'']^1=0$,
and arguing as above, we easily obtain the relations
\begin{align}
[T_j,T'_j] &= [T_j,T_j^*] + [T_j,T_j], \label{eq9.7}\\
[T'_j,T_j] &= [T^*_j,T_j] + [T_j,T_j]-1,\label{eq9.8}\\
[T_j'',T_j] &= [T_j^*, T_j] + [T_j,T_j],\label{eq9.9}\\
[T_j,T_j''] &= [T_j,T_j^*] + [T_j,T_j] -1\label{eq9.10}.
\end{align}
It follows that 
$\sum_{k=1}^r b_{kj}\la_{kj} = 2.$ 
Thus, in general, 
\[
\sum_{k=1}^r b_{kj}\la_{ki} = 2\de_{ij}, 
\]
which proves the proposition.
\cqfd

Let $k\le r-n$, so that the mutation $\mu_k(T)$ is well-defined.
\begin{Prop}\label{LTmut}
We have $\mu_k(L_T) = L_{\mu_k(T)}$. 
\end{Prop}
\proof
Put $\mu_k(L_T)=(\la'_{ij})$.
By definition, $\la'_{ij}=\la_{ij}$ if $i$ or $j$ is different from $k$.
Similarly, since $T$ and $\mu_k(T)$ differ only by the replacement of $T_k$
by $T_k^*$, all entries in $L_{\mu_k(T)}$ not situated on row $k$ or column $k$
are equal to the corresponding entries of $L_T$. 

If $i=k$, we have by definition of $\mu_k(L_T)$ and by Theorem~\ref{th_mutation}~(c) 
\[
\la'_{kj} = \sum_{s=1}^r e_{sk}\la_{sj} =  [T'_k,T_j] - [T_j,T'_k] -[T_k,T_j]+[T_j,T_k] - \la_{kj}. 
\]
Thus, by (\ref{1steq}), we get 
\[
\la'_{kj} = [T_k^*,T_j] - [T_j, T_k^*], 
\]
as claimed. The case $j=k$ follows by skew-symmetry. \cqfd 

Let $H_\ii=(h_{kl})$ be the $r\times r$-matrix with entries
\begin{equation}
h_{kl} =
\begin{cases}
1& \text{ if } l = k^{(-m)} \text{ for some } m\ge 0,\\
0& \text{ otherwise.}
\end{cases}
\end{equation}

\begin{Prop}\label{propLV}
The matrix $L_{V_\ii}$ has the following explicit expression
\[
L_{V_\ii} = H_\ii (R_\ii - R_\ii^T) H_\ii^T. 
\]
\end{Prop}
\proof This follows immediately from Proposition~\ref{prp:HomDim}, if we 
note that the $k$th row of $H_\ii$ is equal to the $M$-dimension vector $\mm(V_{\ii,k})$.  \cqfd

\subsection{The quantum cluster algebra $\AA_{q^{1/2}}(\CC_w)$}

Proposition~\ref{Tcomp} and Proposition~\ref{LTmut} show that the family 
$(L_T,\tB_T)$, where $T$ ranges over all $\CC_w$-maximal rigid modules reachable
from $V_\ii$, gives rise to a quantum analogue of the cluster algebra $\AA(\CC_w)$. 
We shall denote it by $\AA_{q^{1/2}}(\CC_w)$.
For every reduced expression $\ii$ of $w$, an explicit initial quantum seed is given by 
\[
\SC_{V_\ii}:=((X_{V_{\ii,1}},\ldots,X_{V_{\ii,r}}),\ L_{V_\ii},\ \tB_{V_\ii}),
\]
where the matrices $L_{V_\ii}$ and $\tB_{V_\ii}$ have been computed in 
Proposition~\ref{propLV}, and in Section~\ref{subsecGammaT}.

The quantum seed corresponding to a reachable $\CC_w$-maximal rigid module 
$T=T_1\oplus\cdots\oplus T_r$
will be denoted by
\begin{equation}
\SC_T = ((X_{T_1},\ldots,X_{T_r}),\ L_T,\ \tB_T).  
\end{equation}
For every $\aa = (a_1,\ldots, a_r)\in \N^r$, we have a rigid module 
$T^\aa := T_1^{a_1}\oplus \cdots \oplus T_r^{a_r}$ in the additive closure
$\add(T)$ of $T$.
Following (\ref{eqX^a}) and writing $L_T = (\la_{ij})$, we put
\begin{equation}
X_{T^\aa} = X^\aa:=q^{\frac{1}{2}\sum_{i>j}a_ia_j\la_{ij}}X_{T_1}^{a_1}\cdots X_{T_r}^{a_r}.
\end{equation}
Thus, denoting by $\R_w$ the set of rigid modules 
$R$ in the additive closure of some reachable $\CC_w$-maximal rigid module, 
we obtain a canonical labelling 
\begin{equation}
 X_R, \qquad (R\in \R_w).
\end{equation}
of the quantum cluster monomials of $\AA_{q^{1/2}}(\CC_w)$.

\subsection{The elements $Y_R$}\label{eltsY}
It will be convenient in our setting to proceed to a slight rescaling
of the elements $X_R$. For $R\in\R_w$, we define
\begin{equation}\label{XM}
Y_R := q^{-[R,R]/2}\, X_R.
\end{equation}
In particular, writing $R = T^\aa$ as above, an easy calculation gives
\begin{equation}
 Y_{R} = q^{-\a(R)}
Y^{a_1}_{T_1}\cdots Y^{a_r}_{T_r},
\end{equation}
where
\begin{equation}
\a(R) := \sum_{i<j}a_ia_j[T_i,T_j] + \sum_i{a_i\choose2}[T_i,T_i].
\end{equation}
Note that $\a(R)$ is an integer (not a half-integer), so that $ Y_R$
belongs to $\Z[q^{\pm1}][Y_{T_1},\cdots,Y_{T_r}]$.
Moreover, we have the following easy lemma.
\begin{Lem}\label{YMN}
Let $T$ be a reachable $\CC_w$-maximal rigid module.
For any $R, S$ in $\add(T)$ we have
\[
Y_R Y_S = q^{[R,S]} Y_{R\oplus S}.
\]
\end{Lem}
\proof
Write $R = T_1^{a_1} \oplus \cdots \oplus T_r^{a_r}$
and $S = T_1^{b_1} \oplus \cdots \oplus T_r^{b_r}$.
We have
\[
Y_R Y_S = q^{-\a(R)-\a(S)} Y^{a_1}_{T_1}\cdots 
Y^{a_r}_{T_r} Y^{b_1}_{T_1}\cdots Y^{b_r}_{T_r}
= q^{-\a(R)-\a(S)+\sum_{i>j}a_ib_j\la_{ij}}
Y^{a_1+b_1}_{T_1}\cdots Y^{a_r+b_r}_{T_r}.
\]
On the other hand, using the obvious identity
\[
{a+b\choose 2} = {a\choose 2} + {b\choose 2} + ab,
\]
we see that
\[
\a(R\oplus S) = \a(R) + \a(S) +
\sum_{i<j} (a_ib_j + a_jb_i) [T_i,T_j] +
\sum_i a_ib_i[T_i,T_i].
\]
Hence, taking into account (\ref{def_lambda}),
\begin{eqnarray*}
-\a(R)-\a(S)+\sum_{i>j}a_ib_j\la_{ij} &=&
-\a(R\oplus S) + \sum_i a_ib_i[T_i,T_i]
+\sum_{i>j} a_ib_j[T_i,T_j]
+\sum_{i<j} a_ib_j[T_i,T_j]\\
&=&-\a(R\oplus S) + [R,S],
\end{eqnarray*}
and the result follows.
\cqfd 

Note that because of (\ref{CBform}) and the fact that every 
$R\in\add(T)$ is rigid, we have 
\begin{equation}
[R,R]=(\dimv\, R,\ \dimv\, R)/2,
\end{equation}
hence the exponent of $q$ in (\ref{XM}) depends only on
the dimension vector $\dimv\,R$ of $R$.

\subsection{Quantum mutations}\label{mutY}
We now rewrite formulas (\ref{qmut1}) (\ref{qmut2}) for quantum mutations in $\AA_{q^{1/2}}(\CC_w)$,
using the rescaled quantum cluster monomials $Y_R$.

Let $T$ be a $\CC_w$-maximal rigid module, and let $T_k$ be a non-projective
indecomposable direct summand of $T$. Let $\mu_k(T) = (T/T_k)\oplus T_k^*$ be
the mutation of $T$ in direction $k$.
By Theorem~\ref{th_mutation}, we have two short exact sequences
\begin{equation}
0\to T_k\to T'_k \to T_k^* \to 0, \qquad  0\to T_k^*\to T''_k \to T_k \to 0,
\end{equation}
with $T'_k, T''_k \in \add(T/T_k)$.
The quantum exchange relation between the quantum cluster variables $Y_{T_k}$
and $Y_{T_k^*}$ can be written as follows:
\begin{Prop}\label{Propqmut}
With the above notation, we have 
\[
Y_{T_k^*}Y_{T_k} = q^{[T_k^*,T_k]}(q^{-1}Y_{T'_k} + Y_{T''_k}). 
\]
\end{Prop}
\proof
We have 
\begin{eqnarray*}
Y_{T_k^*}Y_{T_k}&=& q^{-([T_k^*,T_k^*]+[T_k,T_k])/2}X_{T_k^*}X_{T_k}\\
&=&q^{-([T_k^*,T_k^*]+[T_k,T_k])/2}
\left(X^{\aa'}+ X^{\aa''}\right)X_{T_k},
\end{eqnarray*}
where, if we write 
\[
T'_k = \bigoplus_{i\not = k} T_i^{a'_i},\qquad
T''_k = \bigoplus_{i\not = k} T_i^{a''_i}, 
\]
the multi-exponents $\aa'$ and $\aa''$ are given by
\[
\aa' = (a'_1\ldots,a'_{k-1},-1,a'_{k+1},\ldots,a'_r),\quad
\aa'' = (a''_1\ldots,a''_{k-1},-1,a''_{k+1},\ldots,a''_r). 
\]
Using (\ref{qT}), (\ref{eqX^a}), and (\ref{def_lambda}), one obtains
easily that
\[
X^{\aa'}X_{T_k} = q^{([T'_k,T_k]-[T_k,T'_k])/2} X_{T'_k},\qquad
X^{\aa''}X_{T_k} = q^{([T''_k,T_k]-[T_k,T''_k])/2} X_{T''_k}. 
\]
Now, using twice (\ref{CBform}), we have
\begin{eqnarray}
[T'_k,T'_k]&=&(\dimv\, T'_k\,,\ \dimv\, T'_k)/2\\
&=&(\dimv\, T_k + \dimv\, T_k^*\,,\ \dimv\, T_k + \dimv\, T_k^*)/2\\
&=&[T_k,T_k]+[T_k^*,T_k^*]+(\dimv\, T_k\,,\ \dimv\, T_k^*)\\
&=&[T_k,T_k]+[T_k^*,T_k^*]+[T_k,T_k^*]+[T_k^*,T_k]-1,\label{relationhom}
\end{eqnarray}
and this is also equal to $[T''_k,T''_k]$.
Therefore, the exponent of $q^{1/2}$ in front of $Y_{T'_k}$ in the product
$Y_{T_k^*}Y_{T_k}$ is equal to
\begin{eqnarray*}
[T'_k,T_k]-[T_k,T'_k]+[T'_k,T'_k]-[T_k^*,T_k^*]-[T_k,T_k]
&=& 
[T_k,T_k^*]+[T_k^*,T_k]-1-[T_k,T'_k]+[T'_k,T_k]\\
&=& 2([T_k^*,T_k]-1),
\end{eqnarray*}
where the second equality follows from (\ref{eq9.7}) and (\ref{eq9.8}).
Similarly, using (\ref{eq9.9}) and (\ref{eq9.10}), we see that 
the exponent of $q^{1/2}$ in front of $Y_{T''_k}$ in the product
$Y_{T_k^*}Y_{T_k}$ is equal to $2[T^*_k,T_k]$, as claimed.
\cqfd

\subsection{The rescaled quantum cluster algebras $\AA_q(\CC_w)$ and $\AA_\A(\CC_w)$}

By definition, $\AA_{q^{1/2}}(\CC_w)$ is a $\Z[q^{\pm1/2}]$-algebra.
It follows from Lemma~\ref{YMN} and Proposition~\ref{Propqmut} that if we replace the
quantum cluster variables $X_{T_i}$ by their rescaled versions $Y_{T_i}$
we no longer need half-integral powers of $q$.
So we are led to introduce the rescaled quantum cluster algebra
$\AA_q(\CC_w)$.
This is defined as the $\Z[q^{\pm1}]$-subalgebra of $\AA_{q^{1/2}}(\CC_w)$
generated by the elements $Y_{T_i}$, where $T_i$ ranges over the
indecomposable direct summands of all reachable $\CC_w$-maximal rigid modules. 

It turns out that in the sequel, in order to have a coefficient ring which 
is a principal ideal domain, it will be convenient to slightly extend coefficients
from $\Z[q^{\pm1}]$ to $\A = \Q[q^{\pm1}]$. We will denote by
\[
\AA_\A(\CC_w) := \A \otimes_{\Z[q^{\pm1}]} \AA_q(\CC_w) 
\]
the corresponding quantum cluster algebra. Of course, the based quantum
tori of $\AA_q(\CC_w)$ (\resp  $\AA_\A(\CC_w)$) will be defined over
$\Z[q^{\pm1}]$ (\resp over $\A$).

\subsection{The involution $\si$} \label{sectinvolsigma}

We now define an involution of $\AA_\A(\CC_w)$, which will turn out to be related
to the involution $\si$ of $U_q(\n)$ defined in \S\ref{dualcan}. For this reason
we shall also denote it by $\si$. This involution is a twisted-bar involution,
as defined by Berenstein and Zelevinsky \cite[\S6]{BZ}.

We first define an additive group automorphism $\si$
of the $\A$-module $\A[Y_R\mid R\in \add(V_\ii)]$ by setting
\begin{eqnarray}
\si(f(q)Y_R)   &=& f(q^{-1})\, q^{[R,R] - \dim R} Y_R, \qquad (f \in \A,\ R \in \add(V_\ii))\label{barYM}.
\end{eqnarray}
Clearly, $\si$ is an involution.
\begin{Lem}\label{baraa}
$\si$ is a ring anti-automorphism.
\end{Lem}
\proof
We have
\begin{eqnarray*}
\si{(Y_R Y_S)}&=&q^{-[R,S]}\si(Y_{R\oplus S})\\
&=&q^{-[R,S]-\dim(R\oplus S)+[R\oplus S,\ R\oplus S]}Y_{R\oplus S}\\
&=&q^{-\dim R - \dim S + [R,R] +[S,S] + [S,R]}Y_{R\oplus S}\\
&=&q^{-\dim R - \dim S + [R,R] + [S,S]} Y_S Y_R\\
&=&\si(Y_S)\si(Y_R).
\end{eqnarray*}
\cqfd

The involution $\si$ extends to an anti-automorphism of the based quantum torus 
\begin{equation}
\V_\ii := \A[Y_{V_i}^{\pm 1}\mid 1\le i\le r]
\end{equation}
which we still denote by $\si$. 
Moreover, we have
\begin{Lem}\label{lemsigma}
For every indecomposable reachable rigid module $U$ we have
\[
\si(Y_U)= q^{[U,U] - \dim U} Y_U. 
\]
\end{Lem}
\proof
By induction on the length of the mutation sequence, we may assume
that the result holds for all indecomposable direct summands $T_j$
of a $\CC_w$-maximal rigid module $T$. We then have to check that
it also holds for $U = T_k^* = \mu_k(T_k)$.
By Proposition~\ref{Propqmut}, we have 
\[
q^{[T_k,T_k] - \dim T_k} Y_{T_k} \si(Y_{T_k^*}) =
q^{-[T_k^*,T_k]}\left(q^{[T'_k,T'_k]-\dim T'_k +1} Y_{T'_k} +
q^{[T''_k,T''_k]-\dim T''_k} Y_{T''_k}\right).
\]
Using that $\dim T'_k=\dim T''_k = \dim T_k + \dim T^*_k$, and 
(\ref{relationhom}), this becomes
\[
Y_{T_k} \si(Y_{T_k^*}) = q^{[T^*_k,T_k^*] - \dim T^*_k + [T_k, T_k^*]}(Y_{T'_k} + q^{-1} Y_{T''_k})
= q^{[T^*_k,T_k^*] - \dim T^*_k} Y_{T_k} Y_{T_k^*},
\]
where the last equality comes again from Proposition~\ref{Propqmut}, since quantum mutations
are involutive. It follows that $\si(Y_{T_k^*}) = q^{[T^*_k,T_k^*] - \dim T^*_k} Y_{T_k^*}$,
as claimed, and this proves the lemma.
\cqfd 

By Lemma~\ref{lemsigma}, $\si$ restricts to an automorphism of $\AA_\A(\CC_w)$, that we again
denote by $\si$. Moreover, for any quantum cluster monomial, that is for every element
$Y_U$ where $U$ is a reachable, non necessarily indecomposable, rigid module in $\CC_w$,
there holds
\begin{equation}\label{sigmaY}
\si(Y_U)= q^{[U,U] - \dim U} \, Y_U = q^{(\dimv\, U,\ \dimv\, U)/2 - \dim U} \, Y_U
= q^{N(\dimv\, U)} Y_U,
\end{equation}
where, for $\b\in Q_+$, $N(\b)$ is defined in (\ref{defN(b)}).

\section{Based quantum tori}

In this section, we fix again a reduced word $\ii:=(i_r,\ldots,i_1)$ for $w$,
and we set 
\begin{equation}
\la_k = s_{i_1}\cdots s_{i_k}(\vpi_{i_k}),\qquad (k\le r).
\end{equation} 
\subsection{Commutation relations}
The next lemma is a particular case of \cite[Theorem~10.1]{BZ}. We include 
a brief proof for the convenience of the reader.
\begin{Lem}\label{LemDe}
For $1\le k<l \le r$, we have in $A_q(\g)$:
\[
\De_{\vpi_{i_k},\,\la_k} \De_{\vpi_{i_l},\,\la_l} = 
q^{(\vpi_{i_k},\vpi_{i_l})-(\la_k,\la_l)}
\De_{\vpi_{i_l},\,\la_l} \De_{\vpi_{i_k},\,\la_k}. 
\]
\end{Lem}
\proof
Since $k<l$, we have $\la_l = s_{i_1}\cdots s_{i_k}(\nu)$, where $\nu = s_{i_{k+1}}\cdots s_{i_l}(\vpi_{i_l})$.
For $x\in U_q(\n)$ and $y\in U_q(\n_-)$, we have 
\[
x\cdot \De_{\vpi_{i_k},\vpi_{i_k}} = \eps(x) \De_{\vpi_{i_k},\vpi_{i_k}},
\qquad
\De_{\vpi_{i_l},\nu}\cdot y = \eps(y)\De_{\vpi_{i_l},\nu}. 
\]
Therefore, using again \cite[Lemma 10.2]{BZ} as in the proof of Lemma~\ref{commute},
we obtain that 
\[
\De_{\vpi_{i_k},\vpi_{i_k}} \De_{\vpi_{i_l},\nu} =
q^{(\vpi_{i_k},\vpi_{i_l}) -(\vpi_{i_k},\nu)}
\De_{\vpi_{i_l},\nu} \De_{\vpi_{i_k},\vpi_{i_k}}. 
\]
Now, using $l-k$ times Lemma~\ref{lem:QMinIdRed}, we deduce from this equality that
\[
\De_{\vpi_{i_k},\, s_{i_1}\cdots s_{i_k}(\vpi_{i_k})}\, \De_{\vpi_{i_l},\, s_{i_1}\cdots s_{i_k}(\nu)} =
q^{(\vpi_{i_k},\vpi_{i_l}) -(\vpi_{i_k},\nu)}
\De_{\vpi_{i_l},\, s_{i_1}\cdots s_{i_k}(\nu)}\, \De_{\vpi_{i_k},\, s_{i_1}\cdots s_{i_k}(\vpi_{i_k})},  
\]
and taking into account that 
$(\vpi_{i_k},\,\nu) = (s_{i_1}\cdots s_{i_k}(\vpi_{i_k}),\,s_{i_1}\cdots s_{i_k}(\nu))= (\la_k,\,\la_l)$,
we get the claimed equality.
\cqfd

\begin{Lem}\label{qtor1}
For $1\le k<l \le r$, we have in $A_q(\n)$:
\[
D_{\vpi_{i_k},\,\la_k} D_{\vpi_{i_l},\,\la_l} = 
q^{(\vpi_{i_k}-\la_k,\ \vpi_{i_l}+\la_l)}
D_{\vpi_{i_l},\,\la_l} D_{\vpi_{i_k},\,\la_k}. 
\]
\end{Lem}
\proof This follows from Lemma~\ref{LemDe} using the same type of calculations as in the proof of Lemma~\ref{qcommute2}.
We leave the easy verification to the reader.
\cqfd

\begin{Lem}\label{qtor3}
For $1\le k<l \le r$, we have:
\[
(\vpi_{i_k}-\la_k,\ \vpi_{i_l}+\la_l) = [V_{\ii,k}, V_{\ii,l}] - [V_{\ii,l}, V_{\ii,k}]. 
\]
\end{Lem}
\proof
We have 
\[
(\vpi_{i_k}-\la_k,\ \vpi_{i_l}+\la_l) = (\dimv\,V_{\ii,k},\,2\vpi_{i_l} -\dimv\,V_{\ii,l}). 
\]
Since $[V_{\ii,k},\,V_{\ii,l}]^1 = 0$, using (\ref{CBform}) we have
\[
(\dimv\,V_{\ii,k},\,\dimv\,V_{\ii,l}) = [V_{\ii,k}, V_{\ii,l}] + [V_{\ii,l}, V_{\ii,k}].
\]
Hence,
\[
(\vpi_{i_k}-\la_k,\ \vpi_{i_l}+\la_l) 
= 2(\dimv\,V_{\ii,k},\,\vpi_{i_l}) - [V_{\ii,k}, V_{\ii,l}] - [V_{\ii,l}, V_{\ii,k}], 
\]
and we are reduced to show that
\begin{equation}\label{eqinterest}
(\dimv\,V_{\ii,k},\, \vpi_{i_l}) = [V_{\ii,k}, V_{\ii,l}],\qquad (k<l).
\end{equation}
By Proposition~\ref{propLV}, the right-hand side is equal to
\[
\sum_{m\le 0}\left( \sum\limits_{\substack{s<k^{(m)}\\ i_s = i_l}} (\b_{k^{(m)}},\,\b_s) + \de_{i_k i_l}\right). 
\]
Since the left-hand side is equal to 
$\sum_{m\le 0} (\b_{k^{(m)}}, \vpi_{i_l})$, it is enough to show that for every $m\ge 0$,
\begin{equation}\label{eqabove}
(\b_{k^{(m)}}, \vpi_{i_l}) = \sum\limits_{\substack{s<k^{(m)}\\ i_s = i_l}} (\b_{k^{(m)}},\,\b_s) + \de_{i_k i_l}. 
\end{equation}
Let $t:=\max\{s < k^{(m)}\mid i_s = i_l\}$. 
Eq. (\ref{eqabove}) can be rewritten
\begin{equation}\label{eqjustabove}
(\b_{k^{(m)}},\,\la_t) = \de_{i_k i_l}. 
\end{equation}
But the left-hand side of (\ref{eqjustabove}) is equal to
\[
(s_{i_1}\cdots s_{i_{k^{(m)}-1}}(\a_{i_k}),\ s_{i_1}\cdots s_{i_t}(\vpi_{i_l})) 
=
(s_{i_1}\cdots s_{i_{k^{(m)}-1}}(\a_{i_k}),\ s_{i_1}\cdots s_{i_{k^{(m)}-1}}(\vpi_{i_l}))
=
(\a_{i_k},\,\vpi_{i_l}) 
= \de_{i_k i_l},
\] 
so (\ref{eqjustabove}) holds, and this proves (\ref{eqinterest}).
\cqfd

\subsection{An isomorphism of based quantum tori}
By Proposition~\ref{flagFq}, the quantum flag minors $D_{\vpi_{i_k},\la_k}=D(0,k)$
belong to $F_q(\n(w))$.
Let $\T_\ii$ be the $\A$-subalgebra of $F_q(\n(w))$ generated by the 
$D_{\vpi_{i_k},\la_k}^{\pm 1}\ (1\le k\le r)$.

\begin{Lem}\label{qtor2} 
The algebra $\T_\ii$ is a based quantum torus over $\A$.
\end{Lem}
\proof
By Lemma~\ref{qtor1}, the generators $D_{\vpi_{i_k},\la_k}$ pairwise $q$-commute.
So we only have to show that the monomials 
\[
D^\aa := \prod_k^{\longrightarrow} D_{\vpi_{i_k},\la_k}^{a_k},\qquad (\aa=(a_1,\ldots,a_r)\in \N^r), 
\]
are linearly independent over $\A$.
Suppose that we have a non-trivial relation
\[
 \sum_{\aa} \ga_\aa(q) D^\aa = 0,
\]
for some $\ga_\aa(q)\in\A$. Dividing this equation (if necessary) by the largest power of $q-1$ which
divides all the coefficients $\ga_\aa(q)$, we may assume that at least one of these
coefficients is not divisible by $q-1$. Using \ref{specq1}, we see that by specializing
this identity at $q=1$ we get a non-trivial $\C$-linear relation between monomials in the
corresponding classical flag minors of $\C[N]$. But all these monomials belong to the dual semicanonical
basis of $\C[N]$ (see \cite[Corollary 13.3]{GLS}), hence they are linearly independent, a contradiction.
\cqfd 

We note that Lemma~\ref{qtor2} also follows from \cite[Theorem 6.20]{Ki}.

\begin{Prop}\label{isomqtori}
The assignment $Y_{V_{\ii,k}}\mapsto  D_{\vpi_{i_k},\,\la_k}\ (1\le k\le r)$ extends to an
algebra isomorphism from $\V_\ii$ to the based quantum torus $\T_\ii$.
\end{Prop}

\proof 
By definition of the cluster algebra $\AA_\A(\CC_w)$, the elements 
$Y_{V_{\ii,k}}^{\pm 1}$ generate a based quantum torus over $\A$,
with $q$-commutation relations
\[
Y_{V_{\ii,k}}Y_{V_{\ii,l}} = q^{[V_{\ii,k},V_{\ii,l}] - [V_{\ii,l},V_{\ii,k}]}
Y_{V_{\ii,l}}Y_{V_{\ii,k}}, \qquad (1\le k<l \le r). 
\]
The proposition then follows immediately from Lemma~\ref{qtor1} and Lemma~\ref{qtor3}.
\cqfd

\section{Cluster structures on quantum coordinate rings}

\subsection{Cluster structures on quantum coordinate rings of unipotent subgroups}

Consider the following diagram of homomorphisms of $\A$-algebras:
\begin{equation}
\AA_\A(\CC_w) \longrightarrow \V_\ii \stackrel{\sim}{\longrightarrow}
\T_\ii \longrightarrow F_q(\n(w)). 
\end{equation}
Here, the first arrow denotes the natural embedding given by 
the (quantum) Laurent phenomenon \cite{BZ}, the second arrow is 
the isomorphism of Proposition~\ref{isomqtori}, and the third 
arrow is the natural embedding.
The composition $\k \colon \AA_\A(\CC_w) \longrightarrow F_q(\n(w))$ of these
maps is therefore injective.
Recall the notation $A_\A(\n(w))$ of \S\ref{specialq1}.

\begin{Prop}\label{surj}
The image $\k\left(\AA_\A(\CC_w)\right)$ contains $A_\A(\n(w))$. 
\end{Prop}

\proof
Since $\k$ is an algebra map, it is enough to show that its image
contains the dual PBW-generators $E^*(\b_k)\ (1\le k \le r)$.

It was shown in \cite[\S13.1]{GLS} that there is an explicit sequence
of mutations of $\CC_w$-maximal rigid modules starting from $V_\ii$ and
ending in $T_\ii$. Each step of this sequence consists of the mutation
of a module $M[d^-,b]$ into a module $M[d,b^+]$, for some 
$1\le b < d \le r$ with $i_b = i_d = i$.
The corresponding pair of short exact sequences is 
\[
0 \to M[d^-,b] \to M[d^-,b^+] \oplus M[d,b] \to M[d,b^+] \to 0, 
\]
\[
0 \to M[d,b^+] \to \bigoplus_{j\not = i} M[d^-(j),(b^-(j))^+]^{-a_{ij}} \to M[d^-,b] \to 0.
\]
Write 
\[
T'_{bd} := M[d^-,b^+] \oplus M[d,b]
\] 
and 
\[
T''_{bd} := \bigoplus_{j\not = i} M[d^-(j),(b^-(j))^+]^{-a_{ij}}.
\]
Then, by Proposition~\ref{Propqmut}, in $\AA_\A(\CC_w)$ we have the 
corresponding mutation relation:
\begin{equation}\label{elementary_mut}
Y_{M[d,b^+]} Y_{M[d^-,b]} = q^{[M[d,b^+],\,M[d^-,b]]}\left(q^{-1} Y_{T'_{bd}} + Y_{T''_{bd}}\right). 
\end{equation}
Moreover,
\[
Y_{T'_{bd}} = q^{-\a(T'_{bd})} Y_{M[d^-,b^+]} Y_{M[d,b]},
\qquad
Y_{T''_{bd}} = q^{-\a(T''_{bd})} \prod_{j\not = i}^{\longrightarrow} Y_{M[d^-(j),(b^-(j))^+]}^{-a_{ij}}
\]
where
$\a(T'_{bd}) = [M[d^-,b^+],\, M[d,b]]$, and
\begin{align*}
\a(T''_{bd}) = & \sum_{j<k}a_{ij}a_{ik} [M[d^-(j),(b^-(j))^+],M[d^-(k),(b^-(k))^+]] \\
&+\ \sum_{j\not = i}{-a_{ij}\choose2} [M[d^-(j),(b^-(j))^+],M[d^-(j),(b^-(j))^+]].
\end{align*}
Note that for $1\le k \le l \le r$, and $i_k = i_l =j$, one has
\begin{equation}
\dimv\,M[l,k] = \mu(k^-,j) - \mu(l,j).
\end{equation}
Therefore, using Proposition~\ref{propcalculhom} and 
the notation of Proposition~\ref{Tsystem}, we see that
\begin{align*}
[M[d,b^+],\,M[d^-,b]] &= (\dimv\, V_{d},\, \dimv\, M[d^-,b]) - (\dimv\, M[d^-,b])_i\\[2mm]
&= (\vpi_i- \mu(d,i),\, \mu(b^-,i)-\mu(d^-,i)) - (\vpi_i,\, \mu(b^-,i)-\mu(d^-,i))\\[2mm]
&= - (\mu(d,i),\, \mu(b^-,i)-\mu(d^-,i))\\[2mm]
&= - A.
\end{align*}
Similarly we obtain that
\[
\a(T'_{bd}) = -B,\qquad  \a(T''_{bd}) = -C.
\]
Therefore (\ref{elementary_mut}) can be rewritten as
\begin{equation}\label{eqYsystem}
q^A Y_{M[d,b^+]} Y_{M[d^-,b]} = q^{-1+B} Y_{M[d^-,b^+]}Y_{M[d,b]} + q^{C}\prod_{j\not = i}^{\longrightarrow} Y_{M[d^-(j),(b^-(j))^+]}^{-a_{ij}}.
\end{equation}
We observe that this has exactly the same form as (\ref{eqTsystem}).
By definition of $\k$ one has
\[
\k(Y_{V_{\ii,k}}) = \k(Y_{M[k,k_{\min}]}) = D(0,k).
\]
Hence all the initial variables of the systems of equations (\ref{eqYsystem}) 
and (\ref{eqTsystem}) are matched under the algebra homomorphism $\k$.
Therefore, by induction, it follows that $\k(Y_{M[d,b^+]}) = D(b,d)$ for
every $b<d$. In particular, $\k(Y_{M[k,k]}) = D(k^-,k)=E^*(\b_k)$, 
by Proposition~\ref{PBWminors}.
\cqfd

\begin{Prop}\label{rank}
The algebra $\AA_\A(\CC_w)$ is a $Q_+$-graded free $\A$-module,
with homogeneous components of finite rank. 
Moreover, 
\[
\rk \left(\AA_\A(\CC_w)_\a\right) = \rk \left(A_\A(\n(w))_\a\right), \qquad (\a\in Q_+).
\]
\end{Prop}
\proof
The $\A$-module $\V_\ii$ is free, hence
its submodule $\AA_\A(\CC_w)$ is projective, and therefore free 
since $\A$ is a principal ideal domain.
It has a natural $Q_+$-grading given by
\[
\deg Y_R := \dimv\, R 
\]
for every indecomposable reachable rigid module $R$ in $\CC_w$. 
The rank of $\AA_\A(\CC_w)_\a$ is equal to the dimension of the
$\C$-vector space $\C\otimes_\A \AA_\A(\CC_w)_\a$, which is  
equal to the dimension of the corresponding 
homogeneous component of the (classical) cluster algebra
$\C\otimes_\Z \AA(\CC_w)$.
Now, by \cite{GLS}, this is equal to the dimension of 
$\C[N(w)]_\a$, that is, by \S\ref{specialq1}, to the rank of
$A_\A(\n(w))_\a$.
\cqfd 

Let $\AA_{\Q(q)}(\CC_w) := \Q(q)\otimes_\A \AA_{\A}(\CC_w)$. 
The $\A$-algebra homomorphism $\k$ naturally extends to a $\Q(q)$-algebra
homomorphism from $\AA_{\Q(q)}(\CC_w)$ to $F_q(\n(w))$,
which we continue to denote by~$\k$.
We can now prove our main result:

\begin{Thm}\label{main}
$\k$ is an isomorphism from the quantum cluster algebra $\AA_{\Q(q)}(\CC_w)$ to 
the quantum coordinate ring $A_q(\n(w))$. 
\end{Thm}

\proof
By construction, $\k$ is injective. By Proposition~\ref{surj},
the image of $\k$ contains $A_q(\n(w))$. 
Finally, since $\k$ preserves the $Q_+$-gradings, and since the
homogeneous components of $A_{\Q(q)}(\CC_w)$ and 
$A_q(\n(w))$ have the same dimensions, by Proposition~\ref{rank},
we see that 
\[
\k\left(\AA_{\Q(q)}(\CC_w)\right) = A_q(\n(w)).
\]
\cqfd

The following Corollary proves the claim made at the end of \S\ref{skewfield}. 
\begin{Cor}\label{minorpol}
All quantum minors $D(b,d)\ (1\le b<d\le r)$ belong to $A_q(\n(w))$,
and therefore are polynomials in the dual PBW generators $E^*(\b_k)=D(k^-,k)=Y_{M_k}$.
\end{Cor}
\proof As shown in the proof of Proposition~\ref{surj}, $D(b,d) = \k(Y_{M[d,b^+]})$,
and so belongs to $\k\left(\AA_{\Q(q)}(\CC_w)\right) = A_q(\n(w))$, by 
Theorem~\ref{main}.
\cqfd

The following Corollary is a $q$-analogue of \cite[Theorem 3.2 (i) (ii)]{GLS}.
\begin{Cor}
\begin{itemize}
\item[(a)] $\AA_{\Q(q)}(\CC_w)$ is an iterated skew polynomial ring. 
\item[(b)] The set $\{Y_M(\aa):= Y_{M_1}^{a_1}\cdots Y_{M_r}^{a_r} \mid \aa =(a_1,\ldots, a_r)\in \N^r\}$
is a $\Q(q)$-linear basis of $\AA_{\Q(q)}(\CC_w)$.
\end{itemize}
\end{Cor}
\proof
This follows via the isomorphism $\k$ from known properties of $A_q(\n(w))$.
\cqfd

\begin{Cor}
Let $\g$ be a simple Lie algebra of type $A_n, D_n, E_n$\,, and let $\n$ be 
a maximal nil\-po\-tent subalgebra of $\g$. Let $\L$ be the corresponding 
preprojective algebra of type $A_n, D_n, E_n$. Then
\begin{itemize}
\item[(a)] The quantum cluster algebra $\AA_{\Q(q)}(\md\L)$
is isomorphic to the quantum enveloping algebra $U_q(\n)$.
\item[(b)] In this isomorphism, the Chevalley generators $e_i\ (1\le i\le n)$
are identified with the quantum cluster variables $Y_{S_i}$ attached to
the simple $\L$-modules $S_i$. 
\end{itemize}
\end{Cor}
\proof
(a) Take $w$ to be the longest element $w_0$ of the Weyl group. Hence 
$\n(w_0) = \n$, so
\[
A_q(\n(w_0)) = A_q(\n) \cong U_q(\n), 
\]
by Proposition~\ref{isomPhi}. On the other hand $\CC_{w_0} = \md(\L)$,
and this proves (a). Property (b) is then obvious. 
\cqfd

We believe that Theorem~\ref{main} can be strengthened as follows.

\begin{Conj}\label{integral}
The map $\k$ restricts to an isomorphism from the integral form $\AA_{\A}(\CC_w)$ of the
quantum cluster algebra to the integral 
form $A_\A(\n(w))$ of the quantum coordinate ring. 
\end{Conj}

\subsection{Example}\label{exampleA3}

We illustrate our arguments on a simple example. We take $\g$ of type $A_3$
and $w=w_0$. We choose the reduced word $\ii = (i_6,i_5,i_4,i_3,i_2,i_1) = (1,2,3,1,2,1)$.
Therefore
\[
\b_1=\a_1,\ \
\b_2 = \a_1+\a_2,\ \
\b_3 = \a_2,\ \
\b_4 = \a_1+\a_2+\a_3,\ \
\b_5 =\a_2+\a_3,\ \
\b_6 = \a_3.
\]
Using the convention of \cite[\S2.4]{GLS} for visualizing $\L$-modules,
we can represent the summands of $V_\ii$ as follows: 
\begin{align*}
V_1 &= {\bsm1\esm} &
V_2 &= {\bsm1\\&2\esm} &
V_3 &= {\bsm&2\\1\esm} \\[2mm]
V_4 &= {\bsm1\\&2\\&&3\esm} &
V_5 &= {\bsm&2\\1&&3\\&2\esm} &
V_6 &= {\bsm&&3\\&2\\1\esm}
\end{align*}
and the summands of $M_\ii$ as follows:
\begin{align*}
M_1 &= {\bsm1\esm} &
M_2 &= {\bsm1\\&2\esm} &
M_3 &= {\bsm2\esm} \\[2mm]
M_4 &= {\bsm1\\&2\\&&3\esm} &
M_5 &= {\bsm2\\&3\esm} &
M_6 &= {\bsm3\esm}
\end{align*}
The sequence of mutations of \cite[\S13]{GLS} consists here in 4 mutations.

\medskip\noindent
\emph{Mutation 1:} One mutates at $V_1$ in the maximal rigid module $V_\ii$.
One has 
\[
\mu_1(V_1) = \mu_1(M[3^-,1]) = M[3,1^+] = M[3,3]=M_3.
\]
This gives rise to the two short exact sequences (which can be read from the
graph $\G_\ii$):
\[
0\to V_1 \to V_3 \to M_3 \to 0,\qquad 0\to M_3 \to V_2 \to V_1\to 0. 
\]
By Proposition~\ref{mutY}, we thus have
\[
 Y_{M_3}Y_{V_1} = q^{[M_3,V_1]}\left(q^{-1} Y_{V_3} + Y_{V_2}\right) = q^{-1} Y_{V_3} + Y_{V_2}, 
\]
since $[M_3,V_1] = 0$. Using the notation $M[l,k]$, this can be rewritten
\begin{equation}\label{eqY1}
Y_{M[3,3]}Y_{M[1,1]} = q^{-1}Y_{M[3,1]} + Y_{M[2,2]}. 
\end{equation}

\medskip\noindent
\emph{Mutation 2:} One mutates at $V_3$ in the maximal rigid module $\mu_1(V_\ii)$.
One has 
\[
\mu_3(V_3) = \mu_3(M[6^-,1]) = M[6,1^+] = M[6,3] = {\bsm&3\\2\esm} .
\]
This gives rise to the two short exact sequences (which can be read from the
graph $\mu_1(\G_\ii)$):
\[
0\to V_3 \to M_3\oplus V_6 \to M[6,3] \to 0,\qquad 0\to M[6,3] \to V_5 \to V_3\to 0. 
\]
By Proposition~\ref{mutY}, we thus have
\[
 Y_{M[6,3]}Y_{V_3} = q^{[M[6,3],V_3]}\left(q^{-1} Y_{M_3\oplus V_6} + Y_{V_5}\right) 
= q^{-1} Y_{M_3}Y_{V_6} + Y_{V_5}, 
\]
since $[M[6,3],V_3] = 0$, and $[M_3,V_6] = 0$.
Using the notation $M[l,k]$, this can be rewritten
\begin{equation}\label{eqY2}
Y_{M[6,3]}Y_{M[3,1]} = q^{-1}Y_{M[3,3]}Y_{M[6,1]} + Y_{M[5,2]}. 
\end{equation}

\medskip\noindent
\emph{Mutation 3:} One mutates at $V_2$ in the maximal rigid module $\mu_3\mu_1(V_\ii)$.
One has 
\[
\mu_2(V_2) = \mu_2(M[5^-,2]) = M[5,2^+] = M[5,5] = M_5.
\]
This gives rise to the two short exact sequences (which can be read from the
graph $\mu_3\mu_1(\G_\ii)$):
\[
0\to V_2 \to V_5 \to M_5 \to 0,\qquad 0\to M_5 \to M_3\oplus V_4 \to V_2\to 0. 
\]
By Proposition~\ref{mutY}, we thus have
\[
 Y_{M_5}Y_{V_2} = q^{[M_5,V_2]}\left(q^{-1} Y_{V_5} + Y_{M_3\oplus V_4}\right) 
= Y_{V_5} + q Y_{M_3}Y_{V_4}, 
\]
since $[M_5,V_2] = 1$, and $[M_3,V_4] = 0$.
Using the notation $M[l,k]$, this can be rewritten
\begin{equation}\label{eqY3}
Y_{M[5,5]}Y_{M[2,2]} = Y_{M[5,2]} + qY_{M[3,3]}Y_{[4,4]}. 
\end{equation}

\medskip\noindent
\emph{Mutation 4:} One mutates at $M_3$ in the maximal rigid module $\mu_2\mu_3\mu_1(V_\ii)$.
One has 
\[
\mu_1(M_3) = \mu_1(M[6^-,3]) = M[6,3^+] = M[6,6] = M_6.
\]
This gives rise to the two short exact sequences (which can be read from the
graph $\mu_2\mu_3\mu_1(\G_\ii)$):
\[
0\to M_3 \to M[6,3] \to M_6 \to 0,\qquad 0\to M_6 \to M_5 \to M_3\to 0. 
\]
By Proposition~\ref{mutY}, we thus have
\[
 Y_{M_6}Y_{M_3} = q^{[M_6,M_3]}\left(q^{-1} Y_{M[6,3]} + Y_{M_5}\right) 
= q^{-1} Y_{M[6,3]} + Y_{M_5}, 
\]
since $[M_6,M_3] = 0$.
Using the notation $M[l,k]$, this can be rewritten
\begin{equation}\label{eqY4}
Y_{M[6,6]}Y_{M[3,3]} = q^{-1}Y_{M[6,3]} + Y_{M[5,5]}. 
\end{equation}
Now, in $A_q(\n)$ we have the following quantum $T$-system given by Proposition~\ref{Tsystem}:
\begin{align}
D(1,3)D(0,1) &= q^{-1}D(0,3) + D(0,2),\label{eqT1}\\[2mm]
D(1,6)D(0,3) &= q^{-1}D(1,3)D(0,6) + D(0,5),\label{eqT2}\\[2mm] 
D(2,5)D(0,2) &= D(0,5) + qD(1,3)D(0,4),\label{eqT3}\\[2mm] 
D(3,6)D(1,3) &= q^{-1}D(1,6) + D(2,5).\label{eqT4}
\end{align}
By definition, the homomorphism $\k \colon \AA_{\Q(q)}(\md\L) \to F_q(\n)$ satisfies
\begin{align*}
&\k(Y_{M[1,1]}) = D(0,1),\quad
\k(Y_{M[3,1]}) = D(0,3),\quad
\k(Y_{M[6,1]}) = D(0,6),
\\[2mm]
&\k(Y_{M[2,2]}) = D(0,2),\quad
\k(Y_{M[5,2]}) = D(0,5),\quad
\k(Y_{M[4,4]}) = D(0,4). 
\end{align*}
Thus, comparing (\ref{eqY1}) and (\ref{eqT1}) we see that $\k(Y_{M[3,3]})=D(1,3)$.
Next, comparing (\ref{eqY2}) and (\ref{eqT2}) we see that $\k(Y_{M[6,3]})=D(1,6)$.
Next, comparing (\ref{eqY3}) and (\ref{eqT3}) we see that $\k(Y_{M[5,5]})=D(2,5)$.
Finally, comparing (\ref{eqY4}) and (\ref{eqT4}) we see that $\k(Y_{M[6,6]})=D(3,6)$.
So in particular, we have
\[
 \k(Y_{M_k}) = D(k^-,k) = E^*(\b_k),\qquad (1\le k \le 6).
\]
This shows that $\k\left(\AA_{\Q(q)}(\md\L)\right)$ contains $A_q(\n)$. 
Comparing dimensions of $Q_+$-homogeneous components, as in Proposition~\ref{rank}, 
shows that $\k\left(\AA_{\Q(q)}(\md\L)\right) = A_q(\n) \cong U_q(\n)$.

In this isomorphism, the Chevalley generators $e_1, e_2, e_3$ of $U_q(\n)$
correspond to the quantum cluster variables $Y_{M_1}, Y_{M_3}, Y_{M_6}$,
respectively.
The quantum Serre relation $e_1e_3 = e_3e_1$ corresponds to the fact that
\[
[M_1, M_6] = [M_6,M_1] = 0,\qquad [M_1, M_6]^1 = 0. 
\]
To recover, for instance, the relation $e_1^2e_2 -(q+q^{-1})e_1e_2e_1 + e_2e_1^2 = 0$,
we start from the mutation relations
\[
Y_{M_3}Y_{M_1} = q^{-1} Y_{V_3} + Y_{V_2},
\qquad
Y_{M_1}Y_{M_3} = q^{-1} Y_{V_2} + Y_{V_3}.  
\]
The first one is our first mutation above, and the second one is the mutation back,
in the opposite direction (mutations are involutive).
Eliminating $Y_{V_2}$ between these two equations, we get
\[
Y_{M_1}Y_{M_3} - q^{-1}Y_{M_3}Y_{M_1} = (1-q^{-2})Y_{V_3}. 
\]
Since $[M_1,V_3]^1 = 0$, $[M_1,V_3] = 1$, and $[V_3,M_1]= 0$, we have
\[
Y_{M_1}Y_{V_3} = qY_{V_3}Y_{M_1},  
\]
which yields
\[
Y_{M_1}^2 Y_{M_3} - q^{-1}Y_{M_1}Y_{M_3}Y_{M_1}
= qY_{M_1}Y_{M_3}Y_{M_1} - Y_{M_3}Y_{M_1}^2,
\]
as expected.

\subsection{Canonical basis of $\AA_{\Q(q)}(\CC_w)$}\label{canonbasACw}

By \cite[\S4.7]{Ki}, the subalgebra $U_q(\n(w))$ is spanned
by a subset $\B^*(w)$ of the dual canonical basis $\B^*$ of $U_q(\n)$.
Using the isomorphism $\k\colon \AA_{\Q(q)}(\CC_w) \stackrel{\sim}{\to} A_q(\n(w)) \cong U_q(\n(w))$,
we can pull back $\B^*(w)$ and obtain a $\Q(q)$-basis of $\AA_{\Q(q)}(\CC_w)$,
which we shall call the \emph{canonical basis of $\AA_{\Q(q)}(\CC_w)$}, and denote
by $\BB(w) = \{b(\aa) \mid \aa\in\N^r\}$.
It may be characterized as follows.
Recall the involution $\si$ of $\AA_\A(\CC_w)$ defined in \S\ref{sectinvolsigma}.
We will also denote by $\si$ its extension to $\AA_{\Q(q)}(\CC_w)$.

\begin{Prop}\label{caracter}
For $\aa=(a_1,\ldots,a_r)\in\N^r$, the vector $b(\aa)$ is uniquely determined by the following
conditions:
\begin{itemize}
\item[(a)] the expansion of $b(\aa)$ on the basis $\{Y_M(\cc)\mid \cc\in \N^r\}$ is of the
form
\[
b(\aa) = Y_M(\aa) + \sum_{\cc\not = \aa} \ga_{\aa,\cc}(q)\,  Y_M(\cc)
\]
where $\ga_{\aa,\cc}(q) \in q^{-1}\Q[q^{-1}]$ for every $\cc\not = \aa$;  
\item[(b)] Let $\b(\aa) := \sum_{1\le k\le r}a_k \b_k$. Then $\si(b(\aa)) = q^{N(\b(\aa))}b(\aa)$.
\end{itemize}
\end{Prop}

\proof
This is a restatement in our setup of \cite[Theorem 4.26]{Ki}. 
The same result was previously obtained in \cite{La} in a particular case.
\cqfd  

It follows from (\ref{sigmaY}) that all quantum cluster monomials satisfy 
condition (b) of Proposition~\ref{caracter}. 
This is similar to \cite[Proposition 10.9 (2)]{BZ}.
Unfortunately, it is not easy 
to prove that quantum cluster monomials satisfy (a), so one can only conjecture

\begin{Conj}\label{conj}
All quantum cluster monomials $Y_R$, where $R$ runs over the set
of reachable rigid modules in $\CC_w$, belong to $\BB(w)$.
\end{Conj}

It follows from the original work of Berenstein and Zelevinsky \cite{BZ0}
that the conjecture holds in the prototypical Example~\ref{exampleA3},
namely, for $\g= \Sl_4$ the dual canonical basis of $U_q(\n)$ is equal to the set of  
quantum cluster monomials. In this case, there are 14 clusters and 
12 cluster variables (including the frozen ones), which all are 
unipotent quantum minors.

We note that the conjecture is satisfied when $R=M[b,a]$ is one
of the modules occuring in the quantum determinantal identities.
Indeed, $Y_{M[b,a]}$ is then a quantum minor, and belongs to $\BB(w)$
by Proposition~\ref{qmincan}.

It is proved in \cite{La,La2} that the conjecture holds for all quantum
cluster variables when $\AA_{\Q(q)}(\CC_w)$ is associated with the
Kronecker quiver or a Dynkin quiver of type $A$, and $w$ is the square of a Coxeter element.

Conjecture~\ref{conj} would imply the open orbit conjecture
of \cite[\S18.3]{GLS} for reachable rigid modules, by specializing~$q$ to $1$.
It also appears as Conjecture~1.1 (2) in \cite{Ki}.

\subsection{Quantum coordinate rings of matrices}

Let $\g$ be of type $A_{n}$. 
Let $j$ be a fixed integer between $1$ and $n$, and set $k=n+1-j$.
Let $w$ be the Weyl group element with
reduced decomposition
\[
 w= (s_j s_{j-1} \cdots s_1)(s_{j+1} s_j \cdots s_2)\cdots (s_n s_{n-1} \cdots s_{k}).
\]
We denote by $\ii=(i_{jk},i_{jk-1},\ldots, i_1)$ the corresponding word.
It is well known (see \cite{MC}) that for this particular choice 
of $\g$ and $w$, $U_q(\n(w))$ is isomorphic to 
the quantum coordinate ring $A_q(\Mat(k,j))$ of the space of $k\times j$-matrices.

In this case, $\CC_w$ is the subcategory of $\mod(\L)$ generated by the indecomposable
projective $P_k$ with simple top $S_k$. 
For $1\le a \le k$ and $1\le b \le j$,
the module $P_k$ has a unique quotient $X_{ab}$
with dimension vector 
$\dimv\, X_{ab} = \a_a+\a_{a+1}+\cdots + \a_{n+1-b}$.
It is not difficult to check that
\[
 M_\ii = \bigoplus\limits_{\substack{1\le a\le k\\[1mm] 1\le b\le j}} X_{ab}.
\]
Then, setting $x_{ab} := Y_{X_{ab}} \in \AA_q(\CC_w)$, we have that the
elements $x_{ab}$ satisfy the defining relations of $A_q(\Mat(k,j))$,
namely
\begin{align}
x_{ab}x_{ac} & =  qx_{ac}x_{ab},&(b <c),\\
x_{ac}x_{bc} & =  qx_{bc}x_{ac},&(a <b),\\
x_{ab}x_{cd} & =  x_{cd}x_{ab}, &(a<c,\ b>d),\\
x_{ab}x_{cd} & =  x_{cd}x_{ab} + (q-q^{-1})x_{ad}x_{cb}, &(a<c,\ b<d).
\end{align}
Thus, Theorem~\ref{main} gives immediately
\begin{Cor}
The quantum coordinate ring $A_q(\Mat({k,j}))$ is isomorphic to the 
quantum cluster algebra $\AA_{\Q(q)}(\CC_w)$. 
\end{Cor}

\begin{example}\label{ex11_8}
{\rm
The quantum coordinate ring $A_q(\Mat({3,3}))$ is isomorphic to 
$U_q(\n(w))$ for $\g$ of type $A_5$ and $w= s_3s_2s_1s_4s_3s_2s_5s_4s_3$.
The category $\CC_w$ has finite representation type $D_4$.
Taking $\ii = (3,2,1,4,3,2,5,4,3)$, the direct indecomposable summands
of $V_\ii$ are
\begin{align*}
V_1 &= {\bsm3\esm} &
V_2 &= {\bsm3\\&4\esm} &
V_3 &= {\bsm3\\&4\\&&5\esm} \\[2mm]
V_4 &= {\bsm&3\\2\esm} &
V_5 &= {\bsm&3\\2&&4\\&3\esm} &
V_6 &= {\bsm&3\\2&&4\\&3&&5\\&&4\esm} \\[2mm]
V_7 &= {\bsm&&3\\&2\\1\esm} &
V_8 &= {\bsm&&3\\&2&&4\\1&&3\\&2\esm} &
V_9 &= {\bsm&&3\\&2&&4\\1&&3&&5\\&2&&4\\&&3\esm}
\end{align*}
The direct indecomposable summands
of $M_\ii$ are
\begin{align*}
M_1 &= {\bsm3\esm} &
M_2 &= {\bsm3\\&4\esm} &
M_3 &= {\bsm3\\&4\\&&5\esm} \\[2mm]
M_4 &= {\bsm&3\\2\esm} &
M_5 &= {\bsm&3\\2&&4\esm} &
M_6 &= {\bsm&3\\2&&4\\&&&5\esm} \\[2mm]
M_7 &= {\bsm&&3\\&2\\1\esm} &
M_8 &= {\bsm&&3\\&2&&4\\1&&\esm} &
M_9 &= {\bsm&&3\\&2&&4\\1&&&&5\esm}
\end{align*}
In the identification of $A_q(\Mat({3,3}))$ with $\AA_{\Q(q)}(\CC_w)$, we have
\begin{align*}
x_{11}& = Y_{M_9},& x_{12} & = Y_{M_8}, &x_{13} &= Y_{M_7}, \\[2mm]
x_{21}& = Y_{M_6},& x_{22} & = Y_{M_5}, &x_{23} &= Y_{M_4}, \\[2mm]
x_{31}& = Y_{M_3},& x_{32} & = Y_{M_2}, &x_{33} &= Y_{M_1}.
\end{align*}
The quantum cluster algebra $\AA_{\Q(q)}(\CC_w)$ has finite cluster type $D_4$.
Its ring of coefficients is the skew polynomial ring in the variables
$Y_{V_3}, Y_{V_6}, Y_{V_7}, Y_{V_8}, Y_{V_9}$.
It has 16 non frozen cluster variables, namely the 9 canonical generators $Y_{M_i}$,
the 4 unipotent quantum minors $Y_{T_i}$ attached to the modules:
\begin{align*}
T_1 &= {\bsm&&3\\&2&&4\\1&&3&&5\esm}&
T_2 &= {\bsm&&3\\&2&&4\\1&&3&&5\\&2\esm} &
T_3 &= {\bsm&&3\\&2&&4\\1&&3&&5\\&&&4\esm}&
T_4 &= {\bsm&&3\\&2&&4\\1&&3&&5\\&2&&4\esm} 
\end{align*}
and the elements $Y_{U_i}$ attached to the 2 modules
\begin{align*}
U_1 &= {\bsm&&33\\&22&&44\\1&&3&&5\esm}&
U_2 &= {\bsm&&3\\&2&&4\\1&&33&&5\\&2&&4\esm} 
\end{align*}
with dimension vector $\a_1+2\a_2+3\a_3+2\a_4+\a_5$.
These cluster variables form 50 clusters.
} 
\end{example}

\subsection{Quantum coordinate rings of open cells in partial
flag varieties}

In this section, we assume that $\g$ is a simple Lie algebra of simply-laced
type.
We briefly review some classical material, using the notation of \cite{GLSflag,GLS}.

Let $G$ be a simple simply connected complex algebraic group 
with Lie algebra $\g$.
Let $H$ be a maximal  torus of $G$, and $B, B^-$ a pair of opposite 
Borel subgroups containing $H$ with unipotent radicals $N, N^-$. 
We denote by $x_i(t)\ (i\in I,\, t\in\C)$ the simple root subgroups
of $N$, and by $y_i(t)$ the corresponding simple root subgroups of
$N^-$.

We fix a non-empty subset $J$ of $I$ and
we denote its complement by $K=I\setminus J$.
Let $B_K$ be the standard parabolic 
subgroup of $G$ generated by $B$
and the one-parameter subgroups 
\[
y_k(t),\qquad (k\in K,\,t\in\C).
\]
We denote by $N_K$ the unipotent radical of $B_K$.
Similarly, let $B_K^-$ be the parabolic subgroup of $G$ generated by $B^-$
and the one-parameter subgroups 
\[
x_k(t),\qquad (k\in K,\,t\in\C).
\]
The projective variety $B_K^-\backslash G$ is called
a \emph{partial flag variety}.
The natural projection map
$$
\pi\colon G \ra  B_K^-\backslash G
$$ 
restricts to an embedding of $N_K$ into 
$B_K^-\backslash G$ as a dense open subset.

Let $W_K$ be the subgroup of the Weyl group $W$ generated by the 
reflections $s_k\ (k\in K)$. This is a finite Coxeter group and
we denote by $w_0^K$ its longest element. The longest element
of $W$ is denoted by $w_0$.
It is easy to check that 
\begin{equation}
N_K = N\left(w_0w_0^K\right), 
\end{equation}
see \cite[Lemma 17.1]{GLS}.
It follows that $\C[N\left(w_0w_0^K\right)]$ can be identified with the 
coordinate ring of the affine open chart $\O_K:=\pi(N_K)$ of $B_K^-\backslash G$. 
Therefore, $A_q(\n(w_0w_0^K))$ can be regarded as the quantum coordinate
ring $A_q(\O_K)$ of $\O_K$, and Theorem~\ref{main} implies:
\begin{Cor}
The quantum coordinate ring $A_q(\O_K)$ is isomorphic to the 
quantum cluster algebra $\AA_{\Q(q)}(\CC_{w_0w_0^K})$. 
\end{Cor}

\begin{example}
{\rm Take $G=SL(6)$, and $J=\{3\}$.
Then $B_K^-\backslash G$ is the Grassmannian $\Gr(3,6)$ of $3$-dimensional subspaces
of $\C^6$, and $N_K = N(w_0w_0^K)$, where $w=w_0w_0^K$ is as in Example~\ref{ex11_8}.
Here, $N(w)$ can be identified with the open cell $\O$ of $\Gr(3,6)$ given by
the non-vanishing of the first Pl\"ucker coordinate.
The quantum cluster algebra $\AA_{\Q(q)}(\CC_w)$ of 
Example~\ref{ex11_8} can be regarded as the quantum coordinate ring of $\O$.
} 
\end{example}

\begin{example}
{\rm Take $G=SO(8)$, of type $D_4$. We label the vertices of the Dynkin diagram 
from $1$ to $4$ so that the central node is $3$. Let $J=\{4\}$.
Then $B_K^-\backslash G$ is a smooth projective quadric ${\mathcal Q}$ of dimension 6,
and $N_K$ can be regarded as an open cell $\O$ in ${\mathcal Q}$. 
Here $N_K=N(w)$ where $w = w_0w_0^K = s_4s_3s_1s_2s_3s_4$.
It is easy to check that the elements $Y_{M_k}\ (1\le k \le 6)$
satisfy
\begin{align}
Y_{M_j}Y_{M_i} & =  q Y_{M_i}Y_{M_j}, \qquad (i<j,\ i+j\not = 7),\\
Y_{M_4}Y_{M_3} & =  Y_{M_3}Y_{M_4},\\
Y_{M_5}Y_{M_2} & =  Y_{M_2}Y_{M_5} + (q-q^{-1})Y_{M_3}Y_{M_4},\\
Y_{M_6}Y_{M_1} & =  Y_{M_1}Y_{M_6} + (q-q^{-1})(Y_{M_2}Y_{M_5} - q^{-1} Y_{M_3}Y_{M_4}),
\end{align}
and that this is a presentation of the quantum coordinate ring $A_q(\O)$.
This shows that $A_q(\O)$ is isomorphic to the quantum coordinate ring of the
space of $4\times 4$-skew-symmetric matrices, introduced by Strickland \cite{St}.
The category $\CC_w$ has finite representation type $A_1\times A_1$.
Hence, there are 4 cluster variables $Y_{M_1}$, $Y_{M_2}$, $Y_{M_5}$, $Y_{M_6}$, together
with 4 frozen ones, namely, $Y_{V_3}=Y_{M_3}$, $Y_{V_4}=Y_{M_4}$, and 
\[
Y_{V_5}  =  Y_{M_2}Y_{M_5} - q^{-1}Y_{M_3}Y_{M_4},\qquad
Y_{V_6}  =  Y_{M_1}Y_{M_6} -q^{-1}Y_{M_2}Y_{M_5} + q^{-2} Y_{M_3}Y_{M_4}.
\]
Observe that $Y_{V_6}$ coincides with the quantum Pfaffian of \cite{St}.
There are 4 clusters
\[
\{Y_{M_1},\ Y_{M_2}\},\quad
\{Y_{M_1},\ Y_{M_5}\},\quad
\{Y_{M_6},\ Y_{M_2}\},\quad
\{Y_{M_6},\ Y_{M_5}\}. 
\]
Note that since all quantum cluster variables belong to the basis
$\{Y_M(\cc)\mid \cc\in \N^r\}$, Conjecture~\ref{conj} is easily verified 
in this case.
} 
\end{example}

\subsection*{Acknowledgments}
We thank the Hausdorff Center for Mathematics in Bonn for organizing
a special trimester on representation theory in the spring of 2011,  
during which this paper was finalized.
The first author acknowledges support from 
the grants PAPIIT IN117010-2 and CONACYT 81948,
and thanks the Max Planck Institute in Bonn for hospitality.
The second author is grateful to the University of Bonn for
several invitations.
The third author also thanks the
Transregio SFB/TR 45 for support.


\bigskip
\small

\noindent
\begin{tabular}{ll}
Christof {\sc Gei{ss}} : &
Instituto de Matem\'aticas,\\ 
&Universidad Nacional Aut\'onoma de M\'exico\\
& 04510 M\'exico D.F., M\'exico \\
&email : {\tt christof@math.unam.mx}\\[5mm]
Bernard {\sc Leclerc} :&
Universit\'e de Caen, LMNO UMR 6139 CNRS,\\
& 14032 Caen cedex, France\\
& Institut Universitaire de France\\
&email : {\tt bernard.leclerc@unicaen.fr}\\[5mm]
Jan {\sc Schr\"oer} :&
Mathematisches Institut, Universit\"at Bonn,\\
&Endenicher Allee 60, D-53115 Bonn, Germany\\
&email : {\tt schroer@math.uni-bonn.de}
\end{tabular}

\end{document}